\documentclass[3p]{elsarticle}

\newtheorem{Theorem}{Theorem}

\newtheorem{Example}[Theorem]{Example}
\newtheorem{Definition}[Theorem]{Definition}
\newtheorem{Remark}[Theorem]{Remark}

\newtheorem{Lemma}[Theorem]{Lemma}

\newenvironment{Proof}[1][Proof]{\noindent\textbf{#1: }}

\usepackage{amsmath, amssymb,  amsfonts, txfonts, pxfonts,extarrows,tensor}
\usepackage{comment}
  
\usepackage{graphics} 
\usepackage{amssymb}  
\usepackage{amsmath}
\usepackage{amscd}
\usepackage{euscript}
\usepackage{enumerate}
\usepackage[all]{xy}
\usepackage{rlepsf}

\CompileMatrices 
\usepackage{paralist,inputenc}
\usepackage{hyperref}
\usepackage{mathtools}

\usepackage{color}

\newcommand{\h}{\mathfrak{h}}
\newcommand{\dd}{\partial}

\newcommand{\stk}[1]{\stackrel{#1}{\longrightarrow}}
\newcommand{\ov}[1]{\overline{#1}}

\def \v {\overrightarrow}
\def \cC {{\cal C}}

\def \s {\sigma}

\def \GL {\mathsf{gl}}

\def \U {\EuScript{U}}
\def \s {\sigma}

\def \C {\mathbb{C}}
\def \Cc {\EuScript{C}}

\def \w {\omega}
\def \id {\mathrm{id}}
\def \g {\mathfrak {g}}

\def \R {\mathbb{R}}

\def \dd {\partial}

\def \N {\mathbb N}

\def \g {\mathfrak{g}}

\def \tn {\otimes}
\def \t {\triangleright}

\def \ra {\xrightarrow}

\def \lra {\longrightarrow} 
\def \le {\mathfrak{e}}

\def \diff {\mathrm{diff}}

\def \an {\textrm{ and }}
\def \gl {\mathsf{gl}}
\def \G {\EuScript{G}}

\def \an {\textrm{ and }}

\def\be{\begin{equation}}
\def\ee{\end{equation}}
\def\bea{\begin{eqnarray}}
\def\eea{\end{eqnarray}}
\def \we {\textrm{ where }}

\def \h {\mathfrak{h}}

\def \ch {\mathsf{ch}}

\def \ss {\tilde{s}}
\def \tt {\tilde{t}}
\def \nss {\hspace*{-.6cm}}
\def \nsl {\hspace*{-1.1cm}}

\def \Ch {{\bf \rm Chain}}
\def \g {\mathfrak{g} }
\def \lg {\mathfrak{g} }
\def \lh {\mathfrak{h}}
\def \lG {\mathfrak{G}}

\def \V {\EuScript{V}}
\def \U {\EuScript{U}}
\def \cW {\EuScript{W}}
\def \gl {\mathfrak{gl}}

\def \ot {\otimes}

\def \bX {\overline{X}}
\def \bY {\overline{Y}}

\def \br {\overline{r}}
\def \bs {\overline{s}}
\def \bt {\overline{t}}

\def \Un {\mathcal{U}^{(n)}}
\def \Ud {\mathcal{U}^{(2)}}

\def \sl {\mathfrak{sl}_2}
\def \Fo {{\mathbb{F}_0}}
\def \Fu {{\mathbb{F}_1}}
\def \fsl {\Fu\rtimes_{\alpha}\sl}

\def \d {\partial}

\def \la {\mathfrak{a}}

\def \F {\mathbb{F}}
\def \R {\mathcal{R}}
\def \st {\mathfrak{String}}
\def \ufo {1_{\Fo}}

\def \CG {\EuScript{C}_{\lG}}

\newcommand{\defT}[6]{ 
\xymatrix @R=1em{ 
    										      &  {#2} \ar@/^1pc/[dr]^{{#4}}  & \\
{#1} \ar@/^1pc/[ur]^{{#3}} \ar@/_1pc/[dr]_{{#5}}&  \quad \Uparrow {#6} \quad               & {#2} \\
                                                  &  {#1} \ar@/_1pc/[ur]_{{#3}}  &                                                
}
}

\makeatother
\newcommand{\morl}[5]{{   {\xymatrix{& \\ & {#1} \ar@/^2pc/[rr]^{{#4} }\ar@/_2pc/[rr]_{{#3}} & \Uparrow {#5} & {#2}\\ &}}}}

\newcommand{\mor}[5]{{   {\xymatrix{ & {#1} \ar@/^2pc/[rr]^{{#4} }\ar@/_2pc/[rr]_{{#3}} & \Uparrow {#5} & {#2}}}}}

\newcommand{\morr}[7] {{ \xymatrix{ & {#1}\ar@/^3pc/[rr]^{{#5}}\ar[rr]|{{#4}} \ar@/_3pc/[rr]_{{#3}} \ar @/^/ @{{}{ }{}} [rr]^{\Uparrow {#7}} \ar @/_/ @{{}{ }{}} [rr]_{\Uparrow{#6}}
 && {#2}}}}
\newcommand{\comp}[2]{ {\begin{array}{c}{#2}\\ {#1} \end{array}}  }  

\newcommand{\rw}[7]{{   {\xymatrix{ & {#1} \ar@/^2pc/[rr]^{{#4} }\ar@/_2pc/[rr]_{{#3}} & \Uparrow {#5} & {#2} \ar[r]^{{#7}} &{#6} }}}}

\newcommand{\lw}[7]{{   {\xymatrix{ &{#6} \ar[r]^{{#7}}   & {#1} \ar@/^2pc/[rr]^{{#4} }\ar@/_2pc/[rr]_{{#3}} & \Uparrow {#5} & {#2}}}}}

\newcommand{\hcomp}[9]{{   {\xymatrix{ & {#1} \ar@/^2pc/[rr]^{{#4} }\ar@/_2pc/[rr]_{{#3}} & \Uparrow {#5} & {#2}  \ar@/^2pc/[rr]^{{#8} }\ar@/_2pc/[rr]_{{#7}} & \Uparrow {#9} & {#6} }}}}


\newcommand{\nat}[4] {\xymatrix{&& #2 \tn #3  \ar@/^1pc/[rd]^{r_{#2,#3}} & \\ &#1 \tn #3 \ar@/_1pc/[rd]_{r_{#1,#3}}\ar@/^1pc/[ru]^{#4 \tn #1 } &\quad \quad \quad \Uparrow T_{(#4,#3)} &#2 \tn #3\\&& #1 \tn #3 \ar@/_1pc/[ru]_{#4 \tn #3}&}}


\newcommand{\tot}[4] {\xymatrix{&& #2 \tn #3  \ar@/^1pc/[rd]^{r_{#2,#3}} & \\ &#1 \tn #3 \ar@/_1pc/[rd]_{r_{#1,#3}}\ar@/^1pc/[ru]^{#4 \tn #1 } &\quad \quad \quad \Uparrow T_{(#4,#3)} &#2 \tn #3\\&& #1 \tn #3 \ar@/_1pc/[ru]_{#4 \tn #3}&}}
\makeatletter

\newcommand{\natr}[4] {\xymatrix{&& #3 \tn y \ar@/^1pc/[rd]^{r_{#3,#1}} & \\ &#2 \tn #1 \ar@/_1pc/[rd]_{r_{x,y}}\ar@/^1pc/[ru]^{#4 \tn #1 } &\quad \quad \quad \Uparrow T_{(#4,#1)} &#3 \tn #1\\&& #2 \tn #1 \ar@/_1pc/[ru]_{#4 \tn #1}&}}

\newcommand{\natcond}[6]{
\xymatrix{
& & & #2 \tn #3 \ar@/^1pc/[rrd]^{r_{#2,#3}} & & \\ 
&#1 \tn #3 \ar@/_1pc/[rrd]_{r_{#1#3}} \ar@/^1pc/[rru]^{#5 \tn #3 } \ar@/_1pc/[rru]_{#4 \tn #3} \ar@{{}{ }{}} [rru]|{ \Uparrow #6 \tn #3}  & & \quad \quad \quad \Uparrow T_{(#4,#3)} & &#2 \tn #3 \\
&&&#1 \tn #3 \ar@/_1pc/[rru]_{#4 \tn #3} & & 
} = \\    
\xymatrix{&&&#2 \tn #3 \ar@/^1pc/[rrd]^{r_{#2,#3}} && \\ &#1 \tn #3 
 \ar@/_1pc/[rrd]_{r_{#1#3}}\ar@/^1pc/[rru]^{#5 \tn #3 }  && \Uparrow T^1_{(#5,#3)} \quad \quad \quad \quad &&#2 \tn #3\\&&&#1 \tn #3 \ar@/^1pc/[rru]^{#5 \tn #3} \ar@{{}{ }{}} [rru]|{ \Uparrow #6 \tn #3}
 \ar@/_1pc/[rru]_{#4 \tn #3}&&
}}

\newcommand{\natcondse}[6]{
\xymatrix{&&#2 \tn #3 \ar@/^1pc/[rd]^{r_{#2,#3}} & \\ &#1 \tn #3 \ar@/_1pc/[rd]_{r_{#1#3}}\ar@/^1pc/[ru]^{#5 \tn #3 } \ar@/_1pc/[ru]_{#4 \tn #3} \ar@{{}{ }{}} [ru]|{ \Uparrow #6 \tn #3} & {\quad \quad \quad \quad \quad \quad\Uparrow T_{(#4,#3)}} &#2 \tn #3\\&&#1 \tn #3 \ar@/_1pc/[ru]_{#4 \tn #3}&}\xymatrix{&\\=\\&} \xymatrix{&&#2 \tn #3 \ar@/^1pc/[rd]^{r_{#2,#3}} & \\ &#1 \tn #3 
 \ar@/_1pc/[rd]_{r_{#1#3}}\ar@/^1pc/[ru]^{#5 \tn #3 }  & \Uparrow T^1_{(#5,#3)} \quad \quad \quad \quad &#2 \tn #3\\&&#1 \tn #3 \ar@/^1pc/[ru]^{#5 \tn #3} \ar@{{}{ }{}} [ru]|{ \Uparrow #6 \tn #3}
 \ar@/_1pc/[ru]_{#4 \tn #3}&}}

\newcommand{\natcondl}[6]{ {\begin{CD} 
\xymatrix{&&&#3\tn #2  \ar@/^1pc/[rrd]^{r_{#3,#2}} && \\ &#3 \tn #1  \ar@/_1pc/[rrd]_{r_{#3,#1}}\ar@/^1pc/[rru]^{#3 \tn #5  } \ar@/_1pc/[rru]_{#3 \tn #4 } \ar@{{}{ }{}} [rru]|{ \Uparrow   #3 \tn #6} &&\quad \quad \quad \Uparrow T_{(#3,#4)}& & #3 \tn #2\\&&&#3 \tn #1 \ar@/_1pc/[rru]_{#3 \tn #4}&&}    \\  \xymatrix{&\\=\\&} \xymatrix{&&&#3 \tn #2 \ar@/^1pc/[rrd]^{r_{#3,#2}} && \\ &#3 \tn #1 
 \ar@/_1pc/[rrd]_{r_{#3,#1}}\ar@/^1pc/[rru]^{#3 \tn #5 }  && \Uparrow T^1_{(#3,#5)} \quad \quad \quad \quad &&#3 \tn #2\\&&&#3 \tn #1 \ar@/^1pc/[rru]^{#3 \tn #5} \ar@{{}{ }{}} [rru]|{ \Uparrow #3 \tn #6}
 \ar@/_1pc/[rru]_{#3 \tn #4}&&}\end{CD} }}
\newcommand{\mixl}[6]{\xymatrix{ &&&&#2\tn #4 \ar@/^1pc/[rdd]^{r_{#2,#4}}     &\\ &&#2 \tn #3 \ar@/^1pc/[rru]^{#2 \tn #6 }\ar@/^1pc/[rd]^{r_{#2,#3}} & &\Uparrow T_{(#2,#6)}\\ &#1 \tn #3 \ar@/_1pc/[rd]_{r_{#1,#3}}\ar@/^1pc/[ru]^{#5 \tn #3 } &\quad \quad \quad \Uparrow T_{(#5,#3)} &#2 \tn #3 \ar@/_1pc/[rr]_{#2 \tn #6}   & &#2 \tn #4   \\&&#1 \tn #3 \ar@/_1pc/[ru]_{#5 \tn #3}&}}
\newcommand{\mixr}[6]{{\xymatrix{ &&&&#2\tn #4 \ar@/^1pc/[rdd]^{r_{#2,#4}}     &\\ &&#1 \tn #4 \ar@/^1pc/[rru]^{#5 \tn #4 }\ar@/^1pc/[rd]^{r_{#1,#4}} & &\Uparrow T_{(#5,#4)}\\ &#1 \tn #3 \ar@/_1pc/[rd]_{r_{#1,#3}}\ar@/^1pc/[ru]^{#1 \tn #6  } &\quad \quad \quad \Uparrow T_{(#1,#6)} &#1 \tn #4 \ar@/_1pc/[rr]_{#5 \tn #4}   & &#2 \tn #4   \\&&#1 \tn #3 \ar@/_1pc/[ru]_{#1 \tn #6}&}}} 
\makeatletter

\def \cB {\EuScript{B}}

\def \Ra {\xRightarrow}

\begin{document}
\author[A]{Lucio Simone Cirio}
\ead{lcirio@uni-math.gwdg.de}

\author[B]{Jo\~{a}o Faria Martins}
\ead{jn.martins@fct.unl.pt}

\address[A]{Mathematisches Institut - Universit\"{a}t G\"{o}ttingen
Bunsenstr. 3-5, 37073 G\"{o}ttingen, Germany
 }

\address[B]{Departamento de Matem\'{a}tica and Centro de Matem\'{a}tica e Aplica\c{c}\~{o}es\\ Faculdade de Ci\^{e}ncias e Tecnologia (Universidade Nova de Lisboa),
Quinta da Torre, 2829-516 Caparica, Portugal}


\title{Infinitesimal 2-braidings and differential crossed modules}

\begin{abstract}
We categorify the notion of an infinitesimal braiding in a linear strict symmetric monoidal category, leading to the notion of a (strict) infinitesimal 2-braiding in a linear symmetric strict monoidal 2-category. We describe the associated categorification of the 4-term relations, leading to six categorified relations. We prove that any infinitesimal 2-braiding gives rise to a flat and fake flat 2-connection in the configuration space of $n$ particles in the complex plane, hence to a categorification of the Knizhnik-Zamolodchikov connection. We discuss infinitesimal 2-braidings in a {certain monoidal} 2-category naturally assigned to every differential crossed module, leading to the notion of a {symmetric} quasi-invariant tensor in a differential crossed module. Finally, we prove that {symmetric} quasi-invariant tensors exist in the differential crossed module associated to {Wagemann's version of} the String Lie-2-algebra. {As a corollary, we obtain a more conceptual proof of the flatness of 
a previously constructed 
categorified Knizhnik-Zamolodchikov connection with values in the String Lie-2-algebra}. 
\end{abstract}
\maketitle
\section*{Introduction}
\noindent A linear (strict)  monoidal category $\Cc$ is a strict monoidal symmetric category $\Cc=(C_0,C_1,\tn,I,B)$, with classes of objects and morphisms ($C_0$ and $C_1$, respectively),  functorial strict tensor product $\tn$ (with strict identity  $I$) and with an involutive braiding $B$. The designation `linear' comes from the fact that we suppose that, given two objects $x,y\in C_0$, the set of morphisms $\hom(x,y)$ is given a vector space structure, and the composition of morphisms is bilinear.

Given a linear strict monoidal category   {$\Cc=(C_0,C_1,\tn,I,B)$}, consider objects $x,y, z \in C_0$. Given a 1-morphism $f\colon x \tn z \to x \tn z$, we define $f^{13}=f^{13}_{x,y,z} \colon x \tn y \tn z \to x \tn y \tn z$ as being (we compose from left to right):
$$ (B_{x,y} \tn \id_z)(\id_y \tn f) (B_{y,x} \tn \id_z)=(\id_x \tn B_{y,z})( f \tn \id_y) (\id_x \tn B_{z,y}). $$ 
Given $g \colon x \tn y \to x \tn y$ we also put $g^{12}=g^{12}_{x,y,z}=g \tn \id_z$. Analogously, given $h\colon y \tn z \to y \tn z$ we define $\id_x \tn h=h^{23}_{x,y,z}=h^{23}.$

An infinitesimal braiding in a linear strict monoidal category $\Cc$  (see \cite[XX.4]{K} and \cite{C}) is given by a family of functorial (natural) isomorphisms $r_{x,y} \colon x \tn y \to x \tn y$, {one} for each pair of objects $x$ and $y$ of $\Cc$, such that:
\begin{enumerate}
 \item For each $x,y \in C_0$ {we have} $\; B_{x,y}r_{y,x}=r_{x,y} B_{x,y} \,$.
\item For each $x,y ,z\in C_0$ we have $\; r_{x,y \tn z}=r^{12}_{x,y} +r^{13}_{x,z} \,$.
\end{enumerate}
The naturality condition yields, in particular, that, given {objects} $x,y,z \in C_0$:
$$r_{y,z}^{23}\,\, r_{x,y\tn z}=r_{x,y\tn z} \,\, r_{y,z}^{23} \, .$$
Combining with the second condition, above, this yields: 
$$r_{y,z}^{23} \,\, r_{x,y}^{12}+r_{y,z}^{23} \,\, r_{x,y}^{13}= r_{x,y}^{12}  \,\,r_{y,z}^{23} + r_{x,y}^{13}\,r_{y,z}^{23} \, .$$
We thus obtain the well known 4-term relations \cite{Bar,Ko,K}, appearing in the theory of universal Vassiliev knot invariants, which we repackage as:
\begin{equation}\label{4t}
 [r_{x,y}^{12}+r_{x,z}^{13},r_{y,z}^{23}]=0.
\end{equation}

{The axioms of an infinitesimally braided monoidal category are {infinitesimal} analogues \cite{C} of the axioms of a braided monoidal category; \cite{JS}.} By using Drinfeld associators \cite{Dr1,Dr}, one can prove that any infinitesimal braiding is the {infinitesimal counterpart (or classical limit)} of a certain {(in general non-strict)} braided monoidal category \cite{C,EK,K}.

{Let $\lg$ be a Lie algebra. Let $\U(\lg)$ be its universal enveloping algebra. Let $\Delta=\Delta^n\colon \U(\lg) \to \U(\lg)^{\otimes n}$, where $\Delta(x)=x\tn 1\tn \ldots \tn 1 +$ cyclic permutations, be the diagonal map.} 
There are two, naturally defined, linear {symmetric}  monoidal categories 
associated to $\lg$. The first is the category $M_\lg$ of $\lg$-modules, with the usual tensor product (with the caveat that it is not strict). {The second is the category $\Cc_\lg$, whose objects are the {non-negative integers} $n \in \N$, and the morphisms are given by pairs $(R,\s)$, where:}
\begin{itemize}
 \item[-] $R \in (\U(\lg)^{\otimes n})_\lg\subset  \U(\lg)^{\otimes n}$ is a $\lg$-invariant tensor, namely $[R,\Delta(X)]=0$ for each $X\in\lg$.
\item[-] $\s\in S_n$  is a permutation of $\{1,\dots,n\}$.
\end{itemize}
(Strictly speaking, the set of 1-morphisms $n \to n$ is the direct sum $\oplus_{s \in S_n}(\U(\lg)^{ \otimes n})_\lg$, and $\hom(m,n)=\{0\}$ if $m \neq n$.)
The composition of 1-morphisms is given by the obvious semi-direct product law, taking into account that the product of $\lg$-invariants tensors in $(\U(\lg)^{\otimes n})_\lg$ is a $\lg$-invariant tensor. {The tensor product is the obvious one, being, on objects $m\tn n=m+n$.} Let $V$ be any representation of $\g$. By sending invariant tensors $R \in  (\U(\lg)^{\otimes n})_\lg$ to the obvious intertwiner $V^{\otimes n} \to V^{\otimes n}$ {(action by $R$)}, and permutations in $S_n$ to the map that permutes factors in the tensor product, we have a {non-strict monoidal} functor $\Cc_\lg \to M_\lg$. {(We will however not use this fact in the following).} Consider a $\lg$-invariant symmetric tensor $r \in \lg \tn \lg$. Then it is easy to see that the family $$\left\{ r_{m,n}=\big((\Delta^m \tn \Delta^n)(r),\id_{S_{m+n}}\big); m,n \in \N\right\},$$
{of maps $r_{m,n} \colon m\tn n \to m \tn n$,} 
is an infinitesimal braiding in the category $\Cc_\g$. Therefore, the 4-term relation $[r^{12}+r^{13},r^{23}]=0$ is satisfied by $r$ {(a well known fact)}.

{Let $n$ be a positive integer.} Recall that  the configuration space of $n$ distinct particles in the complex plane {$\C$} is, by definition:
$$\C(n)=\left\{(z_1,\dots,z_n) \in \C^n\colon z_i \neq z_j \textrm{ if } i \neq j\right\}.$$
This has a left action of the symmetric group $S_n$, by permutation of coordinates. The configuration space of $n$ indistinguishable particles in $\C$ is defined as $\C(n)/S_n$. We recall that the fundamental group of $\C(n)/S_n$ is {isomorphic to the braid group $B_n$, in $n$ strands. Therefore, since the parallel transport of a flat connection is invariant under homotopy,  given a Lie group $G$, the holonomy of a  flat $G$-connection in a vector bundle over  $\C(n)/S_n$ defines a representation of the braid group $B_n$.}

Let $V$ be a vector space. Let $\gl(V)$ be the Lie algebra of linear maps $V \to V$.  Consider a family of endomorphisms $f_{ab} \colon V \to V$, where $1\leq a<b \leq n$, such that $[f_{ab},f_{a'b'}]=0$, if $\{a,b\}\cap \{a',b'\}=\emptyset$.  Define closed 1-forms $\w_{ij}$ in the configuration space $\C(n)$, for $1\leq i,j \leq n$ and $i \neq j$, as:
$$\w_{ij}=\frac{dz_i - dz_j}{z_i-z_j} \, . $$
Then it is well known, and it is proven for example in \cite{Bar,CFM1,CFM2,K,Kh1}, that the $\gl(V)$-valued connection 1-form {in the  configuration space $\C(n)$}:
\begin{equation}\label{defA}
A=\sum_{1 \leq i <j \leq n} \w_{ij}f_{ij}
\end{equation} 
is flat if, and only if, the 4-term relations $[f_{ij}+f_{ik},f_{jk}]=0$ hold {(this is the one of the starting points for constructing the Kontsevich Integral \cite{Bar,Ko})}. In this case {(when the 4-term relations hold)} the family $\{f_{ab} \colon V \to V, \, 1\leq a<b \leq n\}$ {will be} called an infinitesimal Yang-Baxter operator \cite{CFM2}.

What seems very surprising at first sight, and follows from Drinfeld's theory \cite{Dr1,Dr} of quantized universal enveloping algebras of semisimple Lie algebras, is that all quantum group invariants of knots, and therefore also of manifolds, are ultimately determined by the holonomy of the connection form $A$, for  certain naturally defined infinitesimal Yang-Baxter operators.
Given a semsimple Lie algebra $\lg$, let $r=\sum_{i}s_i \tn t_i \in \lg \tn \lg$ be the symmetric invariant tensor in $\lg \tn \lg$ associated to the Cartan-Killing form $\langle,\rangle $ in $\lg$, thus for each $X \in \lg$ we have $\sum_i s_i \langle t_i,X\rangle =X$.  Considering a  $\lg$-module $V$, put $f_{ab}\colon V^{\otimes n} \to  V^{\otimes n}  $ as being the obvious action $r_{ab}$ of  $r$ on $V^{\otimes n}$. This is an infinitesimal Yang-Baxter operator, and for this particular choice of $\{f_{ab}\}$ we call $A$ in \eqref{defA} the Knizhnik-Zamolodchikov (KZ) connection; \cite{KZ,Bar,K}. The KZ-connection is invariant under the natural action of the symmetric group, thus it descends to a flat connection $A$ in the quotient vector bundle, over $\C(n)/S_n$. It is proven in \cite{Dr,Kh1,Kh2} that the representation of the braid-group derived from the quotient KZ connection is equivalent to the representation of the braid group derived from the $\R$-matrix of the quantized universal enveloping 
algebra of $\lg$. This beautiful and deep result is known as the Drinfeld-Kohno Theorem. An even stronger result holds: even though the holonomy of the connection $A$ cannot be immediately extended to links, given that we have particle collisions at local maxima and minima, making some of the the forms $\w_{ab}$ explode at those points, it can easily be regularized, given a choice of a framing, and proven to yield exactly the quantum group invariants of (framed) links; c.f. also the paragraph below.

There is a universal (with respect to all infinitesimal Yang-Baxter operators) KZ-connection in the configuration space $\C(n)$, written in terms of the Lie algebra of horizontal chord diagrams $\ch_n$. Given a positive integer $n$,  $\ch_n$ is the Lie algebra generated by the symbols $r_{ab}$, where $1 \leq a, b \leq n$, subject to the conditions below, normally called infinitesimal braid group relations:
\begin{align*}
& r_{ab} = r_{ba} ,  \\ & [r_{ab},r_{cd}] = 0 \; \mbox{ for } \{a,b\}\cap\{c,d\}=\emptyset \, , \\
& [r_{ab}+r_{ac},r_{bc}]  = 0 = [r_{ab},r_{ac}+r_{bc}]   \, .&  
\end{align*}
{The universal KZ-connection in the configuration space $\C(n)$ has the form below, where $h$ is a formal parameter}: 
$$A_h=\sum_{1\leq a<b \leq n} h\,\w_{ab}\,r_{ab}. $$
This is a flat connection, which is invariant under the obvious action of the symmetric group $S_n$. The holonomy of this connection will define a representation of the braid group in the algebra $\U(\ch(n))[[h]]$ of formal power series in the enveloping algebra $\U(\ch(n))$ of $\ch(n).$ There are well know regularization techniques \cite{LM} for the holonomy of $A_h$ at maximal and minimal points of links (which depend on a choice of a framing). This yields a holonomy with values in the Hopf algebra of chord diagrams in the circle $S^1$, which is, modulo an easily resolvable anomaly, a link invariant, normally called the Kontsevich Integral, being a universal Vassiliev invariant of links; see \cite{Ko,Bar,C,LM}.

{Categorifications of link invariants have been} quite popular since Khovanov's celebrated construction of a categorification of the Jones polynomial in \cite{Kov}, which has been generalized in a plethora of ways (see for example \cite{Bar2,KR}), and proven to be inherently 4-dimensional \cite{Ras}. This also sparked research on the very categorification of the quantum groups themselves \cite{KL}.  Most approaches to the categorification of quantum groups invariants of links use combinatorial frameworks (one notable exception is \cite{Man}). Given the discussion in the previous paragraphs, it seems natural to use tools from Differential Geometry (and here we mean the notion of holonomy, and its higher order versions) as a framework for categorifying quantum knot invariants, and presumably extend them to invariants of braid and  of link cobordisms. This paper is a continuation of previous work in this direction \cite{CFM1,CFM2}, however still on the very heavy algebraic and categorical side of it.  

As indicated by the previous paragraph, in order to use holonomy techniques to define invariants of link and braid cobordisms, we would need a two dimensional notion of holonomy, which make sense for embedded surfaces in a manifold. This can be achieved by considering 2-connections in 2-bundles  \cite{BS1,BS2,BrMe,ACJ}, which should have a structure 2-group (a categorical analogue of a group \cite{BL}), rather than simply a {principal group of structure}. Recall that a strict Lie 2-group can be represented by a Lie crossed module $(\d\colon E \to G,\t)$, where $G$ and $E$ are Lie groups, and $\t$ is a left action of $G$ on $E$ by automorphisms, satisfying natural properties (the Peiffer relations). 

The differential analogue of a strict {Lie} 2-group is a strict Lie 2-algebra \cite{BC}. The latter is uniquely represented by a differential crossed module  $\lG=(\d\colon \lh \to \lg, \t)$, where $\lh$ and $\lg$ are Lie-algebras, $\t$ is a left action of $\lg$ on $\lh$ by derivations and $\d\colon \lh \to \lg$ is a Lie algebra map intertwining the $\g$-action, such that the Peiffer relation $\d(v) \t u=[v,u],$ for each $u,v \in \lh$, is satisfied.
We note that any semisimple Lie algebra $\lg$ gives rise to a family of differential crossed modules (all having the same weak homotopy type), each geometrically realizing the Lie algebra 3-cocycle $\w(X,Y,Z)=\langle X,[Y,Z]\rangle$, where $\langle \, , \rangle$ is the Cartan-Killing form. In the case of $\sl(\C)$, each of these differential crossed modules is called a {String} Lie-2-algebra. For explicit constructions see \cite{BSCS,wag06}. {An explanation of these terms can be found in \ref{sectionstring}.} 

Let $\lG=(\d\colon \lh \to \lg, \t)$ be a differential crossed module, {associated to a Lie-2-group}. A $\lG$-valued 2-connection in a manifold $M$ is locally given by a pair $(A,B)$ of forms, where $A$ is a 1-form with values in $\lg$ and $B$ is a 2-form with values in $\lh$. A 2-connection $(A,B)$ is said to be fake flat if $\d(B)$ agrees with ${\cal F}_A=d A+[A,A]$, the curvature of $A$.  A local {fake flat} 2-connection can be integrated to give a 2-dimensional (surface) holonomy, for smooth maps $\Gamma\colon [0,1]^2 \to M$,  with values in the associated Lie 2-group; see \cite{BS1,BS2,SW1,FMP1,FMP2}.  This surface holonomy is covariant with respect to a gauge 2-group of transformations, and therefore it also makes sense globally for a 2-bundle with a non-local (global) 2-connection;  \cite{SW2,FMP2}. Moreover, if a {fake flat} 2-connection $(A,B)$  is flat, meaning that the 2-curvature 3-form $\G_{(A,B)}=dB+A \wedge B$ of $(A,B)$ vanishes, then this 2-dimensional holonomy is invariant under 
homotopy, relative 
to the boundary of the 
square $[0,1]^2$.

In \cite{CFM1,CFM2} we addressed the categorification of the 4-term relations \eqref{4t} in the language of differential crossed module-valued 2-connections, which led to the definition of an infinitesimal 2-Yang-Baxter operator in a (strict) Lie 2-algebra,  satisfying conditions categorifying the 4-term relation. In this paper we address the categorification of an infinitesimal braiding itself (leading to the definition of a strict infinitesimal 2-braiding, in a linear strict symmetric monoidal 2-category $\CG$) and its applications for defining flat {(and moreover fake flat)} 2-connections in the configuration space of {$n$} particles in the complex plane. We also introduce the definition of a {symmetric} quasi-invariant tensor $(r,c,\xi)$ in a differential crossed module  {$\lG=(\d\colon \lh \to \lg,\t)$}, which categorifies the notion of an invariant tensor in a Lie algebra. Here $r \in \lg \tn \lg$ is a symmetric tensor, {$c$ is in the $\lg$-invariant part of  $\lh$ and $\xi: \lg \to \big(\lg\tn\h\big) \
oplus \big (\lh\tn\lg\big)$ is a linear map, measuring the failure of $r$ 
to be invariant.} Each {symmetric} quasi invariant tensor in $\lG$ leads to the construction of an infinitesimal 2-braiding in a certain  linear monoidal 2-category $\CG$, naturally associated to the differential 
crossed module $\lG$. (We will give explicit details on the construction of $\CG$, which is not at all trivial.)

{Finally, we prove that the {String} Lie-2-algebra (in the version of Wagemann \cite{wag06}) contains a non-trivial {symmetric} quasi-invariant tensor. This will give another proof of the main result of \cite{CFM2}, where we proved that {this particular} {String} Lie-2-algebra contains infinitesimal 2-Yang Baxter operators, providing  a conceptual approach to a proof of the fact that the categorified KZ-connection in \cite{CFM2} is indeed flat}. 

Before describing the plan of the paper, we illustrate in more detail some {possible} applications and related open questions that naturally arise from our construction. We will address these in future publications.

\subsection*{Applications to Lie-2-bialgebras theory and quantum 2-groups}


Our interest in {symmetric} quasi-invariant tensors in Lie 2-algebras {also} stems  from the crucial role we expect them to play in understanding Lie 2-bialgebras as being the classical limit of quantum 2-groups. The analogous statement in the `decategorified` setting involves quasitriangular quasi-Hopf algebras, Drinfeld twists and quasi-Lie bialgebras. We just recall here the main ideas, referring for more details to the original papers by Drinfeld \cite{Dr,Dr1} or to \cite[\S 16]{ChP}. 

A quasitriangular Hopf algebra (QHA) is a Hopf algebra $H$ whose coproduct $\Delta$ is cocommutative up to conjugation by an invertible element $R\in H\tn H$, known as the quasitriangular structure (or ${\cal R}$-matrix). 
Quasitriangular Hopf algebras are exactly those whose module category is strict and braided, 
with triangular Hopf algebras ($R^{21}=R^{-1}$) corresponding to symmetric categories. In a quasi-triangular quasi-Hopf algebra (QqHA) the coproduct $\Delta$ is coassociative up to conjugation by an invertible element $\Phi \in H\tn H\tn H$, known as the Drinfeld associator
(other structures are relaxed as well but they are not relevant for our present discussion),
so that {its} module category {is} still braided but {not} strict, {in an essential way}. The larger class of QqHAs is sometimes more convenient to work with, thanks to a {flexible} notion of equivalence (or `gauge transformations') of QqHAS realized by Drinfeld twists. 
Namely, given a QqHA $H$ and an invertible element $F\in H\tn H$, we can define a new QqHA $H^F$ by `twisting' (suitably conjugating) the structures of $H$ with respect to $F$. 
The two QqHAs $H$ and $H^F$  will {not be isomorphic, in general, nevertheless they} have {braided monoidally equivalent module categories.} 

Among QqHAs, a notable class is formed by deformations of universal enveloping algebras (UEAs). By deformation we mean that all the algebraic structures are defined by formal power series in a deformation parameter $h$, and disregarding terms in $h$ or higher powers (resp. $h^2$ or higher) should correspond to the undeformed structures (resp. taking the classical limit) \cite[\S XVI]{K}.
A quantum universal enveloping algebra (QUEA) is then a QqHA $(A_h,\Delta_h,R_h,\Phi_h)$ such that $A_h/hA_h\simeq \U(\g)$, the universal enveloping algebra of some Lie algebra $\g$ (hence $R_h=1\tn 1$ and $\Phi_h = 1\tn 1\tn 1$, both equalities mod $h$). Drinfeld introduced in \cite{Dr1} quasi-Lie bialgebras as the classical limit of a QUEAs. A quasi-Lie bialgebra $(\g,\delta,\phi)$ is a Lie algebra $\g$ with a compatible Lie cobracket $\delta$ ($\delta$ is a Chevalley-Eilenberg (CE) 1-cocycle on $\g$, $d_{CE}\delta=0$, with values in $\wedge^2\g$)  which satisfies a co-Jacobi identity up to a CE-coboundary $d_{CE}\phi$, with $\phi\in\wedge^3\g$.

Important special cases are Lie bialgebras, corresponding to $\phi=0$, and coboundary Lie bialgebras, where $\delta=d_{CE}\rho:=\delta^{\rho}$ for some $\rho\in\g\tn\g$. In order for $\delta^{\rho}$ to be a well-defined Lie cobracket, the symmetric tensor $\rho_+=\frac{1}{2}(\rho+\rho_{21})$ and $CYB(\rho)=[\rho_{12},\rho_{13}]+[\rho_{12},\rho_{23}]+[\rho_{13},\rho_{23}]$, the classical Yang-Baxter equation on $\rho$, must both be $\g$-invariant. The stronger condition $CYB(\rho)=0$ corresponds to quasitriangular Lie bialgebras (and further adding $\rho_+=0$ we get triangular ones). Now, the classical limit of a QUEA $(A_h,\Delta_h,R_h,\Phi_h)$ with $A_h/hA_h\simeq \U(\g)$ determines a quasi-Lie bialgebra structure on $\g$ with $\delta=\Delta_h - \Delta_h^{21}$ (mod $h^2$) and $\phi=\mathrm{Alt}(\Phi)$ (mod $h^2$). A crucial result proved by Drinfeld is that the element $r=(R^{21}R-1\tn 1)$ (mod $h^2$) is a symmetric invariant tensor {in} $\g$, stable 
under Drinfeld twists. Moreover, by suitably twisting the QUEA, we can make $r = \frac{R-1}{2h}$ (mod $h$) and $\Phi=1\tn 1\tn 1$ (mod $h^2$), so that the classical limit of the QUEA is the quasi-Lie bialgebra with $\delta^{\rho}=0$ (since $\rho=\rho_+=r$) and $\phi=\frac{1}{4}[r^{12},r^{13}]= - \, CYB(r)$. (Note that for $\g$ semisimple $\phi$ is the generator of $H^3(\g,\mathbb{C})$, associated to the Cartan-Killing form, which correspond to the String Lie algebra of $\g$, see Section \ref{sectionstring}.) We denote this quasi-Lie bialgebra as $(\g,r)$, emphasizing that the symmetric invariant tensor $r$ contains all the information of the classical limit of the QUEA. Furthermore, Drinfeld proves that, for $\g$ semisimple,  the  {symmetric $\lg$-invariant tensor} $r$ permits {us} to reconstruct the QUEA, up to twist equivalence, and derives {an} associator precisely by considering the holonomy of the KZ-connection associated to $r$.

A classical limit of Drinfeld twists acts on Lie quasi-bialgebras, and we can twist back $(\g,r)$ to the more familiar form of a quasitriangular Lie bialgebra with $\delta=d_{CE}\rho$, where $\rho=\rho_-+r$ and $\rho_-\in\wedge^2\g$, {a solution} of the modified CYB equation $CYB(\rho_-)= - \, CYB(r)$. We thus see that, given {given an $r$-matrix $\rho = \rho_+ + \rho_-$, the $\lg$-invariant symmetric part
$\rho_+=r$ describes}  the (infinitesimal) braiding associated to the Lie bialgebra, while the antisymmetric part $\rho_-$ generates the Lie cobracket structure (and note that, while $r$ is stable under twists, $\rho_-$ can always be twisted away).

We conjecture that symmetric quasi-invariant tensors {in differential crossed modules}, in addition to describing infinitesimal 2-braidings, can be interpreted in the above setting. Namely, it should be possible to canonically associate a quasi-Lie 2-bialgebra \cite{BSZ13} to the datum of a Lie 2-algebra endowed with a {symmetric} quasi-invariant tensor. We also foresee the existence of Drinfeld twists for quasi-Lie 2-bialgebras, so that the quasi-Lie 2-bialgebra associated to a {symmetric} quasi-invariant tensor can be shown to be twist equivalent to a strict coboundary Lie 2-bialgebra whose {2-$r$-matrix} {(satisfying a categorification of the classical Yang-Baxter equation)} contains the original {symmetric} quasi-invariant tensor as its symmetric part. 

Finally, and most importantly, given a {symmetric} quasi-invariant tensor $(r,c,\xi)$, in a differential crossed module $\lG$, there should exist a quasi-Hopf-2-algebra having $(\lG,r,c,\xi)$ as the classical limit. This `quantization' process will most likely require the categorification of the concept of Drinfeld associators. Examples of Drinfeld 2-associators should be constructible from the {two-dimensional holonomy of the associated categorified KZ-connection.}

\subsection*{Relation with braided monoidal 2-categories}

Our axioms of an infinitesimally braided monoidal 2-category are obtained by directly categorifying Cartier's notion of an infinitesimally braided category \cite{C}, {through} embedding it into the 2-category of 2-vector spaces \cite{BC}, and switching categorical natural transformations into 2-categorical pseudo-natural transformations. Presumably, we could also have followed an alternative approach by directly linearizing the axioms of a braided monoidal 2-category.  {(And in fact, we implicitly used this point of view to derive the  axiom of coherence for infinitesimal 2-braidings, Definition \ref{coherent}, which parallels one analogous axiom appearing in \cite[Definition 6]{BN}: the axiom $S^+=S^-$, imposing the  equality of two hexagons.)} 

{A natural  question to ask} is then how it can be precisely stated that an infinitesimal 2-braiding is the infinitesimal counterpart (or classical limit) of a given braided monoidal 2-category (see also Remark \ref{sd}). {A related issue}  is whether we can always find a braided monoidal 2-category whose classical limit is a given infinitesimal 2-braiding.  Again, this will probably require a categorification of the theory of Drinfeld associators, and for that we expect that the two dimensional holonomy of the {associated categorified Knizhnik-Zamolodchikov connection} may play a major role.

\subsection*{Invariants of braid cobordisms ?}

Braided surfaces will, in this discussion, mean cobordisms $b \xRightarrow {S} b'$ between braids in {$\mathbb{R} \times [-1,1] \times [0,1]$}; \cite{KT}. Braided surfaces form a braided monoidal 2-category $\cB$ \cite{BN,KV}, which follows from the general construction of 2-tangles in \cite{BaLa}. The objects of $\cB$ are non-negative integers $n$, the monoidal structure being given by the sum. The set of 1-morphisms $n \to n'$ is only non-empty when $n=n'$. In the latter case, it is made out of braids $b=(n \ra{b} n) $, with $n$ strands, connecting the set $\{1, \dots, n\}\times \{0,0\}$ {to the set}  $\{1, \dots, n\}\times \{0,1\}$. Braids  are not considered up to isotopy, {but only up to reparametrization}, so {that} they are simply 1-dimensional neat submanifolds of {$\mathbb{R} \times  [-1,1] \times [0,1]$,} with boundary being $\{1,\dots, n\} \times \{0\} \times \{0,1\}$. Braids compose vertically (along $\mathbb{R} \times {[-1,1]} \times \{0, 1\}$) and horizontally in the obvious way, and this 
yields, 
respectively, the composition 
and the tensor product of 1-morphisms. 

The set of 2-morphisms $b \Ra{S} b'$, where $b,b'\colon n \to n$, is given by braided surfaces (cobordisms) connecting $b$ and $b'$. Therefore $S$ is a neat submanifold of $\mathbb{R}\times {[-1,1]} \times [0,1]^2$, with boundary being $\big (b \times \{0\} \big) \cup \big( b' \times\{1\}\big)\cup\big( \{1,\dots, n\} \times \{0\} \times \{ 0,1\} \times [0,1]\big)$. Moreover, the projection $p\colon S \to \{(0,0)\} \times [0,1]^2\cong [0,1]^2$ is to be a branched cover, with simple points only. These (braided surfaces) are considered up to the obvious notion of isotopy, relative to the boundary. Braided surfaces compose in three different ways, giving the {vertical and horizontal compositions, as well as the tensor product of 2-morphisms;} details are  in \cite{BaLa}.

\begin{figure}[h]
\centerline{\relabelbox 
\epsfysize1cm 
\epsfbox{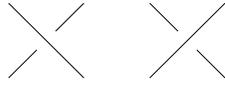}
\endrelabelbox}
\caption{Generating 1-morphisms for the monoidal 2-category of braided surfaces.\label{1mor}}
\end{figure}

\begin{figure}[h]
\centerline{\relabelbox 
\epsfysize6cm 
\epsfbox{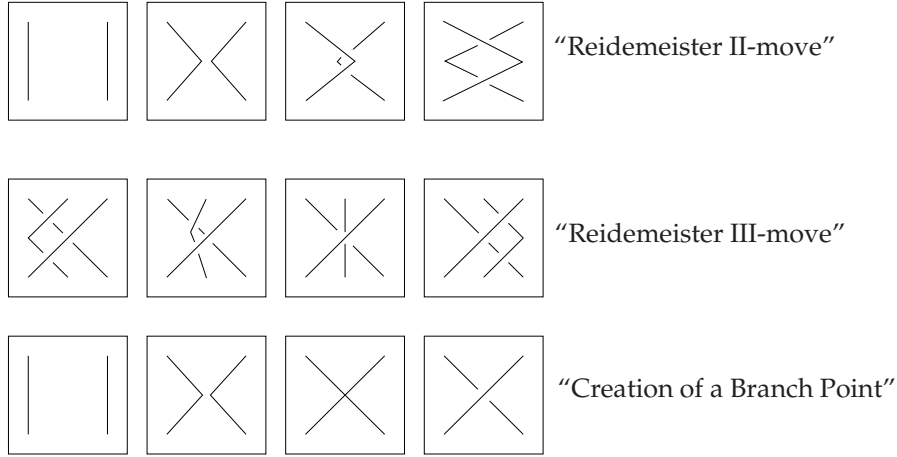}
\relabel{A}{$\textrm{``Reidemeister II-move''} $}
\relabel{B}{$\textrm{``Reidemeister III-move''} $}
\relabel{C}{$\textrm{``Creation of a Branch Point''} $}
 \endrelabelbox}
\caption{{Generating 2-morphisms for the monoidal 2-category of braided surfaces (the first two  called Reidemeister II and III moves). We call the last 2-morphism ``Creation of a Branch Point''.} \label{2mor}}
\end{figure}

The monoidal 2-category $\cB$ can be combinatorially presented, a consequence of the general combinatorial description of knotted surfaces in \cite{CRS,CKS}. Leaving details aside, it is generated by the 1-morphisms of figure \ref{1mor}, and by the 2-morphisms of figure \ref{2mor} (to which mirror images must be added). These satisfy numerous relations (movie moves), compiled in full generality in \cite{CRS,CKS}, and in the case of braided surfaces in  \cite{KV}. In particular 14 relations  are satisfied by these 1- and 2-morphisms, enabling a fully combinatorial description of the monoidal 2-category $\cB$. Of these  {relations}, it plays a primary role the celebrated Zamolodchikov tetrahedron equation, involving only Reidemeister III moves.  

{Only four of the 14 relations between 2-morphisms actually engage the  last generating 2-morphism appearing in figure \ref{2mor} (Creation of a Branch Point).} This morphism is of a very different type of the others. The first two (Reidemeister II and III moves) are generated by isotopies of $S^3$, and the inherent surface is topologically the disjoint union of two disks. In the {latter}, we have a branch point of the projection $p\colon S \to [0,1]^2$, and the surface is topologically the union of two disks, with a 1-handle connecting them (something like a fattened ``H''). 

In \cite{FMP2}, we discussed the issue of assigning a two dimensional holonomy to a surface embedded in a {manifold $M$}, the latter being provided with a 2-connection. Such a holonomy depends on the isotopy class of the choice of embedding, however this ambiguity is {ruled} by the mapping class group of the surface {itself, only.} (For example, in the case of an embedded 2-sphere, the mapping class group being $\{\pm 1\}$ means that a two-dimensional holonomy can be unambiguously assigned to oriented embedded spheres {in $M$}.)

\begin{figure}[h]
\centerline{\relabelbox 
\epsfysize4cm 
\epsfbox{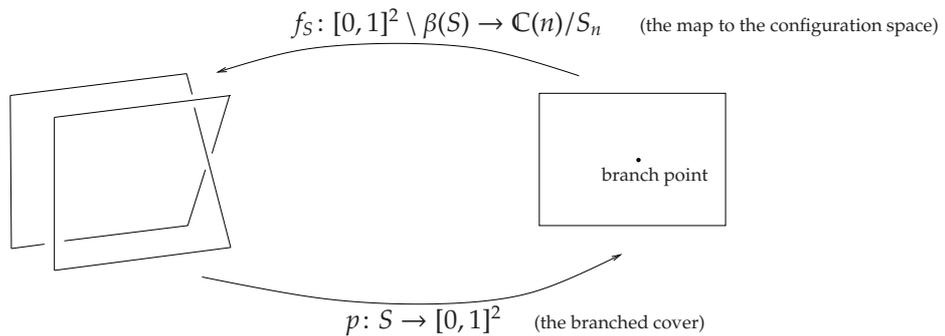}
 \relabel{A}{$\scriptsize{p\colon S \to [0,1]^2 \textrm{\phantom{-----} (the branched cover)}}$}
 \relabel{C}{$\scriptsize{\textrm{branch point}}$}
 \relabel{B}{$\scriptsize{f_S\colon [0,1]^2 \setminus \beta{(S)} \to {\C(n)/S_n} \textrm{  \phantom{-----} (the map to the configuration space) }}$}
\endrelabelbox}
\caption{A braided surface with a branch point.\label{branch}}
\end{figure}

To use the categorified KZ-connection to define invariants of braid cobordisms we must address the following issues. First of all, even though braided surfaces are {\it bona fide} 2-dimensional manifolds (with corners), they do not define embeddings {$[0,1]^2 \to \C(n)/S_n$}, where $\C(n)$ is the configuration space of $n$ particles. Rather, we have a smooth map {$f_S\colon [0,1]^2 \setminus \beta(S) \to \C(n)/S_n$}, where $\beta(S)$ is the set of branch points of the projection $p\colon S \to [0,1]^2$; see figure \ref{branch}. At the branch points, the map $f_S$ has a very special kind of singularities; see \cite{CFM1}. In order to define the two dimensional holonomy of such a singular map $f_S\colon [0,1]^2 \setminus\beta(S)\to {\C(n)/S_n}$, one must find a way to deal with these singularities, as far as their two dimensional holonomy is concerned. This issue is similar to the fact that, for defining invariants of links from the holonomy of the KZ-connection, one must first deal with the 
issue that 
links $L\subset S^3$   do not  define maps {$[0,1] \to \C(n)/S_n$}. Rather, if $L$ is a link, we have a map $f\colon [0,1] \setminus c(L) \to \C(n)$, where $c(L)$ is the set of critical values of the projection of $L$ into the last variable, which we can suppose {to be} a Morse function, with co-domain $[0,1]$. 

In order to define this regularized {two-dimensional} holonomy, we probably will need to address categorified chord diagrams, assigned to general two dimensional manifolds, {rather than merely disjoint union of disks}. (It is  unlikely that the two-dimensional holonomy of the last braided surface of figure \ref{branch} {(Creation of a Branch Point)} could be defined without extending the already defined categorified algebras of diagrams; see \cite{CFM1,CFM2}.)

Nevertheless, the two-dimensional holonomy approach outlined here and in \cite{FMP1,FMP2} can be used to construct algebraic information which can be assigned to the {Reidemeister II and III moves} braided surfaces in figure \ref{2mor}. The flatness of the categorified KZ-connection implies that the associated 2-morphisms do satisfy the movie moves relations, {including the Zamolodchikov tetrahedron equation}, except that we will need to work in a more general (tri-categorical) setting, therefore requiring the {incorporation} of Drinfeld associators.

Since the configuration spaces $\C(n)$ are aspherical, the theory of surface braids without branch points is {uninteresting}. {To have interesting invariants of surface braids we {must  find} the algebraic data to be assigned to the surface braid {``Creation of a Branch Point'' in} figure \ref{branch}. This  fulcral point requires extending the relevant categorified algebras of diagrams, and will be addressed in a future publication.} 

\subsection*{Plan of the paper}
The plan of the paper is the following: 
in Section \ref{2cat} we give the definition of a linear symmetric strict monoidal 2-category (stemming out {from} the definition of a {Baez-Crans} 2-vector space; \cite{BC}), in which we {present} a very strict version of the axioms of a braided monoidal 2-category {(however very well adapted to our purpose)}. Any object $x$ of a  linear monoidal 2-category gives rise to a differential crossed module $\GL(x)$.
In Section \ref{stricti2b}, we discuss the definition of a strict infinitesimal 2-braiding $(r,T)$, which is {a very} {natural and intuitive} categorification of the notion of an infinitesimal braiding $r$, {and is} obtained by imposing that $(r,T)$ is a {(2-categorical) pseudo-natural transformation (see \cite{Pow,SW1,BaH}) rather than a natural transformation  (we will unpack what this means explicitly)}, {together with some technical conditions}. In subsection \ref{c4t} we present the categorified 4-term relations, which {are derived} from the quasi-naturality {(defined later in this paper)} of a strict infinitesimal 2-braiding; see Theorems \ref{main1} and \ref{main2}.  

{Given an infinitesimal 2-braiding $(r,T)$, in a linear symmetric strict monoidal 2-category $\Cc$, and an object of $\Cc$,} we construct a {flat and fake flat} local 2-connection  $(A,B)$, {the Knizhnik-Zamolodchikov 2-connection,} {in the configuration space $\C(n)$, of the form:} 
$$A=\sum_{1\leq a < b \leq n} \w_{ab}r^{ab} \, , \qquad 
B= 2\sum_{1\leq a < b < c \leq n}\w_{bc}\wedge \w_{ca}\,P^{abc} -2\sum_{1\leq a < b < c \leq n}  \w_{ca}\wedge \w_{ab}\,Q^{abc}$$
(see \ref{f4t} for the definition of $P$ and $Q$, which measure the failure of the 4-term {relations \eqref{4t} to hold). In particular, we prove in {Section} \ref{2KZ} that the categorified 4-term relations are equivalent to the vanishing of the 2-curvature {3-form} of $(A,B)$,} {itself one of the strongest results of this paper.}

Finally, in Section \ref{ibdcm}, we work in the framework of differential crossed modules. In subsection \ref{s2c} we categorify the category $\Cc_\lg$ associated to a Lie algebra {$\lg$}, in order to define a linear symmetric monoidal 2-category $\CG$ {(already mentioned above)} associated to a differential crossed module $\lG$. {Similarly to the Lie algebra case, sketched above}, given a (weak) categorical representation of $\lG$ in a chain-complex of vector spaces, we can construct a 2-functor from $\CG$ into the 2-category of weak categorical representations {of $\lG$, these last taking values in the 2-category of chain-complexes of vector spaces, chain maps and chain-homotopies (up to 2-fold homotopy); see \cite{CFM1,CFM2}.} 

In subsection \ref{infbrxmod} we give the definition of a {symmetric} quasi-invariant tensor in a differential crossed module $\lG$, and prove that each of these gives rise to  an infinitesimal 2-braiding in the linear strict {symmetric} monoidal 2-category $\CG$. Finally, we prove in subsection \ref{qitsl2a} that {Wagemann's} {String} differential crossed module \cite{wag06} contains a {non-trivial} {symmetric} quasi-invariant tensor, categorifying the invariant tensor of $\sl(\C)\tn\sl(\C)$ coming from the Cartan-Killing form.
\section{Preliminaries on 2-categories}\label{2cat}

\subsection{Linear 2-categories}
{Recall that a (Baez-Crans) 2-vector space \cite{BC} is a category internal to the category of vector spaces. We therefore have a set of objects $B$ and a set $A$ of morphisms, and the source and target maps $s,t\colon A \to B$, as well as the identity map $i\colon B \to A$, are all to be linear. Therefore the set of composable arrows $A \tensor[_t]{\times}{_s} A$ is a vector space, and we impose that the composition of arrows is linear.  A 2-vector space can uniquely (up to isomorphism) be constructed from its underlying complex $t\colon \ker(s) \subset  A \to B$ of vector spaces. {(We review this important issue in Lemmas \ref{oneone} and \ref{flat}.)}}

{
If we have 2-vector spaces $(s,t\colon A \to B)$ and $(s',t'\colon A' \to B')$ then we have a tensor product 2-vector space of the form: $(s\tn s', t \tn t'\colon A\tn A' \to B \tn B')$, considering the tensor product of identities $i\tn i'\colon B \tn B' \to A \tn A'$; see \cite[Proposition 14]{BC}. With this tensor product we have a symmetric monoidal category.}

{
Recall that a 2-category is nothing more than a category enriched over the monoidal category of small categories, with the cartesian product.}  
\begin{Definition}[Linear 2-category]
 {A linear 2-category is a category enriched over the symmetric monoidal category of 2-vector spaces. 
Unpacking this telegraphic definition, a  linear  2-category $\cC=(C_0,C_1,C_2)$ is given by:}
\begin{enumerate}
 \item A set $C_0$ of objects.
 \item A set $C_1$ of 1-morphisms, together with source and target maps $s,t\colon C_1 \to C_0$, as well as an inclusion map $i\colon C_0 \to C_1$, such that $s \circ i=t\circ i=\id_{C_0}$. 
The elements of $C_1$ will normally be denoted as: $$x \ra{f} y=s(f) \ra{f} t (f).$$
\item A set $C_2$ of 2-morphisms, together with source   and target maps $s',t'\colon C_2 \to C_1$, as well as an inclusion map $i'\colon C_1 \to C_2$. These are such that:
\begin{align*}
    &s' \circ i'=t'\circ i'=\id_{C_1}\,,
    &s \circ s' =s \circ t'\,,
    &&t \circ s' =t \circ t'\,.
   \end{align*}
We define $\ss,\tt\colon C_2 \to C_0$ as $\ss= s \circ s' =s \circ t'$ and $\tt= t \circ s' =t \circ t'$. Each element of $C_2$ {can  be denoted as:}
$$\mor{x}{y}{f}{g}{T} \quad = \hspace*{-.6cm} \mor{\ss(T)}{\tt(T)}{s'(T)}{t'(T)}{T}. $$

\item Given $x,y \in C_0,$ we have a vector space structure on each non-empty hom set $C_1(x,y),$ where:
 $$C_1(x,y)=\{f \in C_1: s(f)=x \an t(f)= y\}.$$
\item Given $x,y,z \in C_0$, if both $C_1(x,y)$ and $C_1(y,z)$ are non-empty, we  have a bilinear map (the 1-composition):
$$(f,g) \in C_1(x,y) \times C_1(y,z) \mapsto fg \in C_1(x,z) $$
This is to be associative and to have units, being given by the map $i$, {namely}: 
$$x \stk{i(x)} x \, .$$
The 1-composition is  denoted in the form:
$$x \stk{f} y \stk{g} z=x\stk{fg} z. $$
Given that $C_1(x,x) \neq \emptyset$, we also have the zero morphism (the null element of the vector space $C_1(x,x)$):
$$x \stk{0_x} x, $$
{and since the composition is bilinear, we have the equation below, for each $f \in C_1(x,x)$:}
$$ x \stk{0_x} x \stk{f}x =x \stk{f} x \stk{0_x} f=x \stk{0_x} x.$$
\item Given $x,y \in C_0$, define:
$$C_2(x,y)=\{T \in C_2: \ss(T)=x \an \tt(T)=y\}.$$
In the case when $C_2(x,y)$ is not empty, we  impose that we are given a vector space structure on {$C_2(x,y)$}, such that  $s',t'\colon C_2(x,y) \to C_1(x,y)$ and also $i'\colon C_1(x,y) \to C_2(x,y)$ each are linear. 
Let $0_{x}^2$ be the null vector of the vector space $C_2(x,x)$. We in particular have that
$i'(0_x)=0^2_x. $

\item Let $x,y \in C_0$. Given $f,g \in C_1$, with $s(f)=s(g)=x$ and $t(f)=t(g)=y$ define the  2-hom set: $$C_2(f,g)=\{T \in C_2: s'(T)=f \an t'(T)=g\}.$$ Given $f,g, h \in C_1(x,y)$ we are to have a vertical composition  map: $$(T,S)\in C_2(f,g)\times C_2(g,h)\mapsto \comp{T}{S} \in C_2(f,h).$$ This will normally be denoted as:
$$\morr{x}{y}{f}{g}{h}{T}{S} \quad = \nss \mor{x}{y}{f}{h}{\comp{T}{S}} \, .$$
This vertical composition is to be associative and to have units, given by $i'\colon C_1 \to C_2$, {namely}:
$$\mor{x}{y}{f}{f}{i'(f)}. $$
 Given $x,y$ with ${C_1(x,y)} \neq \emptyset$, consider the vector space of composable 2-morphisms:
$${\rm Comp}(x,y)=\{(T,S) \in C_2(x,y) \times C_2(x,y): t'(T)=s'(S) \}.$$
{We suppose that the map below is linear:}
$${(T,S) \in {\rm Comp}(x,y) \mapsto \comp{T}{S} \in C_2(x,y)} .$$

\item  Let $x,y,z \in C_0$. Given $f,g \in C_1(x,y)$, we have a  right whiskering map: 
$$(T,h) \in C_2(f,g)\times C_1(y,z)\mapsto Th \in C_2(fh,gh) \, ,$$ 
which will be  visualized as:
$$\rw{x}{y}{f}{g}{T}{z}{h} \quad = \hspace*{-.6cm}\mor{x}{z}{fh}{gh}{Th} \, ,$$
and required to be bilinear as a map $C_2(x,y)\times C_1(y,z) \to C_2(x,z)$. In particular:
$$T0_x=0^2_x \; \an \; 0_x^2f=0_x^2 \, , \qquad \forall \, x\in C_0, f\in C_1(x,x), T\in C_2(x,x) \, .$$

We analogously have a left whiskering map, visualized as:
$$\lw{x}{y}{f}{g}{T}{z}{h} \quad = \hspace*{-.6cm} \mor{z}{y}{hf}{hg}{hT} .$$
These whiskering maps are to be distributive with respect to the 2-composition of 2-morphisms: 
$$\comp{T}{S} f = \comp{Tf}{Sf} \; \an \; f \comp{T}{S} =\comp{fT}{fS}\, , $$
and also associative:
$$f (T g)=(fT) g \, , \quad T(fg)=(Tf)g \, , \quad (fg)T=f(gT) \, .$$
(The previous identities are to hold for any triple of composable morphisms).
In addition, the interchange law is  also to be satisfied: if we are given 1-morphisms:
$$x\stk{f,g} y \; \textrm{ and } \; y \stk{f',g'} z $$
and 2-morphisms:
$$\mor{x}{y}{f}{g}{T} \qquad \an \mor{y}{z}{f'}{g'}{T'}, $$
then we should have:
\begin{equation}
\label{Vcomp}
\comp{fT'}{Tg'}=\comp{Tf'}{g T'} \, . 
\end{equation}
Therefore we can define the horizontal composition of two 2-morphisms, as a particular instance of vertical composition:
$$\hcomp{x}{y}{f}{g}{T}{z}{f'}{g'}{T'} \quad = \hspace*{-.6cm} \mor{x}{z}{ff'}{gg'}{TT'}, $$
where 
\begin{equation}
\label{2hor} 
TT'=\comp{fT'}{Tg'}=\comp{Tf'}{g T'} \, .
\end{equation}
\end{enumerate}
\end{Definition}
\begin{Remark}
Note that we could state conditions 4, 6 and 7 of the definition of a linear 2-category as saying that, given $x,y \in C_0$, then either the category having objects $C_1(x,y)$ and morphisms $C_2(x,y)$ is empty, or it is a Baez-Crans 2-vector space, \cite{BC}. 
\end{Remark}

\begin{Lemma}\label{cwis}
 Let $x,y,z,w \in C_0$. Let $T= \hspace*{-1.1cm}\mor{x}{y}{f}{f'}{T}$, $S= \hspace*{-1.1cm} \mor{z}{w}{g}{g'}{S}$ be 2-morphisms. Given a 1-morphism $y \ra{h} z$ we have:
$$(T h) S=T(hS). $$
\end{Lemma}
\begin{Proof}
By definition of horizontal composition:
 \begin{align*}
  T(h S)=\comp{T(hg)}{f'(hS)}\,, 
 \end{align*}
whereas:
 \begin{align*}
  (Th) S=\comp{(Th)g}{(f'h)S} \, .
 \end{align*}
These coincide due to the associativity of the whiskering. \qed\end{Proof}

Consider an object $x \in C_0$. Let: $$V(x,x)=\{ A\in C_2(x,x): s'(A)=0_{x}\}.$$ 
One has a map:
$$(f,A) \in C_1(x,x)\times V(x,x) \stackrel{\eta}{\mapsto} i'(f)+A \in  C_2(x,x) \, .$$ Since $s'$ and $t'$ are linear maps, we can see that $s'(i'(f)+A)=f$ and $t'(i(f)+A)=f+t'(A)$.
In the opposite direction, given $T$ in $C_2(x,x)$, we can consider the pair: $$\mu(T)=(s'(T),\v{T})\in C_1(x,x)\times V(x,x),$$
where $\v{T}$ is the arrow part of $T$, namely \cite{BC}:
\begin{equation}
\label{defvert}
\v{T}=T-i'(s'(T)) \, .
\end{equation} 
Therefore $T=i(s'(T))+\v{T}$, i.e. we can decompose a 2-morphism into an arrow part plus (the image via $i$ of) the source 1-morphism. Note that $t'(\v{T})=t'(T-i'(s'(T)))=t'(T)-s'(T)$. Pictorially:
$$ \mor{x}{x}{s'(T)}{t'(T)}{T} \quad = \nss \mor{x}{x}{s'(T)}{s'(T)}{i'(s'(T))} \quad + \nss \mor{x}{x}{0_x}{t'(T)-s'(T)}{\v{T}} \, .$$

\begin{Lemma}\label{oneone}
The above correspondence between $C_1(x,x)\times V(x,x)$ and $C_2(x,x)$ is one-to-one.
\end{Lemma}
\begin{Proof}
Taking into account the linearity of $s'$ and $t'$, if $(f,A) \in C_1(x,x)\times V(x,x) $  we have:
\begin{align*} 
 \mu\big(\eta(f,A)\big)&=\mu\big(i'(f)+A\big)= \big (s'(i'(f)+A)), i'(f)+A - i'(s'(i'(f)+A )) \big) \\&= \big (f,i'(f)+A-i'(f)\big) =(f,A) \in C_1(x,x)\times V(x,x)\, .
\end{align*}
On the other hand, if $T \in C_2(x,x)$ we have:
\begin{equation*}
\eta\big (\mu(T)\big)=\eta \big(s(T),T-i'(s'(T)\big)  = i'(s(T))+T-i'(s(T))= T \in C_2(x,x) \, .   
\end{equation*}
\qed\end{Proof}

The next lemma shows how vertical composition of $2$-morphisms behaves under this correspondence.

\begin{Lemma}\label{flat}
If $(A,B)$ is a pair of vertically composable 2-morphisms in $C_2(x,x)$ then:
\begin{equation}
\label{vertvert}
\comp{A}{B} \; = \; \eta(s'(A) \, , \, \v{A}+\v{B}) \; = \; i'(s(A))+\v{A}+\v{B}   
\, . 
\end{equation}
\end{Lemma}
\begin{Proof}
See \cite[pg. 11]{BC}. Since $A$ and $B$ are vertically composable 2-morphisms, $t'(A)=s'(B)$. By linearity of vertical composition we can write:
$$\comp{A}{B}=
\comp{i'(s'(A)) + \v{A}}{i'(s'(B)) + \v{B}} \; = \;  
\comp{i'(s'(A)) + \v{A}}{i'(s'(B))} \; + \; \comp{i'(0_x)}{\v{B}}\, ,
$$
and now it is easy to recognise in the right hand side  $i'(s'(A))+\v{A}$ as the first summand and $\v{B}$ as the second summand. Note that $i'\colon C_1 \to C_2$ is the 2-identity map. 
\qed\end{Proof}

\begin{Remark}
Note that in particular we can see that any 2-morphism is invertible, with respect to the vertical composition. For this reason the inverse of $f \ra{A} g$ will be denoted by $g \ra{-A} f$, thus $\v{(-A)}=- \, (\v{A})$. 
\end{Remark}

\begin{Example}[Chain complexes]The most natural non-trivial example of a linear 2-category is probably the category $\Ch$ of chain complexes of vector spaces, chain maps, and chain homotopies  {up to to fold homotopy (see paragraph below)}, except that we have a class (rather than a set) of objects. As objects of $\Ch$ we take all chain-complexes of vector spaces. Given chain-complexes  $x$ and $y$, the set of 1-morphisms is the vector space of chain-maps $x \to y$, with 1-composition being the composition of chain-maps.

 Given chain complexes $x$ and $y$, of vector spaces, with boundary $\partial$,   recall that a chain-homotopy $s$ is a degree-1 map $x \to y$, and that two chain homotopies $s$ and $t$ are said to be 2-fold homotopic if there exists a degree 2-map $\alpha\colon x \to y$ with $s-t=\alpha \circ \partial-\partial \circ \alpha$. Then given chain maps $f,g \colon x \to y$ the set of two morphisms $f \Rightarrow g$ is the set of pairs $(f,s)$, where $s$ is a homotopy (up to 2-fold homotopy) connecting $f$ and $g$, therefore $f-g=s \circ \partial+\partial \circ s$. The vertical composition of 2-morphisms corresponds to the sum of homotopies, whereas the whiskering is obtained from the composition of homotopies with chain-maps. It is essential to consider homotopies up to {2-fold} homotopy in order that the interchange law holds. For details see \cite{CFM1,CFM2}.
 
\end{Example}

\subsection{Differential crossed modules from linear 2-categories}\label{diffx}
Recall that a differential crossed module \cite{BC,BL} is given by a Lie algebra morphism $\partial \colon \lh \to \lg$ together with a left action $\t$ by derivations of $\lg$ on $\lh$, such that:\begin{itemize}                                                                                                                                                                                                \item[-] $\d(X \t v)=[X,\d(v)]$, for each $v \in \lh$ and $X \in \lg$;
\item[-] $\d(v) \t w=[v,w]$, for each $v,w \in \lh$.                                                                                                                                                                                         \end{itemize}
\noindent
Let $\cC=(C_0,C_1,C_2)$ be a linear 2-category. Given an object $x \in C_0$ one has a differential crossed module:
$$\GL(x)=(\beta\colon \gl^1(x) \to \gl^0(x),\t) \, ,  $$
essentially constructed in \cite{BC,BL}.   The Lie algebra $\gl^0(x)$ {is given by  $C_1(x,x)$,} with bracket: 
$$\{f,g\}=fg-gf $$
(recall that $C_1(x,x)$ is a vector space). The Lie algebra $\gl^1(x)$ is given by all 2-morphisms {$T$ of the form:}
$$\mor{x}{x}{0_x}{f}{T} \,  $$
{(in other words  $T \in V(x,x)$),}
the bracket being: $$\{T,S\}=TS-ST \, .$$
Here $TS$ is the horizontal composition of 2-morphisms. We have an action of $\gl^0(x)$ on $\gl^1(x)$, where:
$$f \t T= f T-Tf. $$
This is an action by derivations due to the associativity of the whiskering. We have a Lie algebra morphism $\beta\colon \gl^1(x) \to \gl^0(x)$, where {$\beta=t'\colon C_2(x,x) \to C_1(x,x)$ is the target map.} Clearly: $$\beta(f \t T)=f \t \beta(T).$$ 
{Let us see that $\beta(T) \t S= \{T,S\}$. Note that $s'(T)=s'(S)=0_x$. Let $f=t'(T)$ and $g=t'(S)$.}
$$
\{T,S\} = TS-ST = \comp{T 0_x}{ f S} - \comp{0_x T}{S f} = \comp{ 0_x^2}{ f S} - \comp{0^2_x}{S f} = 
\comp{ i'(0_x)}{ f S}-\comp{i'(0_x)}{S f}=fS-Sf = \beta(T) \t S \, .
$$
Note that $i'\colon C_1 \to C_2$ is the 2-identity map.

\subsection{Strict monoidal linear 2-categories}
\noindent {In the following definition, the prefix ``strict'' emphasizes the  strictness of the  monoidal structure itself. We note that the notion of a monoidal linear category could be weakened in a plethora of ways. There are some chances that more generality could be needed in order to capture all of the topology of braided surfaces.}
\begin{Definition}
\label{lin2cat}
A strict monoidal linear 2-category $(\cC,\tn,I)=(C_0,C_1,C_2,\tn,I)$ is given by:
\begin{enumerate}
 \item A linear 2-category $\cC=(C_0,C_1,C_2)$.
 \item For any $x,y \in C_0$ an object $x \tn y$ of $C_0$. We are to have that $x \tn (y \tn z)=(x \tn y) \tn z$, {given} objects $x,y,z \in C_0$.
\item An object $I$ of $\cC$. We suppose that $x \tn I=I \tn x=x$, for each object $x$.
 \item For any $x \in C_0$ and any $y \ra{f} z \in C_1$, we have 1-morphisms: 
$$x \tn y  \ra{x \tn f} x \tn z \quad\an\quad  y \tn x  \ra{f \tn x} z \tn x \, .$$ 
Tensoring with an object is to be distributive with respect to composition of morphisms, namely
$$ x\tn (ff')=(x\tn f)(x\tn f') \, ,$$  and 
moreover it is to define a linear map $ C_1(y,z) \to C_1(x\tn y,x\tn z) \, .$
The same is to hold for right tensoring by $x \in C_0$.  Tensoring with $I$ is supposed to not change anything:  $I \tn f=f=f \tn I$ for each 1-morphism $f$. Tensoring morphisms with objects is to be associative, for instance:
$$ x\tn ( f \tn y)= (x \tn f) \tn y \; \an \; (f \tn x) \tn  y=f \tn (x \tn y).$$
\item For any $x \in C_0$ and \hspace*{-1.1cm} $\mor{y}{z}{f}{g}{S} \in C_2$, we are to be given 2-morphisms:
\begin{align*}
x\tn   \left ( \nsl\mor{y}{z}{f}{g}{S}\right) \quad & = \nss \mor{x\tn y}{x\tn z}{x\tn f}{x\tn g}{x \tn S} \, , \\
\left ( \nsl \mor{y}{z}{f}{g}{S}\right)\tn x \quad & = \nss \mor{y\tn x}{z \tn x}{f\tn x}{g\tn x}{S \tn x} \, .
\end{align*}
Moreover, tensoring with an object is to be linear and distributive with respect to vertical composition of 2-morphisms. Tensoring with $I$ is supposed to not change anything. Finally, whiskering commutes with tensoring with objects: given $x,y,z\in C_0$ with $C_2(x,y)$ and $C_1(y,z)$ not empty, we suppose that
$$(A,f) \in C_2(x,y) \times C_1(x,y) \mapsto Af \in C_2(x,y)$$ 
commutes with left and right tensoring. For instance:
$$x\tn (A f)=(x \tn A)(x \tn f) \an (Af) \tn x=(A\tn x) (f \tn x) \, . $$
The same is of course to hold for left whiskering. Tensoring with objects is to be associative.
\item {We now impose some natural interchangeability conditions.}

\begin{enumerate}[(i)]
\item Given $x \ra{f} y$ and $x'\ra{f'} y'$ we impose that:
$$ \big( \, x\tn x' \ra{x \tn f'} x \tn y'  \ra{f \tn y'} y \tn y' \, \big) = 
\big( \, x\tn x' \ra{f \tn x'} y\tn x'  \ra{y \tn f'} x' \tn y' \, \big) \, .$$
This permits us to define the tensor product of morphisms $f\tn f': x\tn x' \rightarrow y\tn y'$ as being one of the previous compositions.
\item Given \nsl $\mor{y}{z}{f}{g}{S}\in C_2$ and $x \ra{h} x'$ we have the following equality of 2-morphisms:
$$\nsl\nsl\rw{y \tn x}{z \tn x} {f \tn x}{g\tn x} {S \tn x}{z \tn x'} {z \tn h} \quad = \nss \lw{y \tn x'}{z \tn x'} {f \tn x'}{g\tn x' } {S \tn x}{y \tn x} {y \tn h}$$
We define:
$$\left( \nsl\mor{y}{z}{f}{g}{S}\right)\tn (x \ra{h} x') \quad = \nss \mor{y\tn  x }{z \tn x'}{f\tn h }{g\tn h}{S \tn h}$$
as being one of the previous compositions. A similar identity is to hold also for left tensoring. 
\item Given:
$$\mor{x}{y}{f}{g}{S} \quad \an \nss \mor{x'}{y'}{f'}{g'}{S'}\,,  $$
we have that:
$$\morr{x \tn x'}{y \tn y'}{f \tn f'}{g \tn f'}{g \tn g'}{S \tn f'}{g \tn S'} \quad = \nss \morr{x \tn x'}{y \tn y'}{f \tn f'}{f \tn g'}{g \tn g'}{f \tn S'}{S  \tn g'} \, . $$
We define:
$$\mor{x \tn x'}{y \tn y'}{f \tn f'}{g \tn g'}{S \tn S'}$$
as being one of the previous compositions.
\end{enumerate}
\end{enumerate}
\end{Definition}
\begin{Example}[Chain-complexes again] The 2-category $\Ch$ of chain complexes is nearly a strict monoidal category, except that the tensor product of objects is not strictly associative; for details, in a different language, see \cite{CFM1,CFM2}. The tensor product of chain-complexes is the usual one; see \cite{D}, as it is the tensor product of chain maps, and the tensor product of chain-maps and homotopies, up to 2-fold homotopy. Note however that given two 2-morphisms $f \stackrel{(f,s)}{\Rightarrow} g$ and $f' \stackrel{(f',s')}{\Rightarrow} g'$ then the tensor product is:
$$\left (f \stackrel{(f,s)}{\Rightarrow} g\right ) \otimes \left (f' \stackrel{(f',s')}{\Rightarrow} g'\right)=f \tn f' \xLongrightarrow{(f \tn f',f\tn s'+s \tn g')} g \tn g' =f \tn f' \xLongrightarrow{(f\tn f', s \tn f'+g \tn s')} g \tn g'.  $$
Again we must consider homotopies up to 2-fold homotopy in order that the last equality holds.
\end{Example}
\subsection{Symmetric strict monoidal linear 2-categories}
This is a very restricted case of the definition of a braided monoidal 2-category in \cite{BN,KV}.
\begin{Definition}
\label{tssml2c} 
A totally symmetric strict monoidal linear 2-category $\cC=(\cC,\tn,I, B)$ is given by:
\begin{enumerate}
\item A strict monoidal linear 2-category $(\cC,\tn,I)$.
 \item For any 2-objects $x$ and $y$ an invertible 1-morphism $x \tn y \ra{B_{x,y}} y \tn x .$ This is to satisfy the following properties:
\begin{enumerate}
\item For any two objects $x$ and $y$ we have: $B_{x,y}B_{y,x}=\id_{x \tn y}.$
\item For any three objects $x,y,z$: 
\begin{equation}
\label{braideq}
\left (x \tn y \tn z \ra{B_{x, y \tn z}}y \tn z \tn x \right)= \left (x \tn y \tn z\ra{ B_{x,y} \tn z} y \tn x \tn z \ra{y \tn B_{x,z}} y \tn z \tn x\right) \, ,
\end{equation}
and analogously for $B_{x \tn y, z}$. Moreover, $B_{I,x}=B_{x.I}=\id_x$.
\item Given {1-morphisms} $x \ra{f} y$ and $x' \ra{f'} y'$ we have:
\begin{equation}
\label{Bfunct1}
\left ( x \tn x' \ra{f \tn f'} y \tn y' \ra{B_{y,y'}} y' \tn y\right)=\left (x \tn x' \ra{B_{x,x'}} x' \tn x \ra{f' \tn f} y' \tn y \right) \, .
\end{equation}
\item Given {2-morphisms} $\nsl\mor{x}{y}{f}{g}{A}$ and $\nsl\mor{x'}{y'}{f'}{g'}{A'}$ we have:
\begin{equation}
\label{Bfunct2}
(A \tn A')B_{(y,y')}=B_{(x,x')} (A' \tn A) \, .
\end{equation}
\end{enumerate}
\end{enumerate}
Note that properties (a), (c), (d) are equivalent to saying that $B$ is an involutive natural transformation between the strict 2-functors  $\otimes$ and $\otimes^{op}$ (the opposite tensor product).  {For the definition of a strict 2-functor and of a natural transformation see \cite[\S 3]{BaH}, or \cite[Definitions 2.3 and 2.4]{Pow}.}
\begin{Example} Forgetting about the associativity morphisms, the {most} basic example of a totally symmetric strict monoidal linear {2-category} is the 2-category of chain complexes, chain-maps, and homotopies up to {2-fold homotopy}. {Given chain complexes $\V$ and $\cW$, the symmetry $\V \tn \cW \to \cW \tn \V$ is the usual graded version of the flip}      \cite{D,Sch}. See \cite{CFM1,CFM2} for details, in a different language. 
\end{Example}

\end{Definition}
The following has an essentially obvious proof, and follows from Mac Lane coherence theorem for braided tensor categories \cite{ML}, in the particular case of symmetric categories.
\begin{Lemma} Consider a totally symmetric strict monoidal linear 2-category $\cC=(\cC,\tn,I, B)$.
Let $a_1,\dots, a_n$ be objects of $\cC$. Given any permutation $\sigma$ of $\{1,\dots ,n\}$ there exists a unique map: $$M_\sigma \colon a_1 \tn \dots \tn a_n \to a_{{\sigma}(1)}\tn \dots \tn a_{{\sigma}(n)}$$ constructed from compositions of the braidings $B_{x,y}$ (where $x$ and $y$ are tensor products of objects in $\{a_1,\dots, a_n\}$) possibly tensored by tensor products of objects in  $\{a_1,\dots, a_n\}.$
\end{Lemma}
\noindent
As an example for $n=3$ and $\s=(13)$ then $$x \tn y \tn z \ra{M_\s} z \tn y \tn x$$
can be defined as (for instance):
$$B_{x, y \tn z} (B_{y \tn z} \tn x )=(x \tn B_{y,z}) B_{x, z \tn y} =(B_{x, y} \tn z) (y \tn B_{x,z}) (B_{y \tn z} \tn x ).$$ 

We now introduce a notation, which will be used in the following section, to  give the definition of a strict infinitesimal 2-braiding. Let $x_1, \dots, x_n$ be objects of $\cC=(\cC,\tn,I, B)$, a totally symmetric strict monoidal linear 2-category. Given distinct indices $a$ and $b$ in $\{1,\dots,n\}$, {consider a morphism $r\colon x_a \tn x_b \to x_a \tn x_b$.} Put:
$$x_1 \tn x_{2} \tn \dots\tn  x_{n} \ra{r^{ab}}  x_{1} \tn x_{2}\tn \dots \tn x_{n}$$
{as being the composition:}
\begin{equation}\label{insdef}
{x_{1} \tn x_{2}\tn \dots \tn x_{n} \ra{M_{\sigma}}  X \tn x_a \tn x_b \tn Y \ra{X \tn r  \tn Y}  X \tn x_a \tn x_b \tn Y  \ra{M_{\sigma^{-1}}}  x_{1} \tn x_{2}\tn \dots \tn x_{n}.}
\end{equation}
Here $\s\colon \{1,\dots,n\} \to  \{1,\dots,n\}$ is a permutation which makes $\s(a)$ and $\s(b)$ contiguous, that is:
$${x_{\s(1)} \tn \dots \tn x_{\s(n)}= x_{\s(1)} \tn \dots\tn x_{\s(i)}\tn x_a\tn x_b \tn x_{\s(i+2)}\tn  \dots\tn x_{\s(n)} .}$$
{We leave the proof of the following easy Lemma to the reader.} 
\begin{Lemma}\label{coher}
The above definition of $r^{ab}$ does not depend on the permutation chosen. 
\end{Lemma}

{This type of notation extends immediately: Given $R\colon x_a \tn x_b \tn x_c \to   x_a \tn x_b \tn x_c  $ we can define:}
$$R^{abc}\colon x_{1} \tn x_{2}\tn \dots \tn x_{n} \to  x_{1} \tn x_{2}\tn \dots \tn x_{n} \, ,$$
{by using the obvious analogue of \eqref{insdef}, and this {too} does not depend on the permutation chosen. 
Analogously, given a 2-morphism:}
$$\mor{ x_a \tn x_b \tn x_c }{ x_a \tn x_b \tn x_c }{R}{R'}{L} \, ,$$
we can unambiguously define a 2-morphism:
$$\mor{ x_{1} \tn x_{2}\tn \dots \tn x_{n} }{   x_{1} \tn x_{2}\tn \dots \tn x_{n} }{R^{abc}}{R'^{abc}}{L^{abc}} \, .$$

\section{Strict infinitesimal 2-braidings}\label{stricti2b}

\subsection{Definition of a strict infinitesimal 2-braiding}
We present a direct categorification of the concept of an infinitesimal braided category; see \cite[XX.4]{K} and \cite{C}, called there infinitesimal symmetric category. 

\begin{Definition}
\label{definf2br}
A strict infinitesimal 2-braiding $(r,T)$ in a symmetric strict monoidal linear 2-category  $\cC=(\cC,\tn,I, B)$ is given by the following data:
\begin{enumerate}
\item For any pair of objects $x$ and $y$, a 1-morphism $r_{x,y}\colon x \tn y \to x \tn y$.  
These are to satisfy the linearity  conditions:
\begin{equation}
\label{rlin}
r_{x,y \tn z} = r_{x,y}^{12}+r_{x,z}^{13} \, , \qquad r_{x \tn y,z}=r_{x,z}^{13}+r_{y,z}^{23}\, ,
\end{equation}
thus for example: \quad \quad 
$r_{x,y \tn z \tn w}   = r_{x,y}^{12}+r_{x,z}^{13}+r_{x,w}^{14} \; , \qquad
r_{x \tn y,  z \tn w} = r_{x,z}^{13}+r_{x,w}^{14}+r_{y,z}^{23}+r_{y,w}^{24} \, .  
$
\item For any morphism $x \ra{f} x'$ and any object $y$, a 2-morphism $T_{(f,y)}$ such that:
\begin{equation}
\label{defT}
\defT{x\tn y}{x'\tn y}{f\tn y}{r_{x',y}}{r_{x,y}}{T_{(f,y)}}   
\end{equation}
{Therefore $T(f,y)$} measures the failure of functoriality of $r$ with respect to the morphism $f\otimes y$. This is to be natural with respect to 2-morphisms. If $S:f\implies g$ then $T_{(f,y)}$ and $T_{(g,y)}$ have to be related via a vertical composing with $S$. {The request is:}
\begin{equation}
\label{natT}
\comp{T_{(f,y)}}{(S\tn y)r_{x',y}} = \comp{r_{x,y}(S\tn y)}{T_{(g,y)}}\, ,
\end{equation}
which is graphically visualized as:
\begin{equation}
\label{natTgr}
\vcenter{\vbox{ 
\xymatrix{
 & & x'\tn y \ar@/^1pc/[rrd]^{r_{x',y}} & & \\ 
 x \tn y \ar@/_1pc/[rrd]_{r_{x,y}} \ar@/^1pc/[rru]^{g \tn y } \ar@/_1pc/[rru]_{f \tn y} \ar@{{}{ }{}} [rru]|{ \Uparrow S \tn y} & & \Uparrow T_{(f,y)} & & x' \tn y \\
 & & x \tn y \ar@/_1pc/[rru]_{f \tn y} & &
}
}} 
\quad =  \vcenter{\vbox{ 
\quad \xymatrix{
 & & x' \tn y \ar@/^1pc/[rrd]^{r_{x',y}} & & \\ 
 x \tn y \ar@/_1pc/[rrd]_{r_{x,y}}\ar@/^1pc/[rru]^{g \tn y }  & & \Uparrow T_{(g,y)} & & x' \tn y\\
& & x \tn y \ar@/^1pc/[rru]^{g \tn y} \ar@{{}{ }{}} [rru]|{ \Uparrow S \tn y} \ar@/_1pc/[rru]_{f \tn y} & & 
} }} \, .
\end{equation}
Moreover we suppose that $T$ satisfies the linearity {condition:}
\begin{equation}
\label{linT}
T_{(f+f',y)}=T_{(f,y)}+T_{(f',y)},
\end{equation}
and that, if $x \ra{f} x'\ra{f'} x''$, we have:
\begin{equation}
\label{compT}
T_{ff',y}=\comp{T_{(f,y)}(f'\tn y)}{(f \tn y) T_{(f',y)}} \, . 
\end{equation}
Graphically this means:
$$
\vcenter{\vbox{ 
\defT{x\tn y}{x''\tn y}{ff'\tn y}{r_{x'',y}}{r_{x,y}}{T_{(ff',y)}}
}} 
\quad = \quad  
\vcenter{\vbox{ 
\xymatrix{
 & & x''\tn y \ar@/^1pc/[rd]^{r_{x'',y}} & \\
x\tn y \ar@/_1pc/[rd]_{r_{x,y}} \ar[r]^{f\tn y} & x'\tn y \ar@/^1pc/[ru]^{f'\tn y} \ar[r]^{r_{x',y}} & 
x'\tn y \ar[r]^{f'\tn y} \ar@{{}{ }{}} [u]|{\Uparrow T_{(f',y)}} & x''\tn y \\
 & x\tn y \ar@{{}{ }{}} [u]|{\Uparrow T_{(f,y)}} \ar@/_1pc/[ru]_{f\tn y} & 
} 
}} \, .
$$
\item
Similarly, for any object $x$ and any 1-morphism $y \ra{f} y'$, a 2-morphism $T_{(x,f)}$:
$$
\defT{x\tn y}{x\tn y'}{x\tn f}{r_{x,y'}}{r_{x,y}}{T_{(x,f)}} \, .$$
As before, this is to be natural with respect to all 2-morphisms $S \colon f \implies g$, in the sense that
\begin{equation}
\label{natTr}
\comp{r_{x,y}(x\tn S)}{T_{(x,f)}} = \comp{T_{(x,g)}}{(x\tn S)r_{x,y'}} \, .
\end{equation} 
Moreover, the  obvious analogues of linearity \eqref{linT} and composition \eqref{compT} properties hold. 
\item
We suppose that that given 1-morphisms $x\ra{f} x'$ and $y \ra{g} y'$, the 2-morphisms $T_{(f,y)}$ and $T_{(x,g)}$ satisfy the following interchange law:
\begin{equation}
\label{Tfg} 
\comp{T_{(f,y)}(x'\tn g)}{(f\tn y)T_{(x',g)}} = \comp{T_{(x,g)}(f\tn y')}{(x\tn g)T_{(f,y')}} \, .
\end{equation}
Graphically: 
$$
\vcenter{\vbox{ 
\xymatrix @=2.5em{
 & & x'\tn y' \ar@/^1pc/[rd]^{r_{x',y'}} & \\
x\tn y \ar@/_1pc/[rd]_{r_{x,y}} \ar[r]^{f\tn y} & x'\tn y \ar@/^1pc/[ru]^{x'\tn g} \ar[r]^{r_{x',y}} & 
x'\tn y \ar[r]^{x'\tn g} \ar@{{}{ }{}} [u]|{\Uparrow T_{(x',g)}} & x'\tn y' \\
 & x\tn y \ar@{{}{ }{}} [u]|{\Uparrow T_{(f,y)}} \ar@/_1pc/[ru]_{f\tn y} & 
} 
}} \quad = \quad 
\vcenter{\vbox{ 
\xymatrix @=2.5em{
 & & x'\tn y' \ar@/^1pc/[rd]^{r_{x',y'}} & \\
x\tn y \ar@/_1pc/[rd]_{r_{x,y}} \ar[r]^{x\tn g} & x\tn y' \ar@/^1pc/[ru]^{f\tn y'} \ar[r]^{r_{x,y'}} & 
x\tn y' \ar[r]^{f\tn y'} \ar@{{}{ }{}} [u]|{\Uparrow T_{(f,y')}} & x'\tn y' \\
 & x\tn y \ar@{{}{ }{}} [u]|{\Uparrow T_{(x,g)}} \ar@/_1pc/[ru]_{x\tn g} & 
} 
}}
$$The above equality allows to define unambiguously $T_{(f,g)}$ as one of the two sides of \eqref{Tfg}. 
\item
Finally we suppose that  given a 1-morphism $f\colon x \to x'$ the following linearity conditions hold:
\begin{align}
\label{Tlinob}
T_{(f,y \tn z)} & = T^{12}_{(f,y)}+T^{13}_{(f,z)} & T_{(y \tn z,f)} & = T^{13}_{(y,f)}+T^{23}_{(z,f)} \\
\label{Tlinmor}
T_{(f \tn y,z)} & = T^{13}_{(f,z)} & T_{(y \tn f,z)} & = T^{23}_{(f,z)}
\end{align}
\end{enumerate}
\end{Definition}
\begin{Remark}
{ Disregarding the linearity conditions, {the previous list of properties means that $(r,T)$ is a pseudo-natural transformation  from the 2-functor} $\tn \colon \Cc \times \Cc\to \Cc$ to itself. Pseudo natural transformations between 2-functors are defined in \cite{Pow,BaH} and \cite[Appendix A1]{SW1}.}
\end{Remark}

\begin{Remark}\label{sd} Equations \eqref{rlin}, \eqref{Tlinob} and \eqref{Tlinmor} could easily be weakened to hold only up to a coherent 2-morphism, which is the type of generality appearing in the definition of a braided-monoidal 2-category; \cite{BN,KV}. We will however not need this generality. This is nevertheless the reason for the term ``strict`` in ``strict infinitesimal 2-braiding''.\end{Remark}

\subsection{The failure of the 4-term relations and coherent infinitesimal 2-braidings}\label{f4t}

Remember that in an infinitesimally symmetric category the natural endomorphism $r:\tn \Rightarrow\tn$ satisfies the so-called 4-term relations or infinitesimal braid relations: for  {every triple of} objects $x,y,z$ we have
$$ [r^{12}_{x,y},r^{13}_{x,z}+r^{23}_{y,z}] = 0 \, , \qquad [r^{12}_{x,y}+r^{13}_{x,z},r^{23}_{y,z}] = 0 \, .$$
These relations express the functoriality of $r_{x,y\tn z}$ (resp. $r_{x\tn y,z}$) with respect to the morphism $r_{x,y}\tn z$ (resp. $x\tn r_{y,z}$).   
Consider now a totally symmetric monoidal linear strict 2-category, with an infinitesimal 2-braiding $(r,T)$. Since $(r,T)$ is a pseudo-natural transformation, $r$ is not functorial and so it does not satisfy the 4-term relations. We can however measure this failure in terms of the 2-morphism $T$ suitably evaluated.  Let $x,y,z \in C_0$. Let us consider the 2-morphism $P_{x,y,z}:=T_{(x,r_{y,z})}\, $. Explicitly: 
\begin{equation}
\label{defP}
r_{x,y \tn z}( x \tn r_{y,z}) \xRightarrow{P_{x,y,z} \, = \, T_{(x,r_{y,z})}}  r_{y,z}  (r_{x,y \tn z}) \, .
\end{equation}
From its definition, and using the linearity of $r$:
\begin{equation}
\label{diagP}
\vcenter{\vbox{
\defT{x\tn y\tn z}{x\tn y\tn z}{x\tn r_{y,z}}{r_{x,y\tn z}}{r_{x,y\tn z}}{T_{(x,r_{y,z})}}
}} \quad = \quad
\vcenter{\vbox{
\xymatrix{
x\tn y\tn z \ar@/^26pt/[rr]^{r_{y,z}^{23} r_{x,y}^{12}+ r_{y,z}^{23}r_{x,z}^{13}} \ar@/_26pt/[rr]_{r_{x,y}^{12}r_{y,z}^{23}+r_{x,z}^{13}r_{y,z}^{23}} & \Uparrow P_{x,y,z} & x\tn y\tn z
}
}}  \, . 
\end{equation}
Similarly we consider the 2-morphism  $Q_{x,y,z}:=T_{(r_{x,y},z)}\,$. Explicitly:
\begin{equation}
\label{defQ}
r_{x\tn y,z} (r_{x,y} \tn z) \xRightarrow{T_{(r_{x,y}\tn z)} \, = \, Q_{x,y,z}} (r_{x,y} \tn z) r_{x \tn y,z} \, . 
\end{equation} 
As before, by linearity of $r$ it is :
\begin{equation}
\label{diagQ}
\vcenter{\vbox{
\defT{x\tn y\tn z}{x\tn y\tn z}{r_{x,y}\tn z}{r_{x\tn y,z}}{r_{x\tn y,z}}{T_{(r_{x,y},z)}}
}} \quad = \quad
\vcenter{\vbox{
\xymatrix{
x\tn y\tn z \ar@/^26pt/[rr]^{r_{x,y}^{12}r_{x,z}^{13} +r_{x,y}^{12}r_{y,z}^{23}} \ar@/_26pt/[rr]_{r_{x,z}^{13} r_{x,y}^{12}+r_{y,z}^{23}r_{x,y}^{12}} & \Uparrow Q_{x,y,z} & x\tn y\tn z
}
}}  \, . 
\end{equation}
We then see that $P_{x,y,z}$ and $Q_{x,y,z}$ are 2-morphisms connecting the left and right sides of the 4-term relations.

Closely related to these, there are three naturally defined 2-morphisms, which connect a combination of 4-term relations sides: 
$$ r^{23}_{yz}r^{13}_{xz}+r^{23}_{yz}r^{12}_{xy}+r^{12}_{xy}r^{13}_{xz}\Longrightarrow r^{13}_{xz} r^{23}_{yz}+ r^{23}_{yz}r^{12}_{xy}+r^{13}_{xz}r^{12}_{xy} \, .$$
They are:
$$P^{2,1,3}_{y,x,z} + i'(r^{23}_{yz}r^{12}_{xy}) \, , \qquad Q^{1,3,2}_{x,z,y}+i'(r^{23}_{yz}r^{12}_{xy}) \, ,$$
{as well as the composition:}
$$r^{23}_{yz}r^{13}_{xz}+r^{23}_{yz}r^{12}_{xy}+r^{12}_{xy}r^{13}_{xz}\xRightarrow{-P_{x,y,z}+i'( r^{12}_{xy}r^{13}_{xz}) } r^{13}_{xz}r^{23}_{yz}+r^{12}_{xy}r^{23}_{yz}+r^{12}_{xy}r^{13}_{xz} \xRightarrow{-Q_{x,y,z}+i'(r^{13}_{xz}r^{23}_{yz}) } r^{13}_{xz}r^{23}_{yz}+r^{23}_{yz}r^{12}_{xy}+r^{13}_{xz}r^{12}_{xy}. $$

\begin{Definition}[Coherent infinitesimal 2-braiding]\label{coherent}
 An infinitesimal 2-braiding is called coherent if the previous three morphisms all coincide. Passing to the arrow parts of each 2-morphism, this is to say that for {each triple of} objects $x,y,z$ we must have:
$$ \v{P^{213}_{y,x,z}}=-\left(\v{P_{x,y,z}}+\v{Q_{x,y,z}}\right)=\v{Q^{132}_{x,z,y}} .$$
\end{Definition}
 {This axiom of coherence is an infinitesimal analogue of an axiom in the definition of braided monoidal 2-categories, which appears in \cite[Definition 6]{BN}. Namely the axiom $S^+=S^-$, giving the  equality of two hexagons. We mention that this axiom is not in the original definition of a braided monoidal 2-category (\cite{KV}) and that it is slightly less natural categorically, however being entirely justified, geometrically.}

\subsection{Totally symmetric infinitesimal 2-braiding}
\begin{Definition}
\label{inf2br}
A {strict} infinitesimal 2-braiding in a symmetric {strict} monoidal linear 2-category $\cC=(\cC,\tn,I, B)$ is called totally symmetric if for any objects $x,x',y,y',z \in C_0$ and morphisms $x\ra{f} x'$, $y\ra{g} y'$ the following conditions hold: 
\begin{enumerate}[(i)]
\item $r$ and $T$ are functorial with respect to the braiding $B$: 
\begin{align}
\label{rsymm}
r_{x,y} & = B_{x,y}r_{y,x}B_{y,x}  \\ 
\label{TsymmB}
T_{(f,g)} & = B_{x,y} T_{(g,f)} B_{y',x'}   
\end{align}
\item $r$ is functorial with respect to $B$ in the second argument:
 \begin{equation} \label{tid}r_{x,y \tn z} ( x \tn  B_{y,z})=( x \tn  B_{y,z}) r_{x,z \tn y}
\end{equation}
and furthermore:
\begin{equation} 
\label{Tid} T_{(x,B_{y,z})}=i'\big({r_{x,y\tn z}\, (x\tn B_{y,z}})\big)\,.
\end{equation}
(in the equation above recall that $i'\colon C_1 \to C_2$ is the 2-identity map.) 
\item Similarly, $r$ is functorial with respect to $B$ in the first argument:
\begin{equation} r_{x \tn y , z} (  B_{x,y} \tn z)=(   B_{x,y} \tn z) r_{y \tn x , z}
\end{equation}
and furthermore
\begin{equation}
T_{(B_{x,y},z)}  =i'\big({(B_{x,y}\tn z)\,r_{x\tn y,z}}\big).
\end{equation}
\end{enumerate}
\end{Definition}

Some comment on the previous properties: the graphical presentation of \eqref{TsymmB} is: 
\begin{equation*}
\label{st1}
\vcenter{\vbox{
\xymatrix @R=1em @C=2em{
&& y' \tn x'  \ar@/^1pc/[rd]^{r_{y',x'}} & \\ 
x \tn y \ar[r]^{B_{x,y}} & y \tn x \ar@/_1pc/[rd]_{r_{y,x}}\ar@/^1pc/[ru]^{g \tn f } &  \quad \Uparrow T_{(g,f)} \quad &y' \tn x'\ar[r]^{B_{y',x'}} &x' \tn y' \\
&& y \tn x \ar@/_1pc/[ru]_{g \tn f}&
}
}} \quad  = \quad 
\vcenter{\vbox{
\defT{x\tn y}{x'\tn y'}{f\tn g}{r_{x',y'}}{r_{x,y}}{T_{(f,g)}}
}} \, .
\end{equation*}
\noindent
A consequence of  \eqref{Tid} is that $s'(T_{(x,B_{y,z})}) = t'(T_{(x,B_{y,z})})$, where
\begin{align*}
s'(T_{(x,B_{y,z})}) & = (r_{x,y\tn z})(x\tn B_{y,z}) = 
(r^{12}_{x,y} + r^{13}_{x,z}) B^{23}_{y,z} \, ,\\
t'(T_{(x,B_{y,z})}) & = (x\tn B_{y,z})(r_{x,y\tn z}) = B^{23}_{y,z}(r^{12}_{x,y} + r^{13}_{x,z}) 
\, ,
\end{align*}
which is exactly equation \eqref{tid}, while \eqref{Tid} expresses the stronger property that $T_{(x,B_{y,z})}$ is in fact the identity 2-morphism. Combining this with \eqref{compT} we get
\begin{equation}
\label{TBcomp}
B^{23}_{y,z}T_{(x,z\tn g)}B^{23}_{z,y'} = T_{(x,g\tn z)} \, .
\end{equation}

\begin{Lemma}
Suppose that the infinitesimal 2-braiding is totally symmetric. We have, for each three objects $x,y,z$:
\begin{equation}
\label{symmPQ}
P_{x,z,y} = P_{x,y.z}^{132} \, , \quad P_{x,y,z} = Q_{z,y,z}^{321} \, , \quad Q_{x,z,y} = Q_{y,x,z}^{213} 
\end{equation}
\end{Lemma}
\begin{Proof}The first equality is proved as:
\begin{align*}
 P_{x,z,y}&=T_{(x,r_{z,y})}=(x \tn B_{z,y})(x \tn B_{y,z}) T_{(x,r_{z,y})}   (x \tn B_{z,y})(x \tn B_{y,z})=(x \tn B_{z,y}) T_{(x, B_{y,z} r_{z,y}  B_{z,y})}(x \tn B_{y,z})\\
&=(x \tn B_{z,y}) T_{(x,  r_{y,z})}(x \tn B_{y,z})=(x \tn B_{z,y}) P_{x,y,z}(x \tn B_{y,z})= P_{x,y,z}^{132}
\end{align*}
where we used \eqref{Tid} in the second step. The  third equality follows similarly. The second equality is proved, using \eqref{TsymmB}, as:
\begin{equation*}
 P_{x,y,z}=T_{(x,r_{y,z})}=T_{(r_{y,z},x)}^{231}=Q_{y,z,x}^{231}=Q_{z,y,x}^{321} \, . 
\end{equation*}  
\qed\end{Proof}
Combining with Definition \ref{coherent} we reach the following result whose easy proof is omitted.
\begin{Lemma}[Jacobi type identity] Given a coherent totally symmetric infinitesimal 2-braiding $\cC=(\cC,\tn,I, B)$, for each $x,y,z \in C_0$ we have:
\begin{equation}
\label{jac}
P_{x,y,z}^{123}+P_{z,x,y}^{312}+P_{y,z,x}^{231} = 0 \, . 
\end{equation}
\end{Lemma}
\subsection{The categorified 4-term relations}\label{c4t}

Let us denote for simplicity the left and right sides of the 4-term relations as follows:
\begin{align*}
r_{x,y}^{12} r_{y,z}^{23}+r_{x,z}^{13}r_{y,z}^{23} & = s'(P_{x,y,z}) \, , & 
r_{y,z}^{23} r_{x,y}^{12}+ r_{y,z}^{23}r_{x,z}^{13} & = t'(P_{x,y,z}) \, ,\\
r_{x,z}^{13} r_{x,y}^{12}+r_{y,z}^{23}r_{x,y}^{12} & = s'(Q_{x,y,z}) \, , &
r_{x,y}^{12}r_{x,z}^{13} +r_{x,y}^{12}r_{y,z}^{23} & = t'(Q_{x,y,z}) \, .
\end{align*}
These  4-term relations are therefore satisfied if $P_{x,y,z}$ and $Q_{x,y,z}$ reduce to the vertical identity 2-morphisms between their equal source and target. When $P_{x,y,z}$ and $Q_{x,y,z}$ are non-trivial, let us  consider the 2-morphisms $T_{(s'(P_{x,y,z}),w)}$ and $T_{(t'(P_{x,y,z}),w)}$. From the definition of an infinitesimal 2-braiding, these 2-morphisms will satisfy five naturality conditions, which we now describe.

Since $s'(P_{x,y,z}) \stackrel{P_{x,y,z}}{\implies} t'(P_{x,y,z})$, the two $T$ morphisms are related by the naturality condition \eqref{natT}{: Every time we have a 2-morphism $ s'(A) \stackrel{A}{\implies} t'(A)$}
the 2-morphisms $T_{(s'(A),w)}$ and $T_{(t'(A),w)}$ are to satisfy:
$$ 
\comp{T_{(s'(A),w)}}{(A\tn w)r_{X,w}} = \comp{r_{X,w}(A\tn w)}{T_{(t'(A),w)}} \, .$$
We are now considering the case $X=x\tn y\tn z$ and $A=P_{x,y,z}$. This leads to the first naturality condition:
\begin{equation}
\label{nat1}
\comp{T_{(s'(P_{x,y,x}),w)}}{(P_{x,y,z}\tn w)r_{x\tn y\tn z,w}} = \comp{r_{x\tn y\tn z,w}(P_{x,y,z}\tn w)}{T_{(t'(P_{x,y,z}),w)}}\,.
\end{equation}
This can be visualized as:
\begin{multline}
\allowdisplaybreaks
\label{nat1diag}
\vcenter{\vbox{
\xymatrix{
 & & x\tn y\tn z\tn w \ar@/^1pc/[rrd]^{r_{x\tn y\tn z,w}} & & \\
x\tn y\tn z\tn w 
\ar@/_10pt/[rrd]_{r_{x\tn y\tn z,w}} 
\ar@/^10pt/[rru]^{(r_{y,z}^{23} r_{x,y}^{12}+ r_{y,z}^{23}r_{x,z}^{13})\tn w} 
\ar@/_10pt/[rru]_{(r_{x,y}^{12} r_{y,z}^{23}+ r_{x,z}^{13}r_{y,z}^{23})\tn w}
\ar@{{}{ }{}} [rru]|{ \Uparrow P_{x,y,z} \tn w} & &
\qquad\qquad \Uparrow T_{((r_{x,y}^{12} r_{y,z}^{23}+ r_{x,z}^{13}r_{y,z}^{23}),w)} & &
x\tn y\tn z\tn w \\
 & & x\tn y\tn z\tn w \ar@/_1pc/[rru]_{(r_{x,y}^{12} r_{y,z}^{23}+ r_{x,z}^{13}r_{y,z}^{23}) \tn w)} & &
}
}} \quad = \\[.4cm]
= \qquad \vcenter{\vbox{
\xymatrix{
 & & x\tn y\tn z\tn w \ar@/^1pc/[rrd]^{r_{x\tn y\tn z,w}} & & \\
x\tn y\tn z\tn w 
\ar@/_10pt/[rrd]_{r_{x\tn y\tn z,w}} 
\ar@/^10pt/[rru]^{(r_{y,z}^{23} r_{x,y}^{12}+ r_{y,z}^{23}r_{x,z}^{13})\tn w} & &
\Uparrow T_{((r_{y,z}^{23} r_{x,y}^{12}+ r_{y,z}^{23}r_{x,z}^{13}),w)} \qquad\qquad & &
x\tn y\tn z\tn w \\
 & & x\tn y\tn z\tn w 
\ar@/_10pt/[rru]_{(r_{x,y}^{12} r_{y,z}^{23}+ r_{x,z}^{13}r_{y,z}^{23}) \tn w)} 
\ar@/^10pt/[rru]^{(r_{y,z}^{23} r_{x,y}^{12}+ r_{y,z}^{23}r_{x,z}^{13})\tn w}
\ar@{{}{ }{}} [rru]|{ \Uparrow P_{x,y,z} \tn w}
& &
}
}}
\end{multline}
We now derive the relations between $r$ and $P_{x,y,z}$ associated to \eqref{nat1}. It is convenient to use the decomposition of 2-morphisms into the source and arrow part, $A=i'(s'(A))+\v{A}$. We then observe that the source parts of the two vertical compositions in \eqref{nat1} agree; hence only relations between the arrow parts of the vertical compositions are left, and they are computed using Lemma \ref{flat}: 
\begin{equation*}
\begin{split}
\v{\comp{T_{(s'(P_{x,y,x}),w)}}{(P_{x,y,z}\tn w)r_{x\tn y\tn z,w}}} 
& = \, \v{ T_{(s'(P_{x,y,x}),w)} } + \v{ (P_{x,y,z}\tn w)r_{x\tn y\tn z,w} }  \,\,\, ,
\\[.3cm]
\v{\comp{r_{x\tn y\tn z,w}(P_{x,y,z}\tn w)}{T_{(t'(P_{x,y,z}),w)}}} 
& = \, \v{r_{x\tn y\tn z,w}(P_{x,y,z}\tn w)} + \v{T_{(t'(P_{x,y,z}),w)}} \,\,\, .
\end{split}
\end{equation*}
Now we compute each term separately. We use linearity of $r$, \eqref{compT} and $\v{fA}=f\,\v{A}$. So, for example, we have that $\v{T_{(ff',y)}} = f\,\v{T_{(f',y)}} + \v{T_{(f,y)}}f'$. We obtain: 
\begin{align*}
\v{T_{(t'(P_{x,y,z}),w)}} = \v{T_{ ( (r_{y,z}^{23} r_{x,y}^{12}+ r_{y,z}^{23}r_{x,z}^{13}),w ) }} & = 
\v{T_{ (r_{y,z}^{23} r_{x,y}^{12},w )}} + \v{ T_{ (r_{y,z}^{23}r_{x,z}^{13},w )}} \\
 & = \v{T_{ (r_{y,z}^{23},w)}} \, r_{x,y}^{12} + \v{T_{ (r_{y,z}^{23},w)}} \, r_{x,z}^{13} + r_{y,z}^{23} \,  
 \v{T_{ ( r_{x,y}^{12},w )}}+ r_{y,z}^{23} \, \v{T_{ ( r_{x,z}^{13},w )}} \\
 & = \v{Q_{y,z,w}^{234}} \, r_{x,y}^{12} + \v{Q_{y,z,w}^{234}} \, r_{x,z}^{13} + r_{y,z}^{23} \, \v{Q_{x,y,w}^{124}} + r_{y,z}^{23} \, \v{Q_{x,z,w}^{134}}\, ,
\end{align*}
and analogously:
$$ 
\v{T_{(s'(P_{x,y,z}),w)}} = \v{T_{ (( r_{x,y}^{12}r_{y,z}^{23}+ r_{x,z}^{13}r_{y,z}^{23}),w )}} = r_{x,y}^{12} \, \v{Q_{y,z,w}^{234}}+ r_{x,y}^{13} \, \v{Q_{y,z,w}^{234}} + \v{Q_{x,y,w}^{124}} \, r_{y,z}^{23}+ \v{Q_{x,z,w}^{134}} \, r_{y,z}^{23} \ . $$
Also:
\begin{align*}
\v{ (P_{x,y,z}\tn w)r_{x\tn y\tn z,w} } & = \v{P^{123}_{x,y,z}} \, (r_{x,w}^{14}+r_{y,w}^{24}+r_{z,w}^{34}) \, , \\[.4cm] 
\v{ r_{x\tn y\tn z,w}(P_{x,y,z}\tn w) } & = (r_{x,w}^{14}+r_{y,w}^{24}+r_{z,w}^{34}) \, \v{P^{123}_{x,y,z}} \, .
\end{align*}
Putting everything together we reach the first naturality condition:
\begin{equation}
\label{nat1eq} 
(r_{x,w}^{14}+r_{y,w}^{24}+r_{z,w}^{34}) \t \v{ P_{x,y,z}^{123}} -
(r_{x,y}^{12} + r_{x,y}^{13})\t \v{Q_{y,z,w}^{234}}
+ r_{y,z}^{23} \t (\v{Q_{x,y,w}^{124}}+ \v{Q_{x,z,w}^{134}})=0 \, . 
\end{equation}
A second naturality condition is provided by the relation between $T_{(w,s'(P_{x,y,z}))}$ and $T_{(w,t'(P_{x,y,z}))}$. Again derived from \eqref{natT}, it reads:
\begin{equation}
\label{nat2}
\comp{ T_{(w,s'(P_{x,y,z}))} }{ (w\tn P_{x,y,z}) \, r_{w,x\tn y\tn z} } \, = \, 
\comp{ r_{w,x\tn y\tn z} \, (w\tn P_{x,y,z}) }{ T_{(w,t'(P_{x,y,z}))} } 
\end{equation}  
and it is graphically visualizes as: 
\begin{multline}
\allowdisplaybreaks
\label{nat2diag}
\vcenter{\vbox{
\xymatrix{
 & & w\tn x\tn y\tn z \ar@/^1pc/[rrd]^{r_{w, x\tn y\tn z}} & & \\
w\tn x\tn y\tn z 
\ar@/_10pt/[rrd]_{r_{w,x\tn y\tn z}} 
\ar@/^10pt/[rru]^{w\tn (r_{y,z}^{34} r_{x,y}^{23}+ r_{y,z}^{34}r_{x,z}^{24})} 
\ar@/_10pt/[rru]_{w\tn (r_{x,y}^{23} r_{y,z}^{34}+ r_{x,z}^{24}r_{y,z}^{34})}
\ar@{{}{ }{}} [rru]|{ \Uparrow w \tn P_{x,y,z} } & &
\qquad\qquad \Uparrow T_{(w,(r_{x,y}^{23} r_{y,z}^{34}+ r_{x,z}^{24}r_{y,z}^{34}))} & &
w\tn x\tn y\tn z \\
 & & w\tn x\tn y\tn z \ar@/_1pc/[rru]_{w\tn (r_{x,y}^{23} r_{y,z}^{34}+ r_{x,z}^{24}r_{y,z}^{34})} & &
}
}} \quad = \\[.4cm]
= \qquad \vcenter{\vbox{
\xymatrix{
 & & w\tn x\tn y\tn z \ar@/^1pc/[rrd]^{r_{w, x\tn y\tn z}} & & \\
w\tn x\tn y\tn z 
\ar@/_10pt/[rrd]_{r_{w,x\tn y\tn z}} 
\ar@/^10pt/[rru]^{w\tn (r_{y,z}^{34} r_{x,y}^{23}+ r_{y,z}^{34}r_{x,z}^{24})} & &
\Uparrow T_{(w,(r_{y,z}^{34} r_{x,y}^{23}+ r_{y,z}^{34}r_{x,z}^{24}))} \qquad\qquad & &
w\tn x\tn y\tn z \\
 & & w\tn x\tn y\tn z 
\ar@/_10pt/[rru]_{w\tn (r_{x,y}^{23} r_{y,z}^{34}+ r_{x,z}^{24}r_{y,z}^{34})} 
\ar@/^10pt/[rru]^{w\tn (r_{y,z}^{34} r_{x,y}^{23}+ r_{y,z}^{34}r_{x,z}^{24})}
\ar@{{}{ }{}} [rru]|{ \Uparrow w\tn P_{x,y,z}}
& &
} 
}} \, .
\end{multline}
With similar computations as before, this yields:
\begin{equation}
\label{nat2eq}
\left( r_{w,x}^{12}+r_{w,y}^{13}+r_{w,z}^{14}\right) \t \v{P_{x,y,z}^{234}} + r_{y,z}^{34} \t \big(\v{P_{w,x,y}^{123}} + \v{P_{w,y,z}^{124}} \big) - \left( r_{x,y}^{23}+r_{x,z}^{24}\right) \t \v{P_{w,y,z}^{134}} \, .
\end{equation}

Two analogous relations are derived from naturality conditions of $T_{(s'(Q_{x,y,z},w)}$ and $T_{(t'(Q_{x,y,z},w)}$ (resp. $T_{(w, s'(Q_{x,y,z})}$ and $T_{(w,t'(Q_{x,y,z})}$) with respect to the 2-morphism $Q_{x,y,z}$. Similarly to \eqref{nat1} we have:
\begin{equation}
\label{nat3}
\comp{T_{(s'(Q_{x,y,x}),w)}}{(Q_{x,y,z}\tn w)r_{x\tn y\tn z,w}} = \comp{r_{x\tn y\tn z,w}(Q_{x,y,z}\tn w)}{T_{(t'(Q_{x,y,z}),w)}}
\end{equation}  
which is displayed as:
\begin{multline}
\allowdisplaybreaks
\label{nat3diag}
\vcenter{\vbox{
\xymatrix{
 & & x\tn y\tn z\tn w \ar@/^1pc/[rrd]^{r_{x\tn y\tn z,w}} & & \\
x\tn y\tn z\tn w 
\ar@/_10pt/[rrd]_{r_{x\tn y\tn z,w}} 
\ar@/^10pt/[rru]^{(r_{x,y}^{12} r_{x,z}^{13} + r_{x,y}^{12} r_{y,z}^{23})\tn w} 
\ar@/_10pt/[rru]_{( r_{x,z}^{13} r_{x,y}^{12} + r_{y,z}^{23} r_{x,y}^{12})\tn w}
\ar@{{}{ }{}} [rru]|{ \Uparrow Q_{x,y,z} \tn w} & &
\qquad\qquad \Uparrow T_{((r_{x,z}^{13} r_{x,y}^{12} + r_{y,z}^{23} r_{x,y}^{12}),w)} & &
x\tn y\tn z\tn w \\
 & & x\tn y\tn z\tn w \ar@/_1pc/[rru]_{(r_{x,z}^{13} r_{x,y}^{12} + r_{y,z}^{23} r_{x,y}^{12}) \tn w)} & &
}
}} \quad = \\[.4cm]
= \qquad \vcenter{\vbox{
\xymatrix{
 & & x\tn y\tn z\tn w \ar@/^1pc/[rrd]^{r_{x\tn y\tn z,w}} & & \\
x\tn y\tn z\tn w 
\ar@/_10pt/[rrd]_{r_{x\tn y\tn z,w}} 
\ar@/^10pt/[rru]^{(r_{x,y}^{12} r_{x,z}^{13} + r_{x,y}^{12} r_{y,z}^{23})\tn w} & &
\Uparrow T_{(r_{x,y}^{12} r_{x,z}^{13} + r_{x,y}^{12} r_{y,z}^{23}),w)} \qquad\qquad & &
x\tn y\tn z\tn w \\
 & & x\tn y\tn z\tn w 
\ar@/_10pt/[rru]_{(r_{x,z}^{13} r_{x,y}^{12} + r_{y,z}^{23} r_{x,y}^{12}) \tn w)} 
\ar@/^10pt/[rru]^{(r_{x,y}^{12} r_{x,z}^{13} + r_{x,y}^{12} r_{y,z}^{23})\tn w}
\ar@{{}{ }{}} [rru]|{ \Uparrow Q_{x,y,z} \tn w}
& &
}
}} \, .
\end{multline}
Following the same line of computations as before, equating the arrow parts of both sides we get the relation:
\begin{equation}
\label{nat3eq}
(r^{14}_{x,w} + r^{24}_{y,w} + r^{34}_{z,w}) \t \v{Q^{123}_{x,y,z}} + r^{12} \t \big( \v{Q^{134}_{x,z,w}} + \v{Q^{234}_{y,z,w}} \big)
- (r^{13}_{x,z} + r^{23}_{y,z}) \t \v{Q^{124}_{x,y,w}} = 0 \, .
\end{equation}
Similarly to \eqref{nat2}, we also know that it holds:
\begin{equation}
\label{nat4}
\comp{ T_{(w,s'(Q_{x,y,z}))} }{ (w\tn Q_{x,y,z}) \, r_{w,x\tn y\tn z} } \, = \, 
\comp{ r_{w,x\tn y\tn z} \, (w\tn Q_{x,y,z}) }{ T_{(w,t'(Q_{x,y,z}))} } 
\end{equation} 
or, graphically:
\begin{multline}
\allowdisplaybreaks
\label{nat4diag}
\vcenter{\vbox{
\xymatrix{
 & & w\tn x\tn y\tn z \ar@/^1pc/[rrd]^{r_{w, x\tn y\tn z}} & & \\
w\tn x\tn y\tn z 
\ar@/_10pt/[rrd]_{r_{w,x\tn y\tn z}} 
\ar@/^10pt/[rru]^{w\tn ( r_{x,y}^{23} r_{x,z}^{24} + r_{x,y}^{23} r_{y,z}^{34} )} 
\ar@/_10pt/[rru]_{w\tn ( r_{x,z}^{24} r_{x,y}^{23} + r_{y,z}^{34} r_{x,y}^{23})}
\ar@{{}{ }{}} [rru]|{ \Uparrow w \tn Q_{x,y,z} } & &
\qquad\qquad \Uparrow T_{(w,( r_{x,z}^{24} r_{x,y}^{23} + r_{y,z}^{34} r_{x,y}^{23}) )} & &
w\tn x\tn y\tn z \\
 & & w\tn x\tn y\tn z \ar@/_1pc/[rru]_{w\tn ( r_{x,z}^{24} r_{x,y}^{23} + r_{y,z}^{34} r_{x,y}^{23}) } & &
}
}} \quad = \\[.4cm]
= \qquad \vcenter{\vbox{
\xymatrix{
 & & w\tn x\tn y\tn z \ar@/^1pc/[rrd]^{r_{w, x\tn y\tn z}} & & \\
w\tn x\tn y\tn z 
\ar@/_10pt/[rrd]_{r_{w,x\tn y\tn z}} 
\ar@/^10pt/[rru]^{w\tn ( r_{x,y}^{23} r_{x,z}^{24} + r_{x,y}^{23} r_{y,z}^{34} )} & &
\Uparrow T_{(w,( r_{x,y}^{23} r_{x,z}^{24} + r_{x,y}^{23} r_{y,z}^{34} ))} \qquad\qquad & &
w\tn x\tn y\tn z \\
 & & w\tn x\tn y\tn z 
\ar@/_10pt/[rru]_{w\tn ( r_{x,z}^{24} r_{x,y}^{23} + r_{y,z}^{34} r_{x,y}^{23})} 
\ar@/^10pt/[rru]^{w\tn ( r_{x,y}^{23} r_{x,z}^{24} + r_{x,y}^{23} r_{y,z}^{34} )}
\ar@{{}{ }{}} [rru]|{ \Uparrow w\tn Q_{x,y,z}}
& &
} 
}} \, .
\end{multline}
This eventually leads to the relation:
\begin{equation}
\label{nat4eq}
(r^{12}_{w,x} + r^{13}_{w,y} + r^{14}_{w,z}) \t \v{Q^{234}_{x,y,z}} + r_{x,y}^{23} \t \big( \v{P^{124}_{w,x,z}} + \v{P^{134}_{w,y,z}} \big) - (r^{24}_{x,z} + r^{34}_{y,z}) \t \v{P^{123}_{w,x,y}} = 0 \, .
\end{equation} 
Finally, a further condition is obtained by considering: 
$$ r_{x,y}:x\tn y\rightarrow x\tn y \, , \quad r_{z,w}:z\tn w \rightarrow z\tn w $$ 
and writing the interchange law \eqref{Tfg}, that assures $T_{(r^{12}_{x,y},r^{34}_{y,z})}$ is well defined. In our present setting it reads:
\begin{equation}
\label{int1} 
\comp{ T_{(r_{x,y}, \, z\tn w)} \, (x\tn y\tn r_{z,w}) }{ (r_{x,y}\tn z\tn w) \, T_{(x\tn y, \, r_{z,w})} } 
\, = \, 
\comp{ T_{(x\tn y, \, r_{z,w})} \, (r_{x,y} \tn z\tn w) }{ (x\tn y \tn r_{z,w}) \, T_{(r_{x,y}, \, z\tn w)} } 
\, .
\end{equation}
Graphically:
\begin{multline}
\label{int1diag}
\vcenter{\vbox{
\xymatrix @C=4em @R=3em{
 & & x\tn y\tn z\tn w \ar@/^1pc/[dr]^{r_{x\tn y,z\tn w}} & \\
x\tn y\tn z\tn w \ar@/_1pc/[dr]_{r_{x\tn y,z\tn w}} \ar[r]^{r^{12}_{x,y}} &
x\tn y\tn z\tn w \ar[r]^{r_{x\tn y,z\tn w}} \ar@/^1pc/[ur]^{r^{34}_{z,w}} &
x\tn y\tn z\tn w \ar[r]^{r^{34}_{z,w}} \ar@{{}{ }{}} [u]|{ \Uparrow T_{(x\tn y,r_{z,w})}} &
x\tn y\tn z\tn w \\
 & x\tn y\tn z\tn w \ar@{{}{ }{}} [u]|{ \Uparrow T_{(r_{x,y},z\tn w)}} \ar@/_1pc/[ur]_{r^{12}_{x,y}} & & 
}
}} \quad = \\ = \quad 
\vcenter{\vbox{
\xymatrix @C=4em @R=3em{
 & & x\tn y\tn z\tn w \ar@/^1pc/[dr]^{r_{x\tn y,z\tn w}} & \\
x\tn y\tn z\tn w \ar@/_1pc/[dr]_{r_{x\tn y,z\tn w}} \ar[r]^{r^{34}_{z,w}} &
x\tn y\tn z\tn w \ar[r]^{r_{x\tn y,z\tn w}} \ar@/^1pc/[ur]^{r^{12}_{x,y}} &
x\tn y\tn z\tn w \ar[r]^{r^{12}_{x,y}} \ar@{{}{ }{}} [u]|{ \Uparrow T_{(r_{x,y},z\tn w)}} &
x\tn y\tn z\tn w \\
 & x\tn y\tn z\tn w \ar@{{}{ }{}} [u]|{ \Uparrow T_{(x\tn y,r_{z,w})}} \ar@/_1pc/[ur]_{r^{34}_{z,w}} & & 
}
}}
\end{multline}
Writing \eqref{int1} along the arrow parts (the source parts of both sides are trivially equal) and using \eqref{Tlinob} {we get:}
\begin{equation}
\label{int1eq}
\begin{split}
& r^{12}_{x,y} \, \v{T_{(x\tn y,r_{z,w})}} + \v{T_{(r_{x,y},z\tn w)}} \, r^{34}_{z,w} - r^{34}_{y,z} \, \v{T_{(r_{x,y},z\tn w)}} - \v{T_{(x\tn y,r_{z,w})}} \, r^{12}_{x,y} = \\[.3cm]
= \; & r^{12}_{x,y} \, \big( \v{T^{134}_{(x,r_{z,w})}} + \v{T^{234}_{(y,r_{z,w})}} \big) + \big(\v{T^{123}_{(r_{x,y},z)}} + \v{T^{124}_{(r_{x,y},w)}} \big) \, r^{34}_{z,w} - r^{34}_{z,w} \, \big(\v{T^{123}_{(r_{x,y},z)}} + \v{T^{124}_{(r_{x,y},w)}} \big) - \big( \v{T^{134}_{(x,r_{z,w})}} + \v{T^{234}_{(y,r_{z,w})}} \big) \, r^{12}_{x,y} = \\[.3cm]
= \; & r_{x,y}^{12} \t \big ( \v{P_{x,z,w}^{134}}+\v{P_{y,z,w}^{234}}\big) - r_{z,w}^{34} \t \big(\v{Q_{x,y,z}^{123}}+\v{Q_{x,y,w}^{124}}\big) = 0 \, . 
\end{split}
\end{equation}
We have finished proving our first main theorem:
\begin{Theorem}\label{main1}
 Consider an infinitesimal 2-braiding $(r,T)$ in a symmetric strict monoidal linear 2-category  $\cC=(\cC,\tn,I, B)$. Let $P_{x,y,z}$ and $Q_{x,y.z}$ be as defined in \eqref{defP} and \eqref{defQ}. The equations below hold:
\begin{align}
 \label{nat1eqt} 
(r_{x,w}^{14}+r_{y,w}^{24}+r_{z,w}^{34}) \t \v{ P_{x,y,z}^{123}} -
(r_{x,y}^{12} + r_{x,y}^{13})\t \v{Q_{y,z,w}^{234}}
+ r_{y,z}^{23} \t (\v{Q_{x,y,w}^{124}}+ \v{Q_{x,z,w}^{134}})&=0 \,, \\
\label{nat2eqt}
\left( r_{w,x}^{12}+r_{w,y}^{13}+r_{w,z}^{14}\right) \t \v{P_{x,y,z}^{234}} + r_{y,z}^{34} \t \big(\v{P_{w,x,y}^{123}} + \v{P_{w,y,z}^{124}} \big) - \left( r_{x,y}^{23}+r_{x,z}^{24}\right) \t \v{P_{w,y,z}^{134}} &=0\, ,\\
\label{nat3eqt}
(r^{14}_{x,w} + r^{24}_{y,w} + r^{34}_{z,w}) \t \v{Q^{123}_{x,y,z}} + r^{12} \t \big( \v{Q^{134}_{x,z,w}} + \v{Q^{234}_{y,z,w}} \big)
- (r^{13}_{x,z} + r^{23}_{y,z}) \t \v{Q^{124}_{x,y,w}} &= 0 \, ,\\
\label{nat4eqt}
(r^{12}_{w,x} + r^{13}_{w,y} + r^{14}_{w,z}) \t \v{Q^{234}_{x,y,z}} + r_{x,y}^{23} \t \big( \v{P^{124}_{w,x,z}} + \v{P^{134}_{w,y,z}} \big) - (r^{24}_{x,z} + r^{34}_{y,z}) \t \v{P^{123}_{w,x,y}} &= 0 \, ,\\
\label{nat5eqt}
 r_{x,y}^{12} \t \big ( \v{P_{x,z,w}^{134}}+\v{P_{y,z,w}^{234}}\big) - r_{z,w}^{34} \t \big(\v{Q_{x,y,z}^{123}}+\v{Q_{x,y,w}^{124}}\big) &= 0\,.
\end{align}

\end{Theorem}

Continuing the notation of the previous Theorem, suppose that  the infinitesimal 2-braiding $(r,T)$ is coherent, Definition \eqref{coherent}. Looking at equation \eqref{nat5eqt}, and applying $B_{y,z}$ to both sides of it,  we reach that for any objects $x,y,z,w$ of $\Cc$:
$$r_{x,y}^{13} \t \big ( \v{P_{x,z,w}^{124}}+\v{P_{y,z,w}^{324}}\big) - r_{z,w}^{24} \t \big(\v{Q_{x,y,z}^{132}}+\v{Q_{x,y,w}^{134}}\big) = 0.$$
If $(r,T)$ is coherent then:
$$r_{x,y}^{13} \t \big ( \v{P_{x,z,w}^{124}}-\v{P_{z,y,w}^{234}}-\v{Q_{z,y,w}^{234}}\big) - r_{z,w}^{24} \t \big(\v{-Q_{x,z,y}^{123}}-\v{P_{x,z,y}^{123}}+\v{Q_{x,y,w}^{134}}\big)  = 0.$$               
Therefore:
\begin{Theorem}\label{main2}
In the conditions of Theorem \ref{main1}, if $(r,T)$ is coherent then furthermore we have:
\begin{equation}\label{nat6eqt}
 r_{x,y}^{13} \t \big ( \v{P_{x,z,w}^{124}}-\v{P_{z,y,w}^{234}}-\v{Q_{z,y,w}^{234}}\big) + r_{z,w}^{24} \t \big(\v{Q_{x,z,y}^{123}}+\v{P_{x,z,y}^{123}}-\v{Q_{x,y,w}^{134}}\big)  = 0
\end{equation}

\end{Theorem}

\noindent{The relations appearing in Theorems \ref{main1} and \ref{main2} are our proposal for a categorification of the 4-term relations.}
\section{A related construction: the categorified Knizhnik-Zamolodchikov connection}\label{2KZ}

We now re-interpret the construction of the categorified Knizhnik-Zamolodchikov (KZ in the following) connection \cite{CFM1,CFM2}, emphasizing the close relation with strict infinitesimal 2-braidings. Let us fix a positive integer $n\in\mathbb{N}$.
Let $\C(n)$ be the configuration space of $n$ distinct particles in the complex plane:
$$\C(n)=\{(z_1,\dots,z_n) \in \C^n\colon z_i \neq z_j \textrm{ if } i \neq j\}.$$
This has an obvious left action $S_n \ni \sigma \mapsto L_\s \in {\rm diff}(\C(n)) $ of the symmetric group $S_n$. The configuration space of $n$ indistinguishable particles in $\C$ is defined as $\C(n)/S_n$.
Define closed 1-forms $\w_{ij}$ in the configuration space $\C(n)$, for $1\leq i,j \leq n$ and $i \neq j$:
$$\w_{ij}=\frac{dz_i - dz_j}{z_i-z_j} \, . $$
Clearly for each $\s\in S_n$: 
\begin{equation}\label{cov3}
L_\s^*(\w_{ij})=\w_{\s^{-1}{(i)}\s^{-1}(j)}.
\end{equation}
We recall the well known Arnold's relation \cite{Ar}, for each distinct indices $i,j,k \in \{1,\dots,n\}$:
\begin{equation}\label{ai}
\w_{ij} \wedge \w_{jk}+\w_{jk} \wedge \w_{ki}+\w_{ki} \wedge \w_{ij}=0 .
\end{equation}
Note that for each distinct indices $i,j \in \{1,\dots,n\}$, we have $\w_{ij}=\w_{ji}$.

Consider an infinitesimal braiding $r$ in a linear monoidal category $\Cc$ (see \cite{K,C} and the Introduction). Choose an object $x \in \Cc$. Given a positive integer $n$ and $1\leq a < b \leq n$ we thus have a morphism $r^{ab} \colon x^{\tn n} \to x^{\tn n}$. Recall the construction of the differential crossed module $\GL(y)$, associated to $\Cc$ and and an object $y$ of $\Cc$, outlined in Section \ref{diffx}. The KZ-connection associated to the triple $(\Cc,r,x)$ is:
$$A=\sum_{1\leq a < b \leq n} \w_{ab}r^{ab}.$$
(The underlying differential crossed module is therefore $\GL(x^{\otimes n})$, in this case collapsing to a Lie algebra).
{Given that $d A=0$, the curvature ${\cal F}_A= dA+A \wedge A$} of $A$ is:
\begin{align*}
{\cal F}_A=A \wedge A= &\quad 2\sum_{a<b<d}\w_{ab}\wedge \w_{ad} \,[r^{ab},r^{ad}]+2\sum_{a<b<d}\w_{ad}\wedge \w_{bd}\, [r^{ad},r^{bd}]+2\sum_{a<b<d}\w_{ab}\wedge \w_{bd}\, [r^{ab},r^{bd}]\,.
\end{align*}
The last term is (by using Arnold relation \eqref{ai}):
                $$ -2\sum_{a<b<d}\w_{bd}\wedge \w_{da}\, [r^{ab},r^{bd}] -2\sum_{a<b<d}\w_{da}\wedge \w_{ab}\, [r^{ab},r^{bd}].$$
By using the fact that $\w_{ab}=\w_{ba}$ we conclude:
$${\cal F}_A= -2\sum_{a<b<d}\w_{bd}\wedge \w_{da}\, [r^{ab}+r^{ad},r^{bd}] -2\sum_{a<b<d}\w_{da}\wedge \w_{ab}\, [r^{ab} ,r^{bd}+r^{ad}].$$
We can see that the curvature of $A$ vanishes exactly because $r$, being an infinitesimal braiding operator, does satisfy the 4-term relations. 

Let us now suppose that what we have is an infinitesimal 2-braiding $(r,T)$ in a linear monoidal symmetric 2-category $\Cc$. The previous calculation of the curvature of $A=\sum_{1\leq a < b \leq n} \w_{ab}r^{ab}$ remains the same. Recalling \eqref{defP} and \eqref{defQ} we set:
$$P^{abc}=\v{P_{x,x,x}^{abc}} \qquad \an \qquad Q^{abc}=\v{Q_{x,x,x}^{abc}} \, .$$ 
Since the differential $\d$ in $\gl(x)$ is the target map, we have:
$$ \d(P^{abc})=[r^{bc},r^{ab}+r^{ac}] \qquad 
\an \qquad \d(Q^{abc})=[r^{ab},r^{ac}+r^{bc}] \, .$$
Therefore we can write:
\begin{equation}
{\cal F}_A= 2\sum_{a<b<c}\w_{bc}\wedge \w_{ca}\,\d(P^{abc}) -2\sum_{a<b<c}\w_{ca}\wedge \w_{ab}\,\d(Q^{abc}).
\end{equation}
{In order that a 2-connection $(A,B)$ has a well defined two-dimensional holonomy its fake curvature  $\partial (B) -{\cal F}_A$ tensor must be zero. We thus set:} 
\begin{equation}
\label{defB}
B= 2\sum_{a<b<c}\w_{bc}\wedge \w_{ca}\,P^{abc} -2\sum_{a<b<c}\w_{ca}\wedge \w_{ab}\,Q^{abc} \, .
\end{equation}

Now we condense our notation: we denote $\alpha \wedge \beta=\alpha\beta$ and we put:
$$ [ab]=\w_{ab} \, , \qquad  [abc]=\w_{ab}\wedge \w_{bc}=[ab][bc] \, , \qquad  W_{ab,cde}=r^{ab} \t P^{cde} \, , \qquad  Z_{ab,cde}=-r^{ab} \t Q^{cde} \, .$$
Note that $r^{ab} \t P^{cde}=0$ and  $r^{ab} \t Q^{cde}=0$ if $\{a,b\} \cap \{c,d,e\}=\emptyset$. 
In this notation we can simply write:
$$B= 2\sum_{a<b<c}[bca]\,P^{abc} -2\sum_{a<b<c}[cab]\,Q^{abc},$$
and we compute the 2-curvature $\G_{(A,B)}=dB+A\wedge B$, of {$(A,B)$}, to be:
\begin{align*}
&\G_{(A,B)}=\\
& \quad \sum_{a<b<c<d} [cd][bca] W_{cd,abc}   +[cd][cab] Z_{cd,abc} +  [bd][bca] W_{bd,abc} + [bd][cab] Z_{bd,abc} +  [ad][bca] W_{ad,abc} + [ad][cab] Z_{ad,abc} \\
& +\sum_{a<b<c<d} [cd][bda] W_{cd,abd}   +[cd][dab] Z_{cd,abd} +  [cb][bda] W_{cb,abd} + [cb][dab] Z_{cb,abd} +  [ca][bda] W_{ca,abd} + [ca][dab] Z_{ca,abd} \\
& +\sum_{a<b<c<d} [bd][cda] W_{bd,acd}   +[bd][dac] Z_{bd,acd} +  [bc][cda] W_{bc,acd} + [bc][dac] Z_{bc,acd} +  [ab][cda] W_{ab,acd} + [ab][dac] Z_{ab,acd} \\
& +\sum_{a<b<c<d} [ab][cdb] W_{ab,bcd}   +[ab][dbc] Z_{ab,bcd} +  [ac][cdb] W_{ac,bcd} + [ac][dbc] Z_{ac,bcd} +  [ad][cdb] W_{ad,bcd} + [ad][dbc] Z_{ad,bcd} 
\end{align*}
Consider now the following closed differential forms: 
$$ \{ \w_{i_aj_a}\wedge \w_{i_bj_b}\wedge \w_{i_cj_c} \; \mbox{ s.t. } \; i_k<j_k \; \an \; j_k< j_{k'} \; \mbox{ for } \; k<k' \}, $$
with all indices running in $\{1,\ldots ,n\}$. These are linearly independent, which follows from the main result of  \cite{Ar}, and can easily be proven directly. In fact these forms are known to span a subspace of the space $\Omega^3(\C(n))$ of differential 3-forms in the configuration space $\C(n)$, generating the Rham cohomology in degree three.
Using our notation, these six linearly independent closed forms in  $\Omega^3(\C(n))$ are:
$$ \left\{ \, [ab][ac][ad] \, , \; [ab]bc][ad] \, ,\; [ab][ac][bd] \, ,\; [ab][bc][bd] \, , \; [ab][ac][cd] \, , \; [ab][bc][cd] \, \right\} \qquad (1\leq a <b<c<d\leq n) \, .$$
In terms of this basis, each of the 3-forms appearing in the 24 terms expressing $\G_{(A,B)}$ is written as a vector with six components. We gather these in the matrix $M$ below, in order of appearance: 
\scriptsize
$$M= \left(
\begin{array}{cccccccccccccccccccccccc}
 0 &  0 &  0 &  0 &  1&  -1 & -1&   1 &  0 &  0&  -1&   1&   1&  -1 &  1 & -1 &  1&  -1 &  0  &0 &  0&   0&  -1&  0\\
  0 &  0 &  0 &  0&  -1 &  0 &  0 &  0&  -1&   1 &  0  & 0 &  0 &  0&  -1 &  1&   0 &  0&   0  & 0 &  0 &  0&   0  & 1\\
  0 &  0  & 1 & -1  & 0&   0&   0 &  0 &  0 &  0 &  1 &  0 &  0 &  1&   0&    0 &  0 &  0 &  0&   0 & -1 &  1 &  0 &  0\\
  0  & 0 & -1&   0&   0&   0&   1&   0 &  1 &  0 &  0&   0&  -1&   0&   0 &  0 &  0 &  0 &  1&  -1 &  1 & -1 &  1 & -1\\
  1 & -1&   0&   0 &  0  & 0 &  1&  -1 &  0  & 0 &  0 &  0 & -1&   0&  -1 &  0&  -1 &  0 &  0&   0 &  1 &  0 &  1  & 0 \\
  -1 &  0 &  0&   0 &  0&   0&  -1 &  0&   0 &  0 &  0 &  0&   1&   0&   1&   0 &  0 &  0 & -1  & 0 &  -1 &  0 & -1&   0
\end{array}
\right)
$$
\normalsize
This is convenient to describe the flatness of the KZ 2-connection in terms of linear algebra.
Indeed consider the vector $V$ whose 24 components are the coefficients along the 3-forms in the previous expression of the 2-curvature $\G_{(A,B)}$:   
\begin{multline} 
V= \left( W_{cd,abc} , Z_{cd,abc},W_{bd,abc},Z_{bd,abc},W_{ad,abc}, Z_{ad,abc}, W_{cd,abd},Z_{cd,abd},W_{cb,abd}, Z_{cb,abd}, W_{ca,abd}. Z_{ca,abd},W_{bd,acd}, \right. \\ \left. Z_{bd,acd}, W_{bc,acd}, Z_{bc,acd},W_{ab,acd}, Z_{ab,acd},W_{ab,bcd},Z_{ab,bcd}, W_{ac,bcd},Z_{ac,bcd}, W_{ad,bcd}, Z_{ad,bcd} \right)
\end{multline}
Then the vanishing of the 2-curvature $\G_{(A,B)}$ is equivalent to the matrix equation $MV=0$. Given that the rank of $M$ is six, this leads to six independent equations in the variables $Z$ and $W$. These 2-flatness conditions have previously appeared in \cite{CFM1,CFM2}, only with a different choice of generators. 

What we want to point out here is that these six 2-flatness conditions are equivalent to the ones appearing in Theorems \ref{main1} and \ref{main2}. To see this we left-multiply $M$  by the matrix $N$ (with rank six): \footnotesize
$$N=\left(
\begin{array}{rrrrrr}
0 & -1 & 0 & -1 & 0 & -1 \\
0 & 0 & 0 & 0 & 0 & -1 \\
1 & 1 & 1 & 1 & 1 & 1 \\
0 & 0 & 0 & 1 & 0 & 1 \\
0 & 0 & 0 & 0 & -1 & -1 \\
0 & 0 & 1 & 0 & 0 & 0
\end{array}
\right) \, . $$
\normalsize
Therefore:
\scriptsize
$$NM=         	
\left(
\begin{array}{rrrrrrrrrrrrrrrrrrrrrrrr}
1 & 0 & 1 & 0 & 1 & 0 & 0 & 0 & 0 &
-1 & 0 & 0 & 0 & 0 & 0 & -1 & 0 & 0
& 0 & 1 & 0 & 1 & 0 & 0 \\
1 & 0 & 0 & 0 & 0 & 0 & 1 & 0 & 0 &
0 & 0 & 0 & -1 & 0 & -1 & 0 & 0 & 0
& 1 & 0 & 1 & 0 & 1 & 0 \\
0 & -1 & 0 & -1 & 0 & -1 & 0 & 0 & 0
& 1 & 0 & 1 & 0 & 0 & 0 & 0 & 0 & -1
& 0 & -1 & 0 & 0 & 0 & 0 \\
-1 & 0 & -1 & 0 & 0 & 0 & 0 & 0 & 1
& 0 & 0 & 0 & 0 & 0 & 1 & 0 & 0 & 0
& 0 & -1 & 0 & -1 & 0 & -1 \\
0 & 1 & 0 & 0 & 0 & 0 & 0 & 1 & 0 &
0 & 0 & 0 & 0 & 0 & 0 & 0 & 1 & 0 &
1 & 0 & 0 & 0 & 0 & 0 \\
0 & 0 & 1 & -1 & 0 & 0 & 0 & 0 & 0 &
0 & 1 & 0 & 0 & 1 & 0 & 0 & 0 & 0 &
0 & 0 & -1 & 1 & 0 & 0
\end{array}
\right) \, .
$$
\normalsize
Then the six equations coming from $NMV=0$ are the equations appearing in Theorems \ref{main1} and \ref{main2}. We have proved the following result, which provides a new interpretation of the flatness conditions of the KZ 2-connection.
\begin{Theorem}\label{main3}
Consider a coherent infinitesimal 2-braiding $(r,T)$ in a linear symmetric strict monoidal 2-category $\Cc$. Choose and object $x$ of $\Cc$. Then the pair $(A,B)$, where:
$$A=\sum_{1\leq a < b \leq n} \w_{ab}r^{ab} \, , \qquad 
B= 2\sum_{a<b<d}\w_{bd}\wedge \w_{da}\,P^{abd} -2\sum_{a<b<d}\w_{da}\wedge \w_{ab}\,Q^{abd}$$
($\w_{ij},P^{ijk},Q^{ijk}$ defined as before) is a flat (and fake flat) 2-connection in the configuration space $\C(n)$, with values in the differential crossed module $\GL(x^{\otimes n})$, associated to the object $x^{n \otimes}$ of $\Cc$.
\end{Theorem}

By using exactly the same techniques as in \cite{CFM1,CFM2} we can prove that:
\begin{Theorem}
 In the conditions of the previous theorem, if the infinitesimal  2-braiding is totally symmetric, then the 2-connection $(A,B)$ is invariant under the obvious action of the symmetric group. 
\end{Theorem}
\section{Infinitesimal 2-braidings for a differential crossed module}\label{ibdcm}

\subsection{Preliminaries}

\subsubsection{Summary of notation and algebra actions}
Given a Lie algebra $\lg$, let $\U(\lg)$ denote its universal enveloping algebra. Let  $n$ be a positive integer. Recall that we have {a Hopf algebra  isomorphism}
{$\U(\lg^{\oplus n}) \to \U(\lg)^{\otimes n} ,$}
being (for $X_i \in \lg$, where $i \in \{1,\dots,n\}$):
\begin{equation}\label{morph}
(X_1,X_2,\dots,X_n) \mapsto X_1 \otimes 1 \otimes \dots \otimes 1 + 1 \otimes X_2 \otimes 1 \otimes \dots \otimes 1 + \dots + 1 \otimes \dots \otimes 1 \otimes X_n. 
\end{equation}
Let $\lG=(\d\colon \lh \to \lg,\t)$ be a differential crossed module, which will be fixed throughout this subsection. Consider the semidirect product $\le= \lh \rtimes_{\t} \lg$. This is a Lie algebra with Lie bracket:
$$[(u,X),(v,Y)]=\big([u,v]+X \t v-Y \t u,[X,Y]\big) \, , \qquad \forall \, X,Y \in \lg \; \an \; u,v \in \lh \, . $$
Note that $\lh$ and $\lg$ (more precisely, their images under the maps $\lh\ni v\mapsto (v,0)\in \le$ and $\lg\ni X \mapsto (0,X) \in \le$) are Lie subalgebras of $\le$. We have a Lie algebra map: $$\beta: \le\ni (v,X) \mapsto X +\d(v) \in \lg \, .$$
By restriction of the adjoint action of $\le$ on itself we recover the action of $\lg$ on $\lh$ and actions (by Lie algebra derivations) of $\lg$ and of $\lh$ on $\le$. All these actions canonically lift to the corresponding enveloping algebras, and to tensor products of them. As an example, since there is an action $\t$ of $\lg$ on $\lh$ by Lie algebra derivations, then $\lg$ acts on $\U(\lh)$ by algebra derivations, explicitly as:
$$X \t (v_1v_2\dots v_n) =\big((X \t v_1)  v_2 \dots v_n\big)+\big(v_1(X\t  v_2)\dots v_n\big)+\dots +\big(v_1  v_2 \dots {(X\t v_n)} \big) \, .$$
The Lie algebra maps $\partial\colon \lh \to \lg$ and $\beta \colon \le \to \lg$ intertwine the $\lg$-actions. We then have algebra maps $\partial\colon \U(\lh) \to \U(\lg)$ and $\beta \colon \U(\le) \to \U(\lg)$, which also intertwine the $\lg$-actions. 

By considering the product action (which is an action by derivations) of $\lg^{\oplus n}$ on $\lh^{\oplus n}$, namely:
$$(X_1,\dots,X_n) \t (v_1,\dots,v_n)=(X_1 \t v_1,\dots,X_n \t v_n) \, ,$$
and the product map $\d\colon \lh^{\oplus n} \to \lg^{\oplus n}$, so that $\d(v_1,\dots,v_n)=(\d(v_1),\dots, \d(v_n))$, we define a differential crossed module $\lG^{\oplus n}$. Therefore the previous description applies to it, and it will be used without much comment.  In particular we have maps $\partial\colon \U(\lh)^{\otimes n} \to \U(\lg)^{\otimes n}$ and $\beta \colon \U(\le)^{\otimes n} \to \U(\lg)^{\otimes n}$, and we have a Lie algebra action by algebra derivations  of $\lg^{\oplus n}$ on $\U(\le^{\oplus n})\cong \U(\le)^{\otimes n}$.

Recall that we have a Lie algebra morphism $\Delta\colon \lg \to \lg^{\oplus n}$, being $\Delta(x)=(x,\dots,x)$. 
This diagonal map, together with similar ones induced on enveloping algebras, permit us to define actions, by algebra derivations, of $\lg$ and of $\lh$ on $\U(\le^{\oplus n})\cong \U(\le)^{\otimes n}$, which will have a prime importance later. 

\subsubsection{An auxiliary vector space}

\noindent
From now on we will make no distinction between $\U(\la^{\oplus n})$ and $\U(\la)^{\otimes n}$, where $\la$ is a Lie algebra. 
Let $\lG=(\d\colon \lh \to \lg,\t)$ be a differential crossed module. As before we put $\le= \lh \rtimes_{\t} \lg$. By the Poincar\'{e}-Birkhoff-Witt Theorem both $\U(\lg)$ and $\U(\lh)$ are subalgebras of $\U(\le)$.
Inside the algebra $\U(\le^{\oplus n})$, we let $A_n$ be the smallest vector space containing all elements of the form $$x\, v \,y, \we x,y \in \U(\lg^{\oplus n}) \an v \in \lh^{\oplus n}.$$ By restriction of the action of $\lg^{\oplus n}$ on $\U(\le^{\oplus n})$ by derivations, $A_n$ is a $\lg^{\oplus n}$-module.

\begin{Definition} Let $n$ be a positive integer. The vector space $\Un$ is given by  the subspace $A_n$, modulo the relations:
$$x\,\d(u)\,y\,v\,z=x\, u \,y \,\d(v)\, z, $$
where $u,v \in \lh^{\oplus n} \subset \U(\le^{\oplus n})$ and $x,y,z \in \U(\lg^{\oplus n})$.
\end{Definition}

\noindent
We provide some examples since the notation may be misleading. If $n=1$ we have, for each $u,v \in \lh$,
$$\d(u) \,v=u \,\d(v) \, .$$ 
With $n=2$ we have, for each $u,v \in \lh$,
$$(\d(u),0)\,\,(0,v)=(u,0)\,\, (0,\d(v))$$
or, by passing from $\U(\le \oplus \le)$ to the isomorphic $\U(\le) \tn\U(\le)$ (see the morphism \eqref{morph}):
$$(\d(u)\tn 1)\,\, (1 \tn v)=(u \tn 1)\,\, (1 \tn \d(v))\, , $$
which is simply written as the following relation, {known to hold, in degree zero, up to homotopy,  in the tensor product of two chain complexes:}
$$\d(u) \tn v=u \tn \d(v) \, .$$
This adds of course to the relation obtained from:
$$\d\big ((u,0)\big)\,\,(v,0)=(u,0)\,\, \d\big((v,0)\big) \, ,$$
which is:
$$(\d(u) \,v)\tn 1=(u \,\d(v)) \tn 1 \, .$$
Similarly when $n=3$ we have, for each $u,v \in \lh$:
$$\d(u,0,0) \,\, (0,0,v)=(u,0,0)\,\,\d(0,0,v) $$
from which we obtain:
$$\d(u) \tn 1\tn v = u\tn 1 \tn \d(v) \, .$$
Therefore if $a\in \U(\lg)$ and $u,v \in \lh$:
$$\d(u) \tn a\tn v = u\tn a \tn \d(v) \, .$$
Of course different relative positions of $u$ and $v$ will give analogous relations. 

{Clearly the} algebra map $\beta \colon \U(\le^{\oplus n}) \to  \U(\lg^{\oplus n})$ descends to a vector space map: $${\beta\colon \Un \to \U(\lg^{\oplus n})}.$$ 

As mentioned before, we  have an action of $\lg$ on $\U(\le)^{\tn n}$ by algebra derivations. This is defined from the obvious action of  $\lg^{\oplus n}$ on $\U(\le)^{\tn n}$ by algebra derivations, and the diagonal map $\Delta \colon \lg \to \lg^{\oplus n}$.
\begin{Lemma}
The action of $\lg^{\oplus n}$ on $\U(\le^{\oplus n})$ descends to an action of $\lg^{\oplus n}$ on the vector space $\Un$.
\end{Lemma}
\begin{Proof}
Given $X \in \lg^{\oplus n}$, $u,v \in \lh^{n \otimes}$ and $x,y,z \in \U(\lg^{n \otimes})$ we have:
\begin{align*}
X \t (x\,u\,y\, \d(v)\,z)&=(X \t x)\,u \, y \, \d(v) \, z+x\, (X \t u)\, y \, \d(v)\, z+x\, u\, y \, (X \t \d(v))\, z+x\, u\, y\,  \d(v)\, (X \t z)\\
&=(X \t x)\, u\, y \, \d(v)\, z+x\, (X \t u)\, y \, \d(v)\, z+x\, u\, y\,  (\d(X \t v))\, z+x\, u\, y \d(v)\, (X \t z)
\end{align*}
\begin{align*}
X \t (x\, \d(u)\, y\,  v\, z)&=(X \t x)\, \d(u)\, y\,  v z+x(X \t \d(u))\, y \, v\, z+x\, \d(u)\, y\,  (X \t v)\, z+x\, \d(u)\, y \, v\, (X \t z)\\&=(X \t x)\, \d(u)\, y\,  v\,  z+x\, (\d (X \t u))\, y \, v\, z+x\, \d(u)\, y \, (X \t v)\, z+x\, \d(u)\, y \, v\, (X \t z).
\end{align*}
Clearly the difference between the previous two expressions is in the kernel of the projection $A_n \to \Un$. \qed\end{Proof}

\noindent
As a consequence of the previous Lemma, $\lg$ acts on $\Un$.

 Note that the maps:
$$(x,a) \in \U(\lg^{\oplus n}) \times A_n \mapsto xa \in A_n \, , \qquad 
(x,a) \in \U(\lg^{\oplus n}) \times A_n \mapsto ax \in A_n $$
clearly directly descend to maps: 
$$\U(\lg^{\oplus n}) \times \Un \to \Un. $$
These will be frequently used. Now the crucial property of all the construction {\cite{CFM1,CFM2}:}
\begin{Lemma}[The defining relations of $\Un$]
 For each $a,b \in \Un$ we have:
\begin{equation} 
\label{dr}
\beta(a)\, b=a \, \beta(b) \, .
\end{equation}
\end{Lemma}
\begin{Proof}
By bilinearity, it suffices proving this result on vector space generators.
Given $x,y,x',y' \in \U(\lg^{\oplus n})$ and $u,u'\in \lh$ we have:
\begin{equation*}
\beta(x\,u\,y)\, (x'\, u'\,y')=x\,\d(u)\,y\,x'\,v'\,y'= x\,u\,y\,x'\,\d(v')\,y'=(x\,u\,y)\,\beta(x'\,u'\, y') \,. 
\end{equation*}
\qed\end{Proof}

Let $m$ and $n$ be positive integers. There exist obvious inclusions $i_m$ and $i_n$ of $\lG^{\oplus m}$ and $\lG^{\oplus n}$ inside  $\lG^{\oplus (m+n)}$, being $i_m(X_1,\dots,X_m)=(X_1,\dots, X_m,0,\dots,0)$ and  $i_n(X_1,\dots,X_n)=(0,\dots,0, X_1,\dots, X_n),$ where $X_i \in \lg$. These extend to the associated enveloping algebras. Clearly $i_m(A_m) \subset A_{m+n}, $ and also this defines an inclusion map $i_m\colon \mathcal{U}^{(m)} \to \mathcal{U}^{(m+n)}$, commuting with $\beta$. Given $a \in \mathcal{U}^{(m)}$ and $l \in \U(\lg)^{\tn n}$ we define 
$$a \tn l=i_m(a) i_n(l)=  i_n(l) i_m(a) \, . $$
As a consequence of the defining relations \eqref{dr} we have, for $a \in \U^{(m)}$ and  $b \in \U^{(n)}$:
\begin{equation}
\label{dr2}
a\tn \beta(b)=\beta(a) \tn b \, . 
\end{equation}
We finish by introducing a map $\lh \otimes \U(\lg^{\tn n}) \to \Un$. By definition:
\begin{equation}
\label{hact}
(v,r) \mapsto v \t r=\Delta(v)r-r\Delta(v) \, .
\end{equation}
Note that given $v\in\lh$ and $r,r' \in \U(\lg^{\oplus n})$:
$$v \t (rr')=\Delta(v)rr'-rr'\Delta(v)=\Delta(v)rr'-r\Delta(v)r'+r\Delta(v)r'-rr'\Delta(v)=(v \t r)r'+r(v \t r') \, .$$

\subsection{{A strict 2-category $\Cc_{\lG}$ derived from  a differential crossed module $\lG$}}\label{s2c}

In \cite{CFM1,CFM2} we lifted the flatness conditions of the KZ 2-connection from the differential crossed module where the KZ connection takes value to an auxiliary short complex of vector spaces $\widehat{\d}:\Un\rightarrow\lg^{\tn n}$ (see e.g \cite[Section 4.3]{CFM1}. This auxiliary complex is where the universal relations expressing the flatness conditions live,  prior to being represented on $\GL(\V)$ for some complex of vector spaces $\V$. This construction mimics what happens with the 4-term relations, whose universal form lives in $\lg^{\tn n}$ and is then represented as KZ flatness conditions in $\gl(V^{\tn n})$ for some $\g$-module $V$.

In this section we aim at giving a categorical interpretation of the auxiliary complex involved in the KZ 2-flatness conditions. We associate to a differential crossed module $\lG$ a linear 2-category $\Cc_{\lG}$ where the universal structure responsible for the KZ 2-flatness, namely a {symmetric} quasi-invariant tensor in $\lG$, defines an infinitesimal 2-braiding.

\subsubsection{Objects and 1-morphisms}\label{o1m}
Let $\lG=(\d\colon \lh \to \lg,\t)$ be a differential crossed module. Given a positive integer $n$, we have a natural action of the permutation group $S_n$ on $\U(\le^{\oplus n})$. This clearly descends to an action on the vector space $\Un$. 

We now define a linear 2-category $\CG$ associated to $\lG$. The objects are given by non-negative integers $n \in \N$.
The set of $1$-morphisms is only non-empty in the case $n \to n$. In this case a 1-morphism is given by:
\begin{enumerate}[(a)]
 \item an element $R \in \U(\lg)^{\tn n}$;
 \item a linear map $\zeta \colon \lg \to \Un$;
\item a permutation $\s\in S_n$;
\item {an integer $k\in \mathbb{Z}$. All of the construction would make sense considering $k\in\C$. We however prefer to think of it as an integer, since it resembles a degree shift in a chain complex, an issue which is very common while categorifying  knot and link invariants.} 
\end{enumerate}
It is depicted as $ n \ra{(R,\zeta,\s,k)} n $. This data is subject to the following conditions: 
\begin{enumerate}[(i)]
\item $ X\t R =\beta(\zeta(X)) \, , \quad \forall \, X \in \lg$;
\item $u\t R=\zeta(\d(u)) \, , \quad \forall \, u \in \lh \;$ (here $\t$ is the map in \eqref{hact});
\item $ \zeta([X,Y])=X \t \zeta(Y)-Y \t \zeta(X) \, , \quad \forall \, X,Y,Z \in \lg$.
\end{enumerate}

\begin{Remark}[Decorated permutations]
To simplify the notation we will sometimes abbreviate $S_n \times \mathbb{Z}\ni (\s,k)=\tau$, and call such a pair a `decorated permutation'. Given any representation $i \mapsto \s(i)$ of $S_n$, $i=(1,\ldots ,n)$, we write $\tau(i)=\s(i)$ if $\tau=(\s,k)$. {The product of two decorated permutations $\tau=(\s,k)$ and $\tau'=(\s',k')$ is $\tau\cdot\tau':=(\s\s',k+k')$.}
\end{Remark}

Note that if $(R,\zeta,\s,k)$ and $(R',\zeta',\s,k)$ are morphisms $n \to n$ then so is $(R+R', \zeta+\zeta',\s,k)$. Therefore the set of morphisms $n \to n$ decomposes as a direct sum  
$$\hom(n,n) = \bigoplus_{(\s,k) \in S_n\times\mathbb{Z}} \hom_{\s,k}(n,n)$$
where $\hom_{\s,k}(n,n)$ is the vector space of morphisms $n \to n$ whose underlying decorated permutation is $(\s,k)$. The composition of 1-morphisms is 
$$n \ra{(R,\zeta,\s,k)} n   \ra{(R',\zeta',\s',k')}n \, = \,  n \ra{(R\sigma(R'),R \s(\zeta')+\zeta \s(R'),\s\s',k+k')} n \, .$$ 
This is associative and has  units. We now check that
conditions (i)-(iii) of the definition of a 1-morphism are satisfied by the composition of two 1-morphisms. For (i) and (ii) we use the fact that on $\U(\lg)^{\tn n}$ the $\lg$ and $\lh$ actions commute with the $S_n$ action. {For every $X\in\lg $ and $ u\in\lh$ we have:}
\begin{align*}
\allowdisplaybreaks
X \t (R\s(R')) & =  
(X \t R )\s(R')+R ( X \t \s(R')) = (X \t R )\s(R')+R \s( X \t R') \\
& = \beta(\zeta(X))\s(R')+R \s(\beta(\zeta'(X)))=\beta(\zeta(X))\s(R')+R \beta(\s(\zeta'(X)))
\\&=\beta\big (\zeta(X) \s(R')+R \s(\zeta'(X)\big) \, ,\\
u \t (R\s(R')) & =(u \t R) \s(R')+R \s(u \t R')=\zeta(\d(u)) \s(R')+R \s(\zeta'(\d(u))\, .
\end{align*}
As for condition (iii), evaluating the lhs on $X,Y \in\lg$:
\begin{align*}
 R \s(\zeta'([X,Y]))+\zeta([X,Y]) \s(R')&=R\sigma(X \t \zeta'(Y)-Y \t \zeta'(X))+(X \t \zeta(Y)-Y \t \zeta(X))\s(R') \, ,
\end{align*}
while, by using the fact that $\g$ acts by derivations on $\U(\le)^{n \otimes}$ and $X \t R=\beta(\zeta(X))$, the rhs reads:
\begin{align*}
X\t &( R \s(\zeta'(Y))+X \t (\zeta(Y) \s(R'))-Y\t ( R \s(\zeta'(X))+Y \t (\zeta(X) \s(R')) = \\ 
&= \beta(\zeta(X))  \s(\zeta'(Y)) +R\,\, ( X \t \sigma( \zeta'(Y))  +(X \t \zeta(Y))\,\,  \s( R'))+\zeta(Y)\,\, \sigma(\beta(\zeta'(X)) + \\ & \quad -\beta(\zeta(Y))  \s(\zeta'(X))  -R\,\, ( Y \t \sigma( \zeta'(X))  -(Y \t \zeta(X))\,\,  \s( R'))-\zeta(X)\,\, \sigma(\beta(\zeta'(Y))\, .
\end{align*}
Subtracting the two sides yields:
$$\beta(\zeta(X))  \s(\zeta'(Y)) +\zeta(Y)\,\, \sigma(\beta(\zeta'(X))-\beta(\zeta(Y))  \s(\zeta'(X))  -\zeta(X)\,\, \sigma(\beta(\zeta'(Y))\, ,$$
which vanishes in $\Un$ by \eqref{dr}. Finally, note that given a positive integer $n$ and $\tau,\tau' \in {S_n}\times \mathbb{Z}$, {the composition map, below, is bilinear:}
$$\hom_\tau(n,n) \times \hom_{\tau'}(n,n) \to \hom_{\tau\tau'}(n,n) .$$

\subsubsection{2-morphisms and whiskering}\label{2m}

The set of 2-morphisms $\big( n\ra{(R,\zeta,\tau)} n \big) \implies \big( n' \ra{(R', \zeta',\tau')} n' \big) $ is non-empty only if $n=n'$ and $\tau=\tau'$. In this case a 2-morphism $(R,\zeta,\tau) \xRightarrow{T} (R', \zeta',\tau)$ is given by an element $T \in \Un$, such that:
\begin{enumerate}[(i)]
\item $R'=R+\beta(T)$;
\item $\zeta'(X)=\zeta(X)+X \t T \, , \quad \forall \, X\in\lg$.
\end{enumerate}
We denote a 2-morphism as:
 $$\xymatrix{ & n\ar@/^1pc/[rr]^{(R',\zeta',\tau)}\ar@/_1pc/[rr]_{(R,\zeta,\tau)} & \Uparrow T & n}.$$
Note that if $n \ra{(R,\zeta,{\tau})} n$ then also  $n \ra{(R+\beta(T),\zeta',\tau)} n$, for any $T \in \Un$. 
This is because, for $X \in \lg$ and $v \in \lh$:
$$X \t (R +\beta(T))=\beta(\zeta(X))+\beta(X \t T)=\beta(\zeta(X))+\beta(\zeta'(X)-\zeta(X))=\beta(\zeta'(X)), $$
$$ v \t (R +\beta(T))=\zeta(\d(v))+ \d(v) \t T=\zeta(\d(v)) +\zeta'(\d(v)) -\zeta(\d(v))=\zeta'(\d(v)),$$
where we used \eqref{dr}. Also given $X,Y \in \lg$:
\begin{align*}
 \zeta'([X,Y])&=\zeta([X,Y])-[X,Y] \t T=X \t \zeta(Y)-Y \t \zeta(X)-X \t (Y \t T)+Y \t (X \t T)\\
&=X \t \zeta'(Y)-Y \t \zeta'(X).
\end{align*}
The vertical composition of 2-morphisms is:
$$\big( \, (R,\zeta,\tau) \stackrel{T} {\implies} { (R',\zeta',\tau)} \stackrel{T'} {\implies} { (R'',\zeta'',\tau)} \, \big) \, = \, \big( \, (R,\zeta,\tau) \stackrel{T+T'} {\implies} { (R'',\zeta'',\tau)} \, \big) \, .$$
It is trivial to check that  $T+T'$ does connect $(R,\zeta,\tau)$ and $(R'',\zeta'',\tau)$.

The set of 2-morphisms connecting 1-morphisms whose underlying decorated permutation is $\tau$ is denoted by $\hom^2_{\tau}(n,n)$. {It} is in one-to-one correspondence with 
$\hom_{\tau}(n,n) \times \Un $ and it is therefore a vector space. With this identification, the source and target maps $s',t'\colon \hom_{\tau}^2(n,n)\to \hom_{\tau}(n,n)$ have the form:
$$
s'((R,\zeta,\tau)\times T)=(R,\zeta,\tau) \, , \quad \an \quad 
t'((R,\zeta,\tau)\times T)=(R+\beta(T),\zeta',\tau) \, ,$$
where $\zeta'(X)=\zeta(X)-X \t T$. Therefore both $s'$ and $t'$ are linear. For each object $n\in\mathbb{N}$ the vector space of 2-morphisms is the direct sum:
$$ \hom^2(n,n) = \bigoplus_{\tau \in S_n\times\mathbb{Z}} \hom^2_{\tau}(n,n) \, .$$
The {whiskerings} are naturally defined as:
\begin{equation*}
\allowdisplaybreaks
\begin{split}
& \xymatrix{ 
n \ar[rr]^{(R'',\zeta'',\tau'')} & & 
n \ar@/^1pc/[rr]^{(R',\zeta',\tau)} \ar@/_1pc/[rr]_{(R,\zeta,\tau)} & \Uparrow T  & n & = & 
n\ar@/^1pc/[rr]^{(R''\tau''(R'),R''\tau''(\zeta')+\zeta''\tau''( R'),\tau''\tau)} \ar@/_1pc/[rr]_{(R''\tau''(R),R''\tau''(\zeta)+\zeta'' \tau''(R),\tau''\tau)} & \Uparrow R'' \tau''(T) &n 
} \, , \\[.4cm]
& \xymatrix{ 
n \ar@/^1pc/[rr]^{(R',\zeta',\tau)} \ar@/_1pc/[rr]_{(R,\zeta,\tau)} & \Uparrow T & 
n \ar[rr]^{(R'',\zeta'',\tau'')} & & n & = &
n \ar@/^1pc/[rr]^{(R'\tau(R''),R'\tau(\zeta')+\zeta'' \tau(R''),\tau\tau'')} \ar@/_1pc/[rr]_{(R\tau(R''),R\tau(\zeta)+\zeta \tau(R''),\tau\tau')} & \Uparrow  T\tau( R'') &n
} \, .
\end{split}
\end{equation*}
That these  2-morphisms have the correct source and target requires using the relations \eqref{dr} for $\Un$. For example, that $X\t T = \zeta'(X)-\zeta(X)$, for the left whiskering, is proved as:
\begin{align*}
X \t (R''\tau''(T)) & =
(X \t R'')  \tau''(T)+R'' (X \t \tau''(T))=\beta(\zeta''(X))\tau''( T)+R''\tau''(\zeta'(X))-R''\tau''(\zeta(X)) \\ 
& = \zeta''(X) \beta(\tau''(T))+R''\tau''(\zeta'(X))-R''\tau''(\zeta(X))\\
& = \zeta''(X) \tau''(R') -\zeta''(X)\tau''(R)+R''\tau''(\zeta'(X))-R''\tau''(\zeta(X))\\
& = R''\tau''(\zeta'(X)) + \zeta''(X)\tau''(R') -R''\tau''(\zeta(X))-\zeta''(X)\tau''(R) \, .
\end{align*}
The relations \eqref{dr2} are also essential to prove the interchange condition, which essentially follows from:
$$S \,\tau(\beta(T))=\beta(S)\, \tau(T) \, , \qquad \forall \,  S,T \in \Un .$$
Finally, note that the left and right whiskering defines bilinear maps:
$$\hom_\tau(n,n)\times \hom^2_{\tau'} (n,n) \to \hom^2_{\tau\tau'} (n,n) \, , \qquad \hom^2_{\tau'} (n,n) \times  \hom_{\tau} (n,n) \to \hom^2_{\tau'\tau} (n,n) \, .$$

\subsubsection{The monoidal structure}
We write $1_n$ for the identity of $\U(\lg)^{\tn n}=\U(\lg^{\oplus n})$. Given $\s '\in S_{n'}$ we put $n \otimes \s'$ as being the permutation of $\{1,\dots, n+n'\}$ which leaves $\{1,\dots,n\}$ static, whereas $(n \otimes \s')(i)=\s'(i-n)+n$ if $i>n$. We analogously define $\s \tn n'$. We write $n \tn (\s',k')=(n \tn \s',k')$, for a decorated permutation $(\s',k')\in {S_{n'}}\times \mathbb{Z}$. 

\begin{Definition}
The monoidal structure $\tn:\CG\times\CG\to\CG$ of $\CG$ is defined as follows:
\begin{enumerate}
\item given two objects $n,n'$ in $\CG$, {we define} $n \otimes n'=n+n'$;
\item given an object $n'$ and a 1-morphism $n \ra{(R,\zeta,tau)} n$, {we define:} 
$$n' \otimes \left (n \ra{(R,\zeta,\tau)} n \right) \, = \, (n+n')\ra{(1_n \tn R,1_n \tn \zeta,n \tn \tau)} (n+n') \, ;$$ 
similarly, {we put}
$\left (n' \ra{(R,\zeta,\tau)} n'\right)\otimes n = (n'+n)\ra{( R\tn 1_n, \zeta \tn 1_n,\tau\otimes n)} (n'+n)$;
\item given an object $n'$ and a 2-morphism 
$T:(R,\zeta,\tau)\implies (R',\zeta',\tau)$, {we define:}
$$\left( \nsl\mor{n}{n}{(R,\zeta,\tau)}{(R',\zeta',\tau)}{T} \right) \otimes n' = 
\nsl\mor{(n+n')}{(n+n')}{(R\tn 1_{n'},\zeta\tn 1_{n'},\tau\tn n')}{(R'\tn 1_{n'},\zeta'\tn 1_{n'},\tau\tn n')}{T \tn 1_{n'}} \, .$$
\end{enumerate}
\end{Definition}
This definition is clearly distributive with respect to compositions of 1- and 2-morphisms, as well as the whiskering. We comment the interchange laws, see Definition \ref{lin2cat}, condition 6. The condition 6.(i) between 1-morphisms and tensor product holds, being
$$(R,\zeta,\tau)\tn (R',\zeta',\tau')=(R\tn R',R \tn \zeta'+\zeta \tn R',\tau \tn \tau') \, . $$
The condition 6.(ii) between 2-morphisms and 1-morphisms also holds, being:
$$\left (\nsl\mor{n}{n}{(R,\zeta,\tau)}{(R',\zeta',\tau)}{T} \right) \otimes (m \ra{(R'',\zeta'',\tau'')} m) \, = \nsl
\mor{(n+m)}{(n+m)}{(R \tn R'',R \tn \eta''+\eta \tn R'',\tau \tn \tau'')}{(R' \tn R'',R' \tn \eta''+\eta' \tn R'',\tau \tn \tau'')}{T \tn R''} $$
(to prove that the latter 2-morphisms has the correct source and target requires using  \eqref{dr2}). Finally, condition 6.(iii) requires that given 2-morphisms:
$$ \mor{n}{n}{(R,\zeta,\tau)}{(R+\beta(T),\zeta',\tau)}{T} \qquad \an \nsl\qquad
\mor{m}{m}{(P,\eta,\nu)}{(P+\beta(S),\eta',\nu)}{S} $$
it {holds:}
$$ \comp{T\tn (P,\eta,\nu)}{(R+\beta(T),\zeta',\tau)\tn S} = 
\comp{(R,\zeta,\tau)\tn S}{T\tn (P+\beta(S),\eta',\nu)} \, .$$
The previous equality follows from:
$$T \tn R'+(R +\beta(T))\tn T')=R \tn T'+T\tn (R+\beta(T'))\,, $$
(both sides of which can be taken as the definition of $T \tn S$) which in turn follows from relations \eqref{dr2}. 

Given objects $n$ and $m$ there is a natural symmetric braiding $B_{n,m}:n+m\to m+n$, essentially given by the permutation $\s_{n,m}$ that exchanges $n$ and $m$. Namely:
$$\s_{n,m}(i)=\left \{
\begin{CD} 
i+m \quad \textrm{ if } \; i\leq n \\ 
i-n \quad \textrm{ if } \; i > n 
\end{CD} \right . \, ,
\qquad i=1,\ldots ,n,n+1,\ldots,n+m \, .
$$
It is a nice  exercise to check the remaining conditions of Definition \ref{tssml2c}, finishing to prove that $\CG$ is a symmetric strict monoidal linear 2-category.

\subsection{Infinitesimal braidings in a differential crossed module: {symmetric} quasi-invariant tensors}
\label{infbrxmod}
Let $\lG=(\d\colon \lh \to \lg,\t)$ be a differential crossed module. The diagonal map $ \U(\lg) \to \U(\lg)^{\otimes n}$, $x\mapsto x\tn 1\ldots\tn 1 + $ cyclic permutations, will be denoted by $\Delta^n$, and similarly for $\le$ and $\lh$, where $\le= \lh \rtimes_{\t} \lg$. The fact that, if $\s \in S_n$, then $\sigma \circ \Delta^n=\Delta^n$ will be used several times below.
\begin{Definition}[{Symmetric quasi-invariant tensor}]
\label{quasi}
A  {symmetric}  quasi-invariant tensor in $\lG$  is a triple {$(r,c,\xi)$} where:
\begin{enumerate}[(a)]
\item  $r = \sum_q s_q \tn t_q \in \lg \tn \lg$ is a symmetric tensor (thus $r=\sum_q t_q \tn s_q$);
\item $\xi: \lg \to \lg\tn\h \oplus \lh\tn\lg \subset \Ud$ is a linear map, whose image is symmetric. We write it as: 
$$\xi(X)= \sum_a \xi_a(X) \mu_a \tn \nu_a = \sum_a \xi_a (X) \mu'_a \tn \mu''_a+ \sum_a \xi_a (X) \mu''_a \tn \mu'_a \, .$$
where $\xi_a: \lg \to \C$ is a linear map, and in the second expression we use the notation $\mu_a' \in \lg$, $\mu_a'' \in \lh$ to distinguish the part in $\lg\tn\lh$ from the one in $\lh\tn\lg$. 
\item $c$ is an $\lg$-invariant element in $\ker(\partial)\subset \lh$.
\end{enumerate}
The conditions that should be verified are:
\begin{enumerate}[(i)]
\item $X \t r= \beta(\xi(X)) \, , \quad \forall \, X \in \lg$;
\item $u\t r = \xi(\d(u)) \, , \quad \forall \, u \in \lh $; 
\item $\xi([X,Y])=X \t \xi(Y)-Y\t \xi(X) \, , \quad \forall \, X,Y \in \lg$.
\end{enumerate}
\end{Definition}

Given a {symmetric} quasi-invariant tensor in $\lG$, we can construct a totally symmetric infinitesimal 2-braiding in $\CG$. This may be seen as the categorification of the infinitesimal braiding (in the category of $\lg$-modules) associated to a symmetric and $\g$-invariant tensor $r\in\lg\tn\lg$; see the Introduction.

\begin{Theorem}
Let $\lG$ be a differential crossed module. Let $(r,\xi,c)$ be a {symmetric} quasi-invariant tensor in $\lG$. 
The following definition of $(r,T)$ provides the linear totally symmetric monoidal 2-category $\CG$ with a strict totally symmetric infinitesimal 2-braiding. Given objects $n,m$ of $\CG$, the morphism $r_{n,m}: n\tn m \to n \tn m$ is defined as:
\begin{equation}
\label{defr}
r_{n,m}=\left(\sum_{\substack{1\leq i \leq n\\ n+1 \leq j \leq m+n }} r^{ij}, \sum_{\substack{1\leq i \leq n\\ n+1 \leq j \leq m+n }} \xi^{ij},\, \id_{S_{n+m}} , \, 1 \right) \, . 
\end{equation}
Here $1 \in \mathbb{Z}$ and $\id_{S_{n+m}}$ denotes the identity of the group $S_{n+m}$. Given a 1-morphism $(R,\zeta,\s,k)\colon n \to n $, we set:
\begin{equation}
\label{defTl} 
T_{((R,\zeta,\s,k),m)}=-\sum_{ n+1 \leq j \leq m+n } \sum_q \zeta(s_q)\tn t_q^j+kR \tn \Delta^m(c)  =-\sum_q \zeta(s_q)\tn \Delta^m(t_q)+kR \tn \Delta^m(c) \, . 
\end{equation}
Analogously we define:
\begin{equation}
\label{defTr}
T_{(m,(R,\zeta,\s,k))} = - \sum_{1\leq i \leq n}  \sum_q s_q^i \tn \zeta(t_q)+k\Delta^m(c) \tn R = - \sum_q \Delta^n(s_q)  \tn \zeta(t_q) + k\Delta^m(c) \tn R \, .
\end{equation}
\end{Theorem}

\begin{Proof}
We need to prove that $r$ and $T$ are well defined 1- and 2-morphisms in $\CG$, and that the properties of an infinitesimal 2-braiding, see Definition \ref{definf2br}, are satisfied. For the computations it is convenient to write (note $r=\sum_q s_q \tn t_q$):
$$\sum_{\substack{1\leq i \leq n\\ n+1 \leq j \leq m+n }} r^{ij}=\sum_{\substack{1\leq i \leq n\\ n+1 \leq j \leq m+n\\q }} s_q^i \tn t_q^j=\sum_q \Delta^n(s_q) \tn \Delta^m(t_q) \, .$$
We easily see that  conditions (i)-(iii) of Section \ref{o1m} are satisfied, so that $r_{n,m}$ is indeed a 1-morphism in $\CG$. Recalling Definition \ref{definf2br}, we then need to prove that $T_{(R,\zeta,\s,k),m}$ is indeed a 2-morphism $r_{n,m}\,\,(R,\zeta,\s,k) \Rightarrow (R,\zeta,\s,)\,\, r_{n,m}$. To start, we have: 
\begin{equation}
\allowdisplaybreaks
\label{difxi}
\begin{split} 
r_{n,m} \, \left( (R,\zeta,\s,k)\tn m \right) & = 
\left( \sum_{\substack{1\leq i \leq n \\ n+1 \leq j \leq m+n \\ q }} s^i_q R \tn t^j_q \, , \sum_{\substack{1\leq i \leq n\\ n+1 \leq j \leq m+n\\q }} s^i_q \zeta \tn t^j_q + \sum_{\substack{1\leq i \leq n \\ n+1 \leq j \leq m+n \\ l}} \xi_l \, \mu_l^i R  \tn \nu_l^j \, ,\,\,\s\tn m,\, k+1 \right) \, , \\
\left( (R,\zeta,\s,k)\tn m \right) \, r_{n,m} & = 
\left(\sum_{\substack{1\leq i \leq n\\ n+1 \leq j \leq m+n\\q }} Rs^i_q  \tn t^j_q \, ,
\sum_{\substack{1\leq i \leq n\\ n+1 \leq j \leq m+n \\ l}} \xi_l \, R \mu_l^i  \tn 
\nu_l^j + \sum_{\substack{1\leq i \leq n\\ n+1 \leq j \leq m+n\\q }}  \zeta \, s^i_q \tn t^j_q \, ,\,\,\s\tn m,\, k+1 \right) \, .
\end{split}
\end{equation}
We check conditions (i) and (ii) of Section \ref{2m}. The first condition follows directly from condition (i) of Section \ref{o1m}. For the second one, we look at the maps $\g \to \U^{(m+n)}$ involved in \eqref{difxi}. Their difference is, for each $X \in \lg$:
\begin{equation}
\label{dif1}
\begin{split} 
\sum_{\substack{1\leq i \leq n\\ n+1 \leq j \leq m+n \\ l}} \xi_l(X) \, R \mu_l^i  \tn \nu_l^j & +\sum_{\substack{1\leq i \leq n\\ n+1 \leq j \leq m+n\\q }}  \zeta(X) s^i_q \tn t^j_q - \sum_{\substack{1\leq i \leq n\\ n+1 \leq j \leq m+n\\q }} s^i_q \zeta(X)\tn t^j_q -\sum_{\substack{1\leq i \leq n\\ n+1 \leq j \leq m+n \\ l}} \xi_l(X) \, \mu_l^i R  \tn \nu_l^j \\
&= - \sum_{\substack{1\leq i \leq n\\ n+1 \leq j \leq m+n \\ l}} \xi_l(X)  \, \mu_l^i \t R  \tn \nu_l^j -\sum_{\substack{1\leq i \leq n\\ n+1 \leq j \leq m+n\\q }}\big( s^i_q \t \zeta(X)\big)\tn t^j_q \\
&= - \sum_l \xi_l(X) \, \mu_l \t R  \tn \Delta^m(\nu_l)- \sum_q s_q \t \zeta(X) \tn \Delta^m(t_q) \, .
\end{split}
\end{equation}
This has to agree with $X\t T_{(R,\zeta,\s,k),m}$. The latter is computed to be:
\begin{align*}
X \t T_{(R,\zeta,\s,k),m} & = X \t \left( -\sum_q \zeta(s_q)\tn \Delta^m(t_q)+kR \tn \Delta^m(c) \right) \\
& = - \sum_q X\t (\zeta(s_q)) \tn \Delta^m(t_q) - \sum_{\substack{n+1\leq j \leq m+n \\ q}} \zeta(s_q) \tn X\t t_q^j \\
& = - \sum_q s_q \t \zeta(X) \tn \Delta^m(t_q) - \sum_q \zeta( [X,s_q] ) \tn \Delta^m(t_q) - \sum_{\substack{n+1\leq j \leq m+n \\ q}} \zeta(s_q) \tn X\t t_q^j \, .
\end{align*}
Here passing from the first to the second line we used that $X\t\Delta^m(c)=0$ since $c$ is $\lg$-invariant and that $\beta(\zeta(X))\tn\Delta^m(c)=\zeta(X)\tn\Delta^m(\d c)$, by \eqref{dr2}, and this vanishes since $c\in\ker(\d)$. Passing from the second to the third line we used $\zeta([X,Y]) = X\t\zeta(Y)-y\t\zeta(X)$. The last two terms now can be rewritten as $- \, (\zeta\tn\Delta^m)(X\t r)$, which then agrees with $- \, (\zeta\tn\Delta^m)\beta(\xi(X))$. Making this last expression explicit:
\begin{equation*}
\begin{split}
- \, (\zeta\tn\Delta^m)\beta(\xi(X)) 
& = - \sum_a \xi_a(X) \zeta(\mu'_a) \tn \Delta^m(\d(\mu''_a)) + \xi_a(X) \zeta(\d\mu''_a) \tn \Delta^m(\mu'_a) \\ 
& = - \sum_a \xi_a(X) \beta(\zeta(\mu'_a)) \tn \Delta^m(\mu''_a) + \xi_a(X) \mu''_a \t R \tn \Delta^m(\mu'_a) \\
& = - \sum_a \xi_a(X) \mu'_a\t R \tn \Delta^m(\mu''_a) + \xi_a(X) \mu''_a \t R \tn \Delta^m(\mu'_a) \\
& = - \sum_a \xi_a(X) \mu_a\t R \tn \Delta^m(\nu_a)   \, ,
\end{split}
\end{equation*}
so that indeed $X \t T_{(R,\zeta,\s,k),m}$ agrees with \eqref{dif1}.
The naturality condition \eqref{natT} is proven similarly. As for condition \eqref{compT} of the definiton of an infinitesimal 2-braiding, note that:
\begin{align*}
T_{((R,\zeta,\s,k)(R',\zeta',\s',k'), \, m)} 
& = T_{((R \s(R'), \zeta \s(R') + R\s(\zeta'),\s\s',k+k'), \, m)} \\[.2cm]
& = - \sum_q \left( \zeta(s_q) \s(R') +  R \s(\zeta'(s_q)) \right) \tn \Delta^m(t_q) + (k+k') R\s(R') \tn \Delta^m(c),
\end{align*}
which clearly is the right hand side of  \eqref{compT}. Note that this is the point where the decorations of the permutations become essential. 

The interchangeability condition \eqref{Tfg} reads as follow: given $(R,\xi,\s,k):n\to n$ and  $(R',\xi',\s',k'):m\to m$ 1-morphisms, then it must hold:
$$
\comp{ T_{( (R,\zeta,\s,n) \, , \, m )} \, \left( n\tn (R',\zeta',\s',k') \right) }
{ \left( (R,\zeta,\s,k)\tn m \right) \, T_{ ( n \, , \, (R',\zeta',\s',k) )}} \, = \,
\comp{ T_{ ( n \, , \, (R',\zeta',\s',k') )} \, \left( (R,\zeta,\s,k)\tn m \right) }
{ \left( n\tn (R',\zeta',\s',k') \right) \, T_{( (R,\zeta,\s,k) \, , \, m )}} \, .
$$  
This means:
\begin{multline*}
(R\tn 1_m)(\s\tn m)(T_{ ( n \, , (R',\zeta',\s',k') )}) +
T_{ ( (R,\zeta,\s,k) \, , \, m ) } \, (\s\tn m)(1_n\tn R) + \\
- (1_n\tn R')(n\tn\s')(T_{ ( (R,\zeta,\s,k) \, , \, m )} - 
T_{ ( n \, , (R',\zeta',\s',k') )} \, (n\tn \s')(R\tn 1_m) = 0 \, .
\end{multline*}
Since $c$ is $\lg$-invariant the lhs of the previous expression reduces to
\begin{equation*}
\begin{split}
\sum_q  \left( 
R\Delta^m(s_q)\tn\zeta'(t_q) + \zeta(s_q)\tn\Delta^m(t_q) R'  \right. & \left. - \, \zeta(s_q) \tn R' \Delta^m(t_q) - \Delta^n(s_q)R \tn \zeta'(t_q) \right) \, = \\
& = \sum_q ( - s_q\t R \tn \zeta'(t_q) + \zeta(s_q)\tn t_q \t R )\\
& = \sum_q ( - \beta(\zeta(s_q)) \tn \zeta'(t_q) + \zeta(s_q) \tn \beta(\zeta'(t_q)) ) = 0\,.
\end{split}
\end{equation*}
The remaining properties of Definition \ref{definf2br} follow from simple calculations.
\qed\end{Proof}
\begin{Definition}[Coherent {symmetric} quasi-invariant tensor]
\label{cohP}
 Given a {symmetric} quasi-invariant tensor $(r,\xi, c)$, let 
$$P=\sum_i s_i \tn \xi(t_i)+ c \tn r.$$
Then $(r,\zeta, c)$, is said to be coherent if $\; P_{123}+P_{231}+P_{312}=0 \, . $
\end{Definition}

We omit the following proof, since it only involves lengthy but straightforward computations.

\begin{Theorem}
Let  $(r,\zeta, c)$ be a coherent {symmetric} quasi-invariant tensor in a differential crossed module $\lG$. Then the associated infinitesimal 2-braiding in the linear monoidal 2-category $\CG$ is coherent.
\end{Theorem}
\subsection{A {symmetric} quasi-invariant tensor in the {String} Lie-2-algebra}\label{qitsl2a}

{We exhibit an explicit example of a {symmetric} quasi-invariant tensor, in the differential crossed module $\st$ (defined by Wagemann \cite{wag06}) representing the String Lie 2-algebra, in the $\mathfrak{sl}(2,\C)$ case. This can be seen as a re-interpretation of the infinitesimal {2-Yang-Baxter operator} in the String Lie 2-algebra that we constructed in \cite[Section 4]{CFM2}, now in the more {powerful} framework of {symmetric} quasi-invariant tensors and infinitesimal 2-braidings.} 
\subsubsection{Differential crossed modules and cohomology}
{Let  $\lG=(\d\colon \lh \to \lg, \t)$ be a differential crossed module. We have an induced representation of ${\rm coker}(\d)$ in $\ker(\d)$. To $\lG$ we can associate a Lie algebra cohomology class $\w \in H^3\big({\rm coker}(\d),\ker(\d)\big)$, normally called the $k$-invariant of $\lG$, and  denoted $\w=k(\lG)$.  In such a case, we say that $\lG$ is a realization, or a model, of $\w$.  Two realizations of a given cohomology class are not necessarily isomorphic, being however weakly equivalent, or, what is the same, equivalent (in a broader sense) in the larger category of weak Lie-2-algebras; \cite{BC}.}

{The construction of $\w$ from $\lG=(\d\colon \lh \to \lg, \t)$ is quite simple.
Choose a section $s\colon {\rm coker}(\d)\to \lg$ of the projection map $p\colon \lg \to {\rm coker}(\d)$. Also choose a section $s'\colon \d(\lh) \to \lh$ of $\d\colon \lh \to \lg$.  Given $Y,Z \in {\rm coker}(\d)$, clearly $s([Y,Z])-[s(Y),s(Z)] \in \ker(p)=\d(\lh).$ Put $b(X,Y)=s'\left (s([Y,Z])-[s(Y),s(Z)] \right)$. {We then have:}}
\begin{equation}\label{kinv}
 \w(X,Y,Z) \doteq b([X,Y],Z)+\textrm{cyclic permutations}-s(X) \t b(Y,Z)+\textrm{cyclic permutations}\in \ker(\d) \subset \lh.
 \end{equation}

{Recall that the Lie algebra cohomology of any simple Lie algebra $\lg$, {considering the trivial representation in $\C$}, is, in dimension three, generated by the  Lie algebra 3-cocycle $(X,Y,Z) \in {\lg\times \lg\times \lg} \mapsto \frac{1}{2}\langle X,[Y,Z]\rangle$, where $\langle,\rangle$ denotes the Cartan-Killing form.
A (generic) string differential crossed module will mean, in this paper, any differential crossed module whose $k$-invariant is $\w$, for a given simple Lie algebra.}
{Two papers where explicit constructions of (generic) string differential crossed modules appear are \cite{BSCS} and \cite{wag06}; the latter only in the $\mathfrak{sl}(2,\C)$ case. These construction do not appear to be isomorphic.}

{In this paper we work with Wagemann's \cite{wag06} explicit realization of the string differential crossed module, henceforth denoted $\st$ and called the String Differential Crossed module. This {model} is completely algebraic, therefore fitting naturally into the framework of this paper. We however note that the {model developed in}  \cite{BSCS} {could also be} made fully algebraic, and, in theory, the construction presented here of an infinitesimal 2-braiding derived from $\st$ could be carried out by using this alternative realization of $\w$, which {has the advantage of generalizing immediately} to all other simple Lie algebras.} 

In the very recently made available \cite{RW} another explicit realization of the string Lie algebra  is constructed for a large class of Lie algebras (containing for instance all simple Lie algebras). Each of these realizations is there proven to likewise possess symmetric quasi-invariant tensors.

\begin{Remark}{
 The nomenclature ``string Lie 2-algebra'' has been common since the arising of the paper \cite{BSCS}. Specifically, given a simple Lie group $G$, with Lie algebra $\lg$, it is proven that {(at least for the particular model of the string Lie 2-algebra appearing there)} there exists a Lie 2-group $STRING(G)$ having the string Lie 2-algebra as the tangent space at the identity. The Lie-2-group $STRING(G)$ is closely related (after taking a simplicial nerve) with the String Group referred to in \cite{MS}, which originally appeared in \cite{Sto}.  The original string group arose as being a 3-connected cover of $G$, for $G={\rm Spin}(n)$. An excellent explanation of all of this appears in \cite{NSW}.}        
\end{Remark}

\subsubsection{The {String} differential crossed module}
\label{sectionstring}
{From now on we abbreviate $\sl=\sl(\C)$.
{We review the construction of the differential crossed module $\st$. The original construction is due to Wagemann \cite{wag06}.}
Let $\lg$ be a Lie algebra and $V$ a $\lg$-module. An anti-symmetric map $w: \wedge^2(\lg) \to V$ is called a 2-cocycle if its coboundary $\delta^V(w):\wedge^3(\lg)\to V$, defined as: }
\begin{equation}
\label{cocy}
\delta^V(w)(X,Y,Z) =X \t w(Y,Z) +Y\t w(Z,X)+Z \t w(X,Y)+w(X,[Y,Z])+w(Y,[Z,X])+w(Z,[X,Y])
\end{equation}
vanishes for all $X,Y,Z \in \lg$.
Let $W_1$ be the Lie algebra of polynomial vector fields in one variable $x$, with Lie bracket given by the commutator of vector fields (here $f'$ denotes the derivative of $f$): 
$$
\big[ f(x) \, \frac{d}{d x} \, , \, g(x) \, \frac{d}{d x} \big] =  \big(f \, \frac{d g}{d x} - \frac{d f}{d x} \, g \big)(x) \, \frac{d}{d x} = \big(fg' - f'g\big)(x) \, \frac{d}{d x} \, , \qquad \forall \, f(x)\,\frac{d}{d x} \, , \, g(x)\,\frac{d}{d x} \in W_1 \, . $$
Let $\Fo$ be the space of {complex} polynomials in {a formal} variable $x$, and $\Fu$ the space of formal 1-forms  $f(x)d x$, with $f(x)$ {a complex}  polynomial. These are $W_1$-modules via the Lie derivative:
\begin{equation}
\label{w1act}
f(x)\frac{d}{d x} \, \t \, g(x) = (fg')(x) \, , \quad f(x)\frac{d}{d x} \, \t \, g(x) d x = (fg' + f'g)(x) \, d x \, , \quad \, f(x) \, \frac{d}{d x} \in W_1 \, , g(x) \in \Fo \, , g(x)\, d x \in \Fu \, .
\end{equation}
The formal de Rham differential $d\colon \Fo \to \Fu$ is a map of $W_1$-modules. We have a short exact sequence of $W_1$-modules, namely $\{0\} \to \C \to \Fo \ra{d} \Fu \to\{0\}$, where $\C$ is given the trivial action. This induces a long exact sequence in Lie algebra cohomology. Therefore we have a  connecting homomorphism $S\colon H^2(W_1,\Fu) \to H^3(W_1,\C)$.

{Simple calculations prove that we have a Lie algebra 2-cocycle ${\alpha}$ with the form:} $$\left (f(x)\frac{d}{d x} , g(x) \frac{d}{d x} \right)\in W_1 \times W_1\stackrel{\overline{\alpha}}{\longmapsto} \frac{1}{2}\begin{vmatrix} f'(x) & g'(x)\\ f''(x) & g''(x)\end{vmatrix} \, dx \in \Fu. $$ Put $\overline{\w}=S(\overline{\alpha}) \in H^3(W_1,\C)$. {Noting that} given three polynomials {$f,g$ and $h$} we have:
$$\frac{d }{ d x} \begin{vmatrix}
                   f(x) & g(x) & h(x) \\
                   f'(x) & g'(x) & h'(x) \\
                   f''(x) & g''(x) & h''(x)
                  \end{vmatrix} = 
                  \begin{vmatrix}
                 \begin{vmatrix} f(x) & g(x)\\ f'(x) & g'(x) \end{vmatrix}' &  \begin{vmatrix} f(x) & g(x)\\ f'(x) & g'(x) \end{vmatrix}'' \\ h'(x)  & h''(x)
                  \end{vmatrix}+\textrm{ cyclic permutations of } (f,g,h),
 $$
{the well known general form of the connecting homomorphism $S\colon H^2(W_1,\Fu) \to H^3(W_1,\C)$ gives that $\overline{\w}$} is the 3-cocycle $(f,g,h) \in W_1 \times W_1 \times W_1\mapsto 
 \frac{1}{2} \begin{vmatrix}    f(0) & g(0) & h(0) \\
                    f'(0) & g'(0)  & h'(0) \\
                    f''(0) & g''(0) & h''(0)
                   \end{vmatrix}\in \C$,
normally called Godbillon-Vey 3-cocycle; \cite{GW}. 

We have an exact sequence of Lie algebras:
$$\{0\} \to \C \to \Fo \ra{\d} \Fu \rtimes_{\overline{\alpha}} W_1 \to W_1 \to \{0\}.$$
Here $\C$, $\Fo$ and $\Fu$ are given the trivial (abelian) Lie algebra structure. In addition $\Fu \rtimes_{\overline{\alpha}} W_1$ denotes semidirect product, twisted by $\overline{\alpha}$. Explicitly: $\big[(a,y),(b,z)\big] := \big( y\t b - z\t a + \overline{\alpha}(y,z) \, , [y,z] \, \big) \, , \quad \forall \, a,b \in \Fu \, , \, y,z \in W_1$. The boundary map $d\colon \Fo \ra{\d} \Fu \ltimes_{\overline{\alpha}} W_1$ explicitly reads $\d(f)=(df,0)$. Together with the action of $ \Fu \rtimes_{\overline{\alpha}} W_1$ on $\Fo$, where $(a,y)\t f=y \t f$, we have a differential crossed module. 
It can be easily seen  {from \eqref{kinv}, and it is proved more conceptually in \cite{wag06},} that this crossed module realizes the cohomology  class $\overline{\w}=S(\overline{\alpha})$.

{Let us put $f=\begin{pmatrix} 0 & -1\\ 0 &0 \end{pmatrix}$, \quad $e=\begin{pmatrix} 0 & 0\\ 1 &0 \end{pmatrix}$  \quad and  $k=\begin{pmatrix} -1/2 & 0\\ 0 &1/2 \end{pmatrix}$, so that the $\sl$ commutation relations read:} 
\begin{align} &[f,e]=2k;&&[k,e]=e; &[k,f]=-f. \end{align}
{We have a group homomorphism $m\colon {\rm SL}(2,\mathbb{R}) \to \diff(\mathbb{R})$, the diffeomorphism group of $\mathbb{R}$, obtained by considering M\"{o}bius transformations, namely $m\left(\begin{pmatrix} a & b\\ c &d \end{pmatrix}\right)(x)=\displaystyle{\frac{ax+b}{cx+d}}$.
The (complexified) infinitesimal form of this group morphisms yields a Lie algebra map $m\colon \sl \to W_1$, where $W_1$ is the  Lie algebra of polynomial vector fields, in the usual way: $m(Y)(x)=\frac{d}{d t} m\big(\exp (- t Y)\big)(x)_{t=0}$. Therefore we can identify $\sl$ as being the Lie subalgebra of $W_1$  generated by:}
\begin{equation}\label{incl}
f = \frac{d}{d x} \, , \qquad\qquad k = x \, \frac{d}{d x} \, , \qquad\qquad e = x^2 \, \frac{d}{d x}.
\end{equation}

The 2-cocyle $\overline{\alpha}\colon W_1 \times W_1 \to \Fu$ restricts to a Lie algebra 2-cocycle $\alpha:\wedge^2(\sl)\to\Fu$, which explicitly is:
\begin{equation}\label{dalpha}
\alpha(k,e)= -\alpha(e,k) = dx   \, , \qquad \textrm{ and zero otherwise.}
\end{equation}
 The differential crossed module 
$$\st=(\dd:\Fo\lra\fsl,\t)$$ 
representing the {String Lie 2-algebra} fits inside the  four-terms exact sequence:
\begin{equation}
\label{mainseq} 
0 \lra \mathbb{C} \stackrel{i}{\lra} \Fo \stackrel{\dd}{\lra} \fsl \stackrel{\pi}{\lra} \sl \lra 0 \, . 
\end{equation}  
where $\fsl$ is the vector space $\Fu\oplus\sl$,  endowed with the Lie bracket:
\begin{equation}
\label{fslbr}
\big[(a,y),(b,z)\big] := \big( y\t b - z\t a + \alpha(y,z) \, , [y,z] \, \big) \, , \quad \forall \, a,b \in \Fu \, , \, y,z \in \sl \, .
\end{equation}
And $\Fo$  is a   $\fsl$-module, with:
\begin{equation}
\label{fslmod}
(a,y)\t f := y\t f \, , \qquad (a,y)\in\fsl , f \in F. 
\end{equation}
(All actions of $\sl$ are restrictions of actions of $W_1$.)
By construction, $\st$ realizes the restriction of the Godbillon-Wey cocyle to $\sl$. Explicit computations prove {that this restriction is $\w(X,Y,Z)=\frac{1}{2}\langle X,[Y,Z]\rangle$.} (The Cartan-killing form $\langle \, , \, \rangle:\sl\ot\sl\to\mathbb{C}$ is taken with the normalization: $\langle x,y \rangle =  \mathrm{Tr}\,(\mathrm{ad}_x\circ \mathrm{ad}_y)$.)

We repeat the relevant structures in $\st=(\dd:\Fo\lra\fsl,\t)$: the Lie bracket in $\fsl$ is given in \eqref{fslbr}, while in $\Fo$ is trivial. The $(\fsl)$-action on $\Fo$  is of the type \eqref{fslmod}, trivially extending the one of $\sl$. The maps $\dd$ and $\pi$ are explicitly $\dd(f)=(df,0)$ and $\pi(\omega,y)= y$, for $f\in\Fo \, , \, (\omega,y)\in\fsl$; hence $\mathrm{ker} (\dd) = \mathbb{C}$ and $\mathrm{coker} (\dd)=\sl$.
Therefore, {since $d\colon \Fo \to \Fu$ is the exterior derivative,}  the kernel of $\partial\colon \Fo\lra\fsl$ is generated by the constant function $1_{\F_0}$ in $\C$. 
{From \eqref{w1act},  we see that $\mathrm{ker} (\dd)$ is $\sl$-invariant, thus also $(\fsl)$-invariant. This fact will play a prime role later.}

\subsubsection{The construction of the {symmetric} quasi-invariant tensor}

Consider the {String} differential crossed module $\st=(\d\colon \Fo\to \fsl,\t)$. Given $X \in \sl$, we put $\overline{X}=(0,X) \in \fsl$. We will also write $\ov{h}=(h,0)$ for $h\in\Fu$. Let $Q\colon \Fu \to \Fo$ be a linear primitive of the exterior differentiation: {$d(Q(w))=w$, for each 1-form $w \in \Fu$.} {We will choose it so that  $Q(dx)=x$.} We use $Q$ to introduce the linear map: 
$$ \w:\wedge^2(\sl)\to\Fo \quad\textrm{ such that }\quad  \w(X,Y)=Q(\alpha(X,Y)) \, .$$
Recall \eqref{dalpha}. As a consequence of $\eqref{dr2}$, we can see that:
\begin{equation*}
\begin{split}
1_\Fo \tn [\ov{X},\ov{Y}] 
& = 1_\Fo \tn (\alpha(X,Y),[X,Y]) = 
1_\Fo \tn \ov{\alpha(X,Y)} + 1_\Fo \tn \ov{[X,Y]} \\
& = 1_\Fo \tn \ov{d(\w(X,Y))} + 1_\Fo \tn \ov{[X,Y]} =
1_\Fo \tn \d(\w(X,Y)) + 1_\Fo \tn \ov{[X,Y]} \\
& = \d(1_\Fo) \tn \w(X,Y) + 1_\Fo \tn \ov{[X,Y]} =1_\Fo \tn \overline{[X,Y]}\, .
\end{split}
\end{equation*}
We therefore have the following absolutely crucial relation in $\Ud$, valid for all $X,Y \in \sl$:
\begin{equation}
\label{key}
1_\Fo \tn [\ov{X},\ov{Y}] = 1_\Fo \tn \overline{[X,Y]}.
\end{equation}

\noindent
Let $r=\sum_{i} s_i \tn t_i \in \sl \tn \sl$ be the infinitesimal braiding of $\sl$ associated to the Cartan-Killing form. Explicitly:
$$  r = \, f\tn e + \, e\tn f - 2 \, k \tn k \, =\sum_{i} s_i \tn t_i \, .
$$
The following identity (the $\sl$-invariance of $r$) will be used time after time below:
\begin{equation}
\label{invr}
\sum_i [s_i,X] \tn t_i+s_i \tn [t_i,X]=0 \, , \qquad 
\textrm{ for all } X \in \sl \, . 
\end{equation}
Let $\Phi:=\delta^{\Fo}(\w):\wedge^3{\sl}\to\Fo$ be the coboundary of $\w$ (an explicit expression is as in \eqref{cocy}).
\begin{Lemma}Given $X,Y \in \sl$ we have:
\begin{equation}
\label{idufo}
\sum_i \Phi(s_i,X,Y) \tn t_i=\ufo\tn [X,Y] \, , \qquad
\sum_i s_i\tn\Phi(t_i,X,Y) = [X,Y]\tn\ufo
\end{equation}
\end{Lemma}
\begin{Proof} Clearly the map $$(X,Y) \mapsto \sum_{i} \Phi(s_i,X,Y) \tn t_i $$
is bilinear and  antisymmetric, so we only need to check \eqref{idufo}  for the pairs $(f,e)$, $(f,k)$ and $(k,e)$. This follows from an easy calculation.
\qed\end{Proof}

We are now ready to construct the {symmetric} quasi-invariant tensor in $\st$. Recalling Definition \ref{quasi}, we consider the following objects:
\begin{enumerate}[(a)]
\item $\br=\sum_i \bs_i \tn \bt_i \in (\fsl)\tn(\fsl)$ (the obvious lifting of $r$ to $\fsl$);
\item the map $\xi:(\fsl) \to \big ((\fsl) \tn \Fo\big) \oplus \big (\Fo\tn (\fsl) \big)\subset\Ud$  defined as $\xi = - \, \xi_0 + C$, where
\begin{align}
\label{defxi1}
\xi_0(\ov{X}) & = \sum_i \w(s_i,X) \tn \bt_i +\bs_i \tn \w(t_i,X) \, , \\
\xi_0(\ov{h}) & = \sum_i s_i \t Q(h) \tn \bt_i + \bs_i \tn t_i \t Q(h) \, , \\
C((h,X)) & = 1_{\Fo} \tn \bX + \bX\tn 1_{\Fo} \, ;
\end{align} 
\item $c \in \ker(\d)=\C \subset \Fo \, .$  
\end{enumerate}
We note that since $\ker(\d)$ is $\sl$-invariant, $\xi_0(\bar{h})$ does not depend on the chosen primitive $Q\colon \Fu \to \Fo$ of the exterior derivation.
\begin{Theorem}
The triple  $(\br,\xi,c)$ is a {symmetric} quasi-invariant tensor in the {String} differential crossed module $\st$.
\end{Theorem}
\begin{Proof}We need to check conditions (i)-(iii) of Definition \ref{quasi}. As for condition (i): 
\begin{align*}
\bX \t \br&=\sum_i \left( (\alpha(X,s_i),[X,s_i]) \tn \bt_i +\bs_i \tn (\alpha(X,t_i),[X,t_i])\right) \\
&=\sum_i \left( \, \ov{\alpha(X,s_i)} \tn \bt_i + \bs_i \tn \ov{\alpha(X,t_i)} \, \right) = - \, \beta(\xi_0(\bX)) = \beta(\xi(\bX))\, ,\\[0.2cm] 
\ov{h} \t \br & = \sum_i - \left( \, \ov{s_i \t h} \, \tn \bt_i + \bs_i \tn \, \ov{t_i \t h} \, \right) \\
& =  \sum_i - \beta \left( s_i \t Q(h) \tn \bt_i + \bs_i \tn t_i \t Q(h) \right) = 
- \, \beta(\xi_0(\ov{h})) = \beta(\xi(\ov{h}))\,.
\end{align*}
Now condition (ii); if $a \in \Fo$ (remember that $Q:\Fu\to\Fo$ is a primitive of $d$):
$$ a \t \br  -\xi(\d(a)) =  \sum_i \left( s_i \t (Qd(a) - a) \tn \bt_i + \bs_i \tn t_i \t (Qd(a)-a) \right)\,, $$
and this vanishes since  $Qd(a)-a$ is in $\ker(d)=\C$, which is $\sl$-invariant. 

Finally, condition (iii); we start by computing:
\begin{align*}
\xi_0([\bX , \bY]) & = \xi_0(\ov{[X,Y]}) - \xi_0(\ov{\alpha(X,Y)}) \\
& = \sum_i \left( \w(s_i,[X,Y])\tn\bt_i + \bs_i\tn\w(t_i,[X,Y]) + s_i\t(Q\alpha(X,Y))\tn\bt_i +
\bs_i\tn t_i\t(Q\alpha(X,Y)) \right)\,, \\
\bX \t \xi_0(\bY) & = \sum_i \left( X\t\w(s_i,Y)\tn\bt_i + \w(s_i,Y)\tn [\bX,\bt_i] + 
[\bX,\bs_i]\tn\w(t_i,Y) - \bs_i\tn X\t\w(t_i,Y) \right) \, ,\\
\bY\t\xi_0(\bX) & =  \sum_i \left( Y\t\w(s_i,X)\tn\bt_i + \w(s_i,X)\tn [\bY,\bt_i] +
[\bY,\bs_i]\tn\w(t_i,X) + \bs_i\tn Y\t\w(t_i,X) \right)\, .
\end{align*}
We now claim that:
\begin{equation}
\label{claim}
\begin{split}
\xi_0([\bX , \bY]) - \bX \t \xi_0(\bY) + \bY\t\xi_0(\bX) & = \sum_i \left( \Phi(s_i,X,Y)\tn\bt_i + \bs_i\tn\Phi(t_i,X,Y) \right) \\
& = \ov{[X,Y]}\tn\ufo + \ufo\tn \ov{[X,Y]}\, ,
\end{split}
\end{equation}
where the left hand side has been rewritten using (the lift of) \eqref{idufo}. To prove the claim one uses the explicit expression of $\Phi$, see \eqref{cocy}, and checks that:
\begin{equation*}
\begin{split}
\xi_0([\bX , \bY]) - \bX \t \xi_0(\bY) + & \bY\t\xi_0(\bX) -  \sum_i \left( \Phi(s_i,X,Y)\tn\bt_i + \bs_i\tn\Phi(t_i,X,Y) \right) = \\
& \sum_i \left( - \, \w(s_i,Y)\tn [\bX,\bt_i] - [\bX,\bs_i] \tn\w(t_i,Y) + \w(s_i,X)\tn [\bY,\bt_i] 
+ [\bY,\bs_i]\tn\w(t_i,X) \, + \right. \\
& \qquad\qquad \left. - \, \w(X,[Y,s_i])\tn\bt_i - \w(Y,[s_i,X])\tn\bt_i - \bs_i\tn\w(X,[Y,t_i]) -
\bs_i\tn\w(Y,[t_i,X]) \right) \, .
\end{split}
\end{equation*}
We know expand each term of the form $[\bX,\bt_i]$ into $\ov{\alpha(X,t_i)}+\ov{[X,t_i]}$. The terms involving the cocycle $\alpha$ are:
\begin{equation*}
\sum_i \left( - \, \w(s_i,Y)\tn \ov{\alpha(X,t_i)} - \ov{\alpha(X,s_i)}\tn\w(s_i,X) + 
\w(s_i,X)\tn \ov{\alpha(Y,t_i)} + \ov{\alpha(Y,s_i)}\tn\w(t_i,X) \right) \, ,
\end{equation*}
and they cancel pairwise by a suitable use of \eqref{dr2} and from the definition of $\w$; for example:
$$ \w(s_i,Y)\tn\ov{\alpha(X,t_i)} = \w(s_i,Y)\tn\d Q\alpha(X,t_i) = \d\w(s_i,Y)\tn\w(X,t_i)=\ov{\alpha(s_i,Y)}\tn\w(X,t_i).$$
The remaining eight terms not containing $\alpha$ cancel pairwise thanks to the $\sl$-invariance of $r$; for example:
$$ - \, \w(s_i,Y)\tn \ov{[X,t_i]} - \w(Y,[s_i,X])\tn\bt_i = 0 \, ,$$
and similarly for the other pairs. This proves the claim. We go on with:
\begin{equation*}
\begin{split}
C([\bX,\bY]) - \bX\t C(\bY) + \bY\t C(\bX) 
& = C(\alpha(X,Y),\ov{[X,Y]}) - \bX\t C(\bY) + \bY\t C(\bX) \\
& = \ufo\tn\ov{[X,Y]} + \ov{[X,Y]}\tn\ufo\, , 
\end{split}
\end{equation*}
where we used \eqref{key}, so that summing up all the contributions we have eventually proved condition (iii) for $\xi([\bX,\bY])$. Now the same property for $\xi([\bX,\ov{h}])$. We start by computing:
\begin{align*}
\xi_0([\bX,\ov{h}]) & = \xi_0(\ov{X\t h}) = s_i \t Q(\ov{X\t h}) \tn \bt_i + \bs_i\tn t_i\t Q(\ov{X\t h})\,,  \\
\bX\t\xi_0(\ov{h}) & = \bX\t (s_i\t Q(h)) \tn\bt_i + s_i\t Q(h)\tn [\bX,\bt_i] + [\bX,\bs_i]\tn t_i\t Q(h) + \bs_i\tn X\t(t_i\t Q(h))\,, \\
\ov{h}\t\xi_0(\bX) & = \, - \, \w(s_i,X)\tn\ov{t_i\t h} - \ov{s_i\t h}\tn\w(t_i,X) \, .
\end{align*}
{Now note that for each $X,Y \in \sl$ then $X\t(Y\t Q(h)) = X\t(Q(Y\t h))$. This follows from:}
$$ {\d\big(Y\t Q(h) - Q(Y\t h)\big) = Y\t (\d Q(h)) - Y\t h =  Y\t h- Y\t h=0\, ,} $$
together with the fact that $\ker(\d)$ is $\sl$-invariant. All this is used to write:
\begin{equation*}
\begin{split}
\xi_0([\bX,\ov{h}]) & - \bX\t\xi_0(\ov{h}) + \ov{h}\t\xi_0(\bX) = \\
& = \, - \, [X,s_i]\t Q(h)\tn\bt_i - \bs_i\tn [X,t_i]\t Q(h) - s_i\t Q(h)\tn [\bX,\bt_i] - [\bX,\bs_i]\tn t_i\t Q(h)  \, + \\
& \quad \,  - \, \w(s_i,X)\tn \ov{t_i\t h} - \ov{s_i\t h}\tn \w(t_i,X) \\
& = \, - \, \left( \, [X,s_i]\t Q(h)\tn\bt_i + s_i\t Q(h)\tn \ov{[X,t_i]} \, \right) - \left( \,
\bs_i\tn [X,t_i]\t Q(h) + \ov{[X,s_i]} \tn t_i\t Q(h) \, \right) + \\
& \quad \, - s_i\t Q(h)\tn\ov{\alpha(X,t_i)} - \ov{\alpha(X,s_i)}\tn t_i\t Q(h) - \w(s_i,X)\tn \ov{t_i\t h} - \ov{s_i\t h}\tn \w(t_i,X) \, .
\end{split}
\end{equation*}
In the two parenthesis we read the $\sl$-invariance of $r$, so they vanish. The four terms in the last line cancel pairwise, again by a suitable use of \eqref{dr2} and from the definition of $\w$; for example:
$$ - s_i\t Q(h)\tn\ov{\alpha(X,t_i)} = s_i\t Q(h) \tn \d Q\alpha(X,t_i) = \d(s_i\t Q(h)) \tn \w(X,t_i) = \ov{s_i\t h} \tn \w(X,t_i) \, .$$
This completes the proof.
\qed\end{Proof}

We now address the coherence of the {symmetric} quasi-invariant tensor $(\br,\xi,c)$. Recalling Definition \ref{cohP} this amounts to $P_{123} + P_{231} + P_{312} = 0$, where now $P$ reads:
\begin{equation*}
\begin{split} 
P & = \sum_i \left( \bs_i\tn\xi(\bt_i) \right) + c\tn \br = 
\sum_i \left( - \bs_i\tn\xi_0(\bt_i) + \bs_i \tn C(\bt_i) \right) + c\tn \br \\
 & = \sum_{i,j} \left( \, - \, \bs_i\tn\w(s_j,t_i)\tn\bt_j - \bs_i\tn\bs_j\tn\w(t_j,t_i) \, \right) 
 + \sum_i \bs_i \tn C(\bt_i) + c\tn\br \, .
\end{split}
\end{equation*}
\begin{Theorem}
The {symmetric} quasi-invariant tensor $(\br,\xi,c)$ on $\st$ is coherent {if, and only if,} $c= \, -2 \cdot \ufo$.
\end{Theorem}
\begin{Proof}By direct computation. We explicit: 
\begin{equation*}
\begin{split}
P_{123} & = \, -2 \, \ov{k}\tn x \tn\ov{f} + 2 \, \ov{f}\tn x\tn\ov{k} + 2\, \ov{f}\tn\ov{k}\tn x - 
 2 \, \ov{k}\tn\ov{f}\tn x + 
 \sum_i \left( \bs_i\tn\ufo\tn\bt_i + \bs_i\tn\bt_i\tn\ufo + c\tn\bs_i\tn\bt_i \right) \\
P_{231} & = \, -2 \, \ov{f}\tn\ov{k}\tn x + 2\, \ov{k}\tn\ov{f}\tn x + 2\, x\tn\ov{f}\tn\ov{k} 
-2 \, x\tn\ov{k}\tn\ov{f} + \sum_i \left( \bt_i\tn\bs_i\tn\ufo + \ufo\tn\bs_i\tn\bt_i + 
\bt_i\tn c\tn \bs_i \right) \\
P_{312} & = \, -2\, x\tn\ov{f}\tn\ov{k} + 2\, x\tn\ov{k}\tn\ov{f} + 2\, \ov{k}\tn x\tn\ov{f} 
- 2\, \ov{f}\tn x\tn \ov{k} + 
 \sum_i \left( \ufo\tn\bt_i\tn\bs_i + \bt_i\tn\ufo\tn\bs_i + \bs_i\tn\bt_i\tn c \right)
\end{split}
\end{equation*}
and note that in $P_{123} + P_{231} + P_{312}$ the terms without summation on $i$ cancel directly, while the sum over $i$ vanishes {if, and only if,} $c= \, -2 \cdot \ufo$. \qed\end{Proof}
\section*{Acknowledgements}
L.S. Cirio was partially supported by the German Research Foundation DFG through the grant SFB 878
``Groups, Geometry \& Actions''. J. Faria Martins was partially supported by CMA/FCT/UNL, under the
project UID/MAT/00297/2013.
This work was partially supported by FCT (Portugal) through the grant  {``Geometry and Mathematical Physics'' (FCT EXCL/MAT-GEO/0222/2012).}


\end{document}